# W. B. VASANTHA KANDASAMY

# SMARANDACHE LOOPS

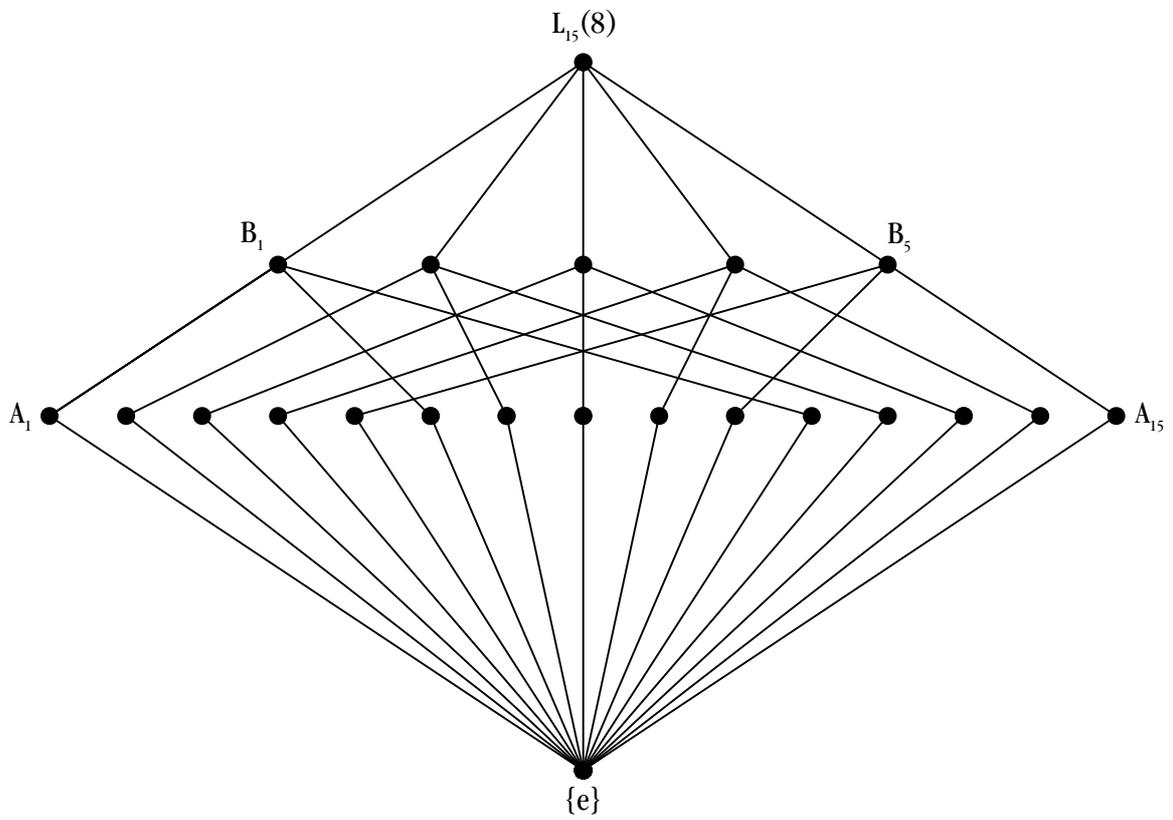

**2002**

# Smarandache Loops

**W. B. Vasantha Kandasamy**

Department of Mathematics
Indian Institute of Technology, Madras
Chennai – 600036, India
e-mail: *vasantha@iitm.ac.in*
web: *http://mat.iitm.ac.in/~wbv*

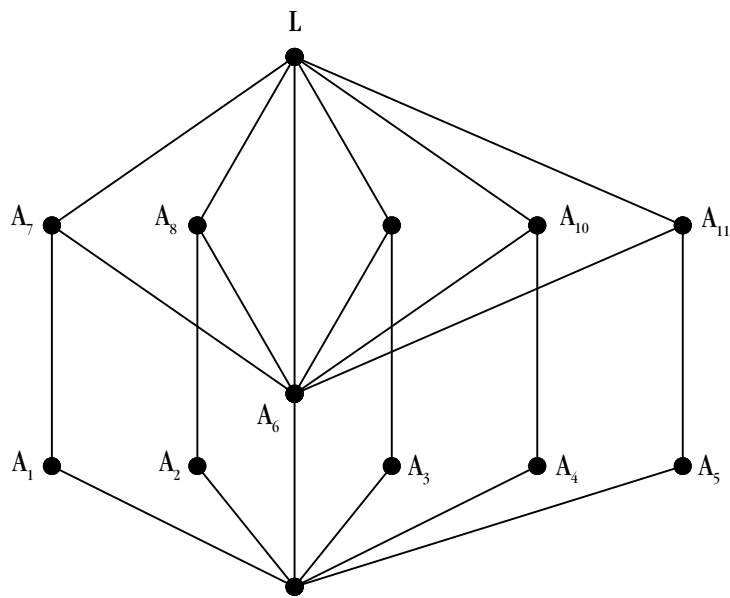

2002

The picture on the cover represents the lattice of subgroups of the Smarandache loop $L_{15}(8)$. The lattice of subgroups of the commutative loop $L_{15}(8)$ is a non-modular lattice with 22 elements. This is a Smarandache loop which satisfies the Smarandache Lagrange criteria. But for the Smarandache concepts one wouldn't have studied the collection of subgroups of a loop.

# CONTENTS

# PREFACE

The theory of loops (groups without associativity), though researched by several mathematicians has not found a sound expression, for books, be it research level or otherwise, solely dealing with the properties of loops are absent. This is in marked contrast with group theory where books are abundantly available for all levels: as graduate texts and as advanced research books.

The only three books available where the theory of loops are dealt with are: R. H. Bruck, *A Survey of Binary Systems*, Springer Verlag, 1958 recent edition (1971); H. O. Pflugfelder, *Quasigroups and Loops: Introduction*, Heldermann Verlag, 1990; Orin Chein, H. O. Pflugfelder and J.D.H. Smith (editors), *Quasigroups and Loops: Theory and Applications*, Heldermann Verlag, 1990. But none of them are completely devoted for the study of loops.

The author of this book has been working in loops for the past 12 years, and has guided a Ph.D. and 3 post-graduate research projects in this field of loops feels that the main reason for the absence of books on loops is the fact that it is more abstract than groups. Further one is not in a position to give a class of loops which are as concrete as the groups of the form $S_n$, $D_{2n}$ etc. which makes the study of these non-associative structures much more complicated. To overcome this problem in 1994 the author with her Ph.D. student S. V. Singh has introduced a new class of loops using modulo integers. They serve as a concrete examples of loops of even order and it finds an application to colouring of the edges of the graph $K_{2n}$.

Several researchers like Bruck R. H., Chibuka V. O., Doro S., Giordano G., Glauberman G., Kunen K., Liebeck M.W., Mark P., Michael Kinyon, Orin Chein, Paige L.J., Pflugfelder H.O., Phillips J.D., Robinson D. A., Solarin A. R. T., Tim Hsu, Wright C.R.B. and by others who have worked on Moufang loops and other loops like Bol loops, A-loops, Steiner loops and Bruck loops. But some of these loops become Moufang loops. Orin Chein, Michael Kinyon and others have studied loops and the Lagrange property.

The purpose of this book entirely lies in the study, introduction and examination of the Smarandache loops. As a result, this book doesn't give a full-fledged analysis on loops and their properties. However, for the sake of readers who are involved in the study of loop theory we have provided a wide-ranging list of papers in the reference. We expect the reader to have a good background in algebra and more specifically a strong foundation in loops and number theory.

The study of Smarandache loops was initiated by the author in the year 2002. This book introduces over 75 Smarandache concepts on loops, and most of these concepts are illustrated by examples. In fact several of the Smarandache loops have classes of loops which satisfy the Smarandache notions.

This book is structured into five chapters. Chapter one which is introductory in nature covers some notions about groups, graphs and lattices. Chapter two gives some basic properties of loops. The importance of this chapter lies in the introduction of a new class of loops of even order. We prove that the number of different representations of right alternative loop of even order (2n), in which square of each element is identity is equal to the number of distinct proper (2n – 1) edge colourings of the complete graph $K_{2n}$.

In chapter three we introduce Smarandache loops and their Smarandache notions. Except for the Smarandache notions several of the properties like Lagrange's criteria, Sylow's criteria may not have been possible. Chapter four introduces Smarandache mixed direct product of loops which enables us to define a Smarandache loops of level II and this class of loops given by Smarandache mixed direct product gives more concrete and non-abstract structure of Smarandache loops in general and loops in particular. The final section gives 52 research problems for the researchers in order to make them involved in the study of Smarandache loops. The list of problems provided at the end of each section is a main feature of this book.

I deeply acknowledge the encouragement that Dr. Minh Perez extended to me during the course of this book. It was because of him that I got started in this endeavor of writing research books on Smarandache algebraic notions.

I dedicate this book to my parents, Worathur Balasubramanian and Krishnaveni for their love.

**Chapter one**

# GENERAL FUNDAMENTALS

In this chapter we shall recall some of the basic concepts used in this book to make it self-contained. As the reader is expected to have a good knowledge in algebra we have not done complete justice in recollecting all notions. This chapter has three sections. In the first section we just give the basic concepts or notion like equivalence relation greatest common divisor etc. Second section is devoted to giving the definition of group and just stating some of the classical theorems in groups like Lagrange's, Cauchy's etc. and some basic ideas about conjugates. Further in this section one example of a complete graph is described as we obtain an application of loops to the edge colouring problem of the graph $K_{2n}$. In third section we have just given the definition of lattices and its properties to see the form of the collection of subgroups in loops, subloops in loops and normal subloops in loops. The subgroups in case of Smarandache loops in general do not form a modular lattice.

Almost all the proofs of the theorem are given as exercise to the reader so that the reader by solving them would become familiar with these concepts.

## 1.1 Basic Concepts

The main aim of this section is to introduce the basic concepts of equivalence relation, equivalence class and introduce some number theoretic results used in this book. Wherever possible the definition when very abstract are illustrated by examples.

**DEFINITION 1.1.1**: *If a and b are integers both not zero, then an integer d is called the greatest common divisor of a and b if*

  i.   *d > 0*
  ii.  *d is a common divisor of a and b and*
  iii. *if any integer f is a common divisor of both a and b then f is also a divisor of d.*

*The greatest common divisor of a and b is denoted by g.c.d (a, b) or simply (g.c.d). If a and b are relatively prime then (a, b) = 1.*

**DEFINITION 1.1.2**: *The least common multiple of two positive integers a and b is defined to be the smallest positive integer that is divisible by a and b and it is denoted by l.c.m or [a, b].*

**DEFINITION 1.1.3**: *Any function whose domain is some subset of set of integers is called an arithmetic function.*

**DEFINITION 1.1.4**: *An arithmetic function f(n) is said to be a multiplicative function if f(mn) = f(m) f(n) whenever (m, n) = 1.*

**Notation**: If x ∈ R (R the set of reals), Then [x] denotes the largest integer that does not exceed x.

***Result 1***: If d = (a, c) then the congruence ax ≡ b (mod c) has no solution if d ∤ b and it has d mutually incongruent solutions if d/b.

***Result 2***: ax ≡ b (mod c) has a unique solution if (a, c) = 1.

## 1.2 A few properties of groups and graphs

In this section we just recall the definition of groups and its properties, and state the famous classical theorems of Lagrange and Sylow. The proofs are left as exercises for the reader.

**DEFINITION 1.2.1**: *A non-empty set of elements G is said to form a group if on G is defined a binary operation, called the product and denoted by ' • ' such that*

1. *a, b ∈ G implies a • b ∈ G (closure property).*
2. *a, b, c ∈ G implies a • (b • c) = (a • b) • c (associative law).*
3. *There exist an element e ∈ G such that a • e = e • a = a for all a ∈ G (the existence of identity element in G).*
4. *For every a ∈ G there exists an element $a^{-1}$ ∈ G such that a • $a^{-1}$ = $a^{-1}$ • a = e (the existence of inverse element in G)*

**DEFINITION 1.2.2**: *A group G is said to be abelian (or commutative) if for every a, b ∈ G, a • b = b • a.*

A group which is not commutative is called non-commutative. The number of elements in a group G is called the order of G denoted by o(G) or |G|. The number is most interesting when it is finite, in this case we say G is a finite group.

**DEFINITION 1.2.3**: *Let G be a group. If a, b ∈ G, then b is said to be a conjugate of a in G if there is an element c ∈ G such that b = $c^{-1}$ac. We denote a conjugate to b by a ~ b and we shall refer to this relation as conjugacy relation on G.*

**DEFINITION 1.2.4**: *Let G be a group. For $a \in G$ define $N(a) = \{x \in G \mid ax = xa\}$. $N(a)$ is called the normalizer of a in G.*

**THEOREM** (Cauchy's Theorem For Groups): *Suppose G is a finite group and $p/o(G)$, where p is a prime number, then there is an element $a \neq e \in G$ such that $a^p = e$, where e is the identity element of the group G.*

**DEFINITION 1.2.5**: *Let G be a finite group. Then*

$$o(G) = \sum \frac{o(G)}{o(N(a))}$$

*where this sum runs over one element a in each conjugate class, is called the class equation of the group G.*

**DEFINITION 1.2.6**: *Let $X = (a_1, a_2, \ldots, a_n)$. The set of all one to one mappings of the set X to itself under the composition of mappings is a group, called the group of permutations or the symmetric group of degree n. It is denoted by $S_n$ and $S_n$ has n! elements in it.*

*A permutation $\sigma$ of the set X is a cycle of length n if there exists $a_1, a_2, \ldots, a_n \in X$ such that $a_1 \sigma = a_2, a_2 \sigma = a_3, \ldots, a_{n-1} \sigma = a_n$ and $a_n \sigma = a_1$ that is in short*

$$\begin{pmatrix} a_1 & a_2 & \ldots & a_{n-1} & a_n \\ a_2 & a_3 & \ldots & a_n & a_1 \end{pmatrix}.$$

*A cycle of length 2 is a transposition. Cycles are disjoint, if there is no element in common.*

***Result***: Every permutation $\sigma$ of a finite set X is a product of disjoint cycles.

The representation of a permutation as a product of disjoint cycles, none of which is the identity permutation, is unique up to the order of the cycles.

**DEFINITION 1.2.7**: *A permutation with $k_1$ cycles of length 1, $k_2$ cycles of length 2 and so on, $k_n$ cycles of length n is said to be a cycle class $(k_1, k_2, \ldots, k_n)$.*

**THEOREM** (Lagrange's): *If G is a finite group and H is a subgroup of G, then $o(H)$ is a divisor of $o(G)$.*

The proof of this theorem is left to the reader as an exercise.

It is important to point out that the converse to Lagrange's theorem is false-a group G need not have a subgroup of order m if m is a divisor of o(G). For instance, a finite group of order 12 exists which has no subgroup of order 6. Consider the symmetric group $S_4$ of degree 4, which has the alternating subgroup $A_4$, of order 12. It is easily established 6/12 but $A_4$ has no subgroup of order 6; it has only subgroups of order 2, 3 and 4.

We also recall Sylow's theorems which are a sort of partial converse to Lagrange's theorem.

**THEOREM** (First Part of Sylow's theorem): *If p is a prime number and $p^\alpha/o(G)$, where G is a finite group, then G has a subgroup of order $p^\alpha$.*

The proof is left as an exercise for the reader.

**DEFINITION 1.2.8**: *Let G be a finite group. A subgroup of G of order $p^m$, where $p^m/o(G)$ but $p^{m+1} \nmid o(G)$ is called a p-Sylow subgroup of G.*

**THEOREM** (Second part of Sylow's Theorem): *If G is a finite group, p a prime and $p^\alpha/o(G)$ but $p^{\alpha+1} \nmid o(G)$, then any two subgroups of G of order $p^\alpha$ are conjugate.*

The assertion of this theorem is also left for the reader to verify.

**THEOREM** (Third part of Sylow's theorem): *The number of p-sylow subgroups in G, for a given prime is of the from 1 + kp.*

Left for the reader to prove. For more about these proofs or definitions kindly refer I.N.Herstein [27] or S.Lang [34].

**DEFINITION 1.2.9**: A s*imple graph in which each pair of disjoint vertices are joined by an edge is called a complete graph.*

*Upto isomorphism, there is just one complete graph on n vertices.*

**Example 1.2.1**: Complete graph on 6 vertices that is $K_6$ is given by the following figure:

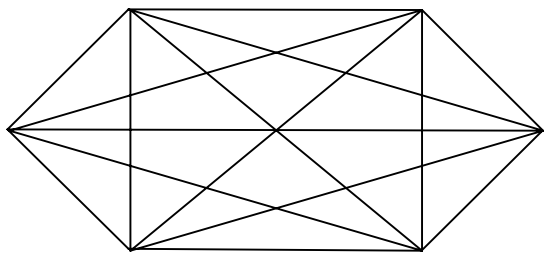

**DEFINITION 1.2.10**: *A edge colouring η of a loopless graph is an G assignment of k colour 1, 2, …, k to the edges of the graph G. The colouring η is proper if no two adjacent edges have the same colour in the graph G.*

For more about colouring problems and graphs please refer [3].

## 1.3 Lattices and its properties

The study of lattices has become significant, as we know the normal subgroups of a group forms a modular lattice. A natural question would be what is the structure of the set of normal subloops of a loop? A still more significant question is what is the structure of the collection of Smarandache subloops of a loop? A deeper question is what is the structure of the collection of all Smarandache normal subloops? It is still a varied study to find the lattice structure of the set of all subgroups of a loop as in the basic definition of Smarandache loops we insist that all loops should contain subgroups to be a Smarandache loop. In view of this we just recall the definition of lattice, modular lattice and the distributive lattice and list the basic properties of them.

**DEFINITION 1.3.1**: *Let A and B be non-empty sets. A relation R from A to B is a subset of $A \times B$. Relations from A to B are called relations on A, for short. If (a, b) $\in$ R then we write aRb and say that 'a is in relation R to b'. Also, if a is not in relation R to b we write $a\not{R}b$.*

*A relation R on a non-empty set A may have some of the following properties:*

    *R is reflexive if for all a in A we have aRa.*
    *R is symmetric if for all a and b in A, aRb implies bRa.*
    *R is antisymmetric if for all a and b in A, aRb and bRa imply a = b.*
    *R is transitive if for a, b, c in A, aRb and bRc imply aRc.*

*A relation R on A is an equivalence relation if R is reflexive, symmetric and transitive. In this case [a] = {b $\in$ A / aRb} is called the equivalence class of a for any a $\in$ A.*

**DEFINITION 1.3.2**: *A relation R on a set A is called a partial order (relation) if R is reflexive, anti-symmetric and transitive. In this case (A, R) is called a partial ordered set or a poset.*

**DEFINITION 1.3.3**: *A partial order relation $\leq$ on A is called a total order (or linear order) if for each a, b $\in$ A either a $\leq$ b or b $\leq$ a. (A, $\leq$) is then called a chain, or totally ordered set.*

**DEFINITION 1.3.4**: *Let $(A, \leq)$ be a poset and $B \subseteq A$.*

i)     *$a \in A$ is called a upper bound of $B \Leftrightarrow \forall b \in B, b \leq a$.*
ii)    *$a \in A$ is called a lower bound of $B \Leftrightarrow \forall b \in B, a \leq b$.*
iii)   *The greatest amongst the lower bounds whenever it exists is called the infimum of B, and is denoted by inf B.*
iv)   *The least upper bound of B, whenever it exists is called the supremum of B and is denoted by sup B.*

**DEFINITION 1.3.5**: *A poset $(L, \leq)$ is called lattice ordered if for every pair x,y of elements of L the $\sup(x, y)$ and $\inf(x, y)$ exist.*

**DEFINITION 1.3.6**: *An (algebraic) lattice $(L, \cap, \cup)$ is a non empty set L with two binary operations $\cap$ (meet) and $\cup$ (join) (also called intersection or product and union or sum respectively) which satisfy the following conditions for all $x, y, z \in L$.*

$(L_1)$ $x \cap y = y \cap x$,                $x \cup y = y \cup x$
$(L_2)$ $x \cap (y \cap z) = (x \cap y) \cap z$,    $x \cup (y \cup z) = (x \cup y) \cup z$.
$(L_3)$ $x \cap (x \cup y) = x$,           $x \cup (x \cap y) = x$.

*Two applications of $(L_3)$ namely $x \cap x = x \cap (x \cup (x \cap x)) = x$ lead to the additional condition $(L_4)$ $x \cap x = x$, $x \cup x = x$.*

         *$(L_1)$ is the commutative law*
         *$(L_2)$ is the associative law*
         *$(L_3)$ is the absorption law and*
         *$(L_4)$ is the idempotent law.*

**DEFINITION 1.3.7**: *A partial order relation $\leq$ on A is called a total order (or linear order) if for each $a, b \in A$ either $a \leq b$ or $b \leq a$. $(A, \leq)$ is then called a chain or totally ordered set. A totally ordered set is a lattice ordered set and $(A, \leq)$ will be defined as a chain lattice.*

**DEFINITION 1.3.8**: *Let L and M be lattices. A mapping $f: L \to M$ is called a*

i)     *Join homomorphism if $x \cup y = z \Rightarrow f(x) \cup f(y) = f(z)$.*
ii)    *Meet homomorphism if $x \cap y = z \Rightarrow f(x) \cap f(y) = f(z)$.*
iii)   *Order homomorphism if $x \leq y \Rightarrow f(x) \leq f(y)$.*

*f is a homomorphism (or lattice homomorphism) if it is both a join and a meet homomorphism. Injective, surjective or bijective lattice homomorphism are called lattice monomorphism, epimorphism, isomorphism respectively.*

**DEFINITION 1.3.9**: *A non-empty subset S of a lattice L is called a sublattice of L, if S is a lattice with respect to the restriction of $\cap$ and $\cup$ on L onto S.*

*Example 1.3.1*: Every singleton of a lattice L is a sublattice of L

**DEFINITION 1.3.10**: *A lattice L is called modular if for all x, y, z ∈ L.*

(M) $x \leq z \Rightarrow x \cup (y \cap z) = (x \cup y) \cap z$         *(modular equation)*

**_Result 1:_** The lattice is non-modular if even a triple x, y, z ∈ L does not satisfy modular equation. This lattice (given below) will be termed as a pentagon lattice.

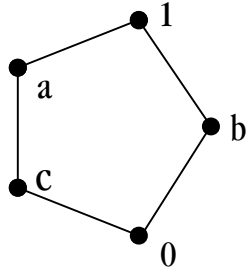

**_Result 2:_** A lattice L is modular if none of its sublattices is isomorphic to a pentagon lattice.

**DEFINITION 1.3.11**: *A lattice L is called distributive if either of the following conditions hold for all x, y, z in L*

$x \cup (y \cap z) = (x \cup y) \cap (x \cup z)$
$x \cap (y \cup z) = (x \cap y) \cup (x \cap z)$         *(distributive equations)*

**_Result 3:_** A modular lattice is distributive if and only if none of its sublattices is isomorphic to the diamond lattice given by the following diagram

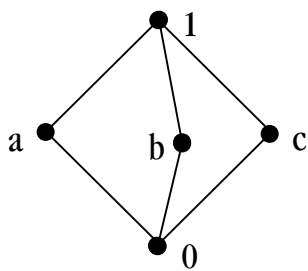

***Result 4:*** A lattice is distributive if and only if none of its sublattices is isomorphic to the pentagon or diamond lattice.

For more about lattices the reader is requested to refer Birkhoff[2] and Gratzer [26].

***Result 5:*** Let G be any group. The set of all normal subgroups of G forms a modular lattice.

***Result 6:*** Let G be a group. The subgroup of G in general does not form a modular lattice.

It can be easily verified that the 10 subgroups of the alternating group $A_4$ does not form a modular lattice.

***Result 7:*** Let $S_n$ be the symmetric group of degree n, n ≥ 5. The normal subgroups of $S_n$ forms a 3 element chain lattice.

**Chapter two**

# LOOPS AND ITS PROPERTIES

This chapter is completely devoted to the introduction of loops and the properties enjoyed by them. It has 7 sections. The first section gives the definition of loop and explains them with examples. The substructures in loops like subloops, normal subloops, associator, commutator etc are dealt in section 2. Special elements are introduced and their properties are recalled in section 3. In section 4 we define some special types of loops. The representation and isotopes of loops is introduced in section 5. Section 6 is completely devoted to the introduction of a new class of loops of even order using integers and deriving their properties. Final section deals with the applications of these loops to the edge-colouring of the graph $K_{2n}$.

## 2.1 Definition of loop and examples

We at this juncture like to express that books solely on loops are meagre or absent as, R.H.Bruck deals with loops on his book "*A Survey of Binary Systems*", that too published as early as 1958, [6]. Other two books are on "*Quasigroups and Loops*" one by H.O. Pflugfelder, 1990 [50] which is introductory and the other book co-edited by Orin Chein, H.O. Pflugfelder and J.D. Smith in 1990 [16]. So we felt it important to recall almost all the properties and definitions related with loops. As this book is on Smarandache loops so, while studying the Smarandache analogous of the properties the reader may not be running short of notions and concepts about loops.

**DEFINITION 2.1.1**: *A non-empty set L is said to form a loop, if on L is defined a binary operation called the product denoted by '•' such that*

  a. *For all a, b ∈ L we have a • b ∈ L (closure property).*
  b. *There exists an element e ∈ L such that a • e = e • a = a for all a ∈ L (e is called the identity element of L).*
  c. *For every ordered pair (a, b) ∈ L × L there exists a unique pair (x, y) in L such that ax = b and ya = b.*

Throughout this book we take L to be a finite loop, unless otherwise we state it explicitly, L is an infinite loop. The binary operation '•' in general need not be associative on L. We also mention all groups are loops but in general every loop is not a group. Thus loops are the more generalized concept of groups.

*Example 2.1.1*: Let (L, ∗) be a loop of order six given by the following table. This loop is a commutative loop but it is not associative.

| * | e | $a_1$ | $a_2$ | $a_3$ | $a_4$ | $a_5$ |
|---|---|---|---|---|---|---|
| e | e | $a_1$ | $a_2$ | $a_3$ | $a_4$ | $a_5$ |
| $a_1$ | $a_1$ | e | $a_4$ | $a_2$ | $a_5$ | $a_3$ |
| $a_2$ | $a_2$ | $a_4$ | e | $a_5$ | $a_3$ | $a_1$ |
| $a_3$ | $a_3$ | $a_2$ | $a_5$ | e | $a_1$ | $a_4$ |
| $a_4$ | $a_4$ | $a_5$ | $a_3$ | $a_1$ | e | $a_2$ |
| $a_5$ | $a_5$ | $a_3$ | $a_1$ | $a_4$ | $a_2$ | e |

Clearly (L, *) is non-associative as $(a_4 * a_3) * a_2 = a_4$ and $a_4 * (a_3 * a_2) = a_4 * a_5 = a_2$. Thus $(a_4 * a_3) * a_2 \neq a_4 * (a_3 * a_2)$.

***Example 2.1.2***: Let L = {e, a, b, c, d} be a loop with the following composition table. This loop is non-commutative.

| • | e | a | b | c | d |
|---|---|---|---|---|---|
| e | e | a | b | c | d |
| a | a | e | c | d | b |
| b | b | d | a | e | c |
| c | c | b | d | a | e |
| d | d | c | e | b | a |

L is non-commutative as a • b ≠ b • a for a • b = c and b • a = d. It is left for the reader as an exercise to verify L is non-associative.

**Notation**: Let L be a loop. The number of elements in L denoted by o(L) or |L| is the order of the loop L.

***Example 2.1.3***: Let L be a loop of order 8 given by the following table. L = {e, $g_1$, $g_2$, $g_3$, $g_4$, $g_5$, $g_6$, $g_7$} under the operation ' • '.

| • | e | $g_1$ | $g_2$ | $g_3$ | $g_4$ | $g_5$ | $g_6$ | $g_7$ |
|---|---|---|---|---|---|---|---|---|
| e | e | $g_1$ | $g_2$ | $g_3$ | $g_4$ | $g_5$ | $g_6$ | $g_7$ |
| $g_1$ | $g_1$ | e | $g_5$ | $g_2$ | $g_6$ | $g_3$ | $g_7$ | $g_4$ |
| $g_2$ | $g_2$ | $g_5$ | e | $g_6$ | $g_3$ | $g_7$ | $g_4$ | $g_1$ |
| $g_3$ | $g_3$ | $g_2$ | $g_6$ | e | $g_7$ | $g_4$ | $g_1$ | $g_5$ |
| $g_4$ | $g_4$ | $g_6$ | $g_3$ | $g_7$ | e | $g_1$ | $g_5$ | $g_2$ |
| $g_5$ | $g_5$ | $g_3$ | $g_7$ | $g_4$ | $g_1$ | e | $g_2$ | $g_6$ |
| $g_6$ | $g_6$ | $g_7$ | $g_4$ | $g_1$ | $g_5$ | $g_2$ | e | $g_3$ |
| $g_7$ | $g_7$ | $g_4$ | $g_1$ | $g_5$ | $g_2$ | $g_6$ | $g_3$ | e |

Now we define commutative loop.

**DEFINITION 2.1.2**: *A loop (L, •) is said to be a commutative loop if for all a, b ∈ L we have a • b = b • a.*

The loops given in examples 2.1.1 and 2.1.3 are commutative loops. If in a loop (L, •) we have at least a pair a, b ∈ L such that a • b ≠ b • a then we say (L, •) is a non-commutative loop.

The loop given in example 2.1.2 is non-commutative.

*Example 2.1.4*: Now consider the following loop (L, •) given by the table:

| •     | e     | $g_1$ | $g_2$ | $g_3$ | $g_4$ | $g_5$ |
|-------|-------|-------|-------|-------|-------|-------|
| e     | e     | $g_1$ | $g_2$ | $g_3$ | $g_4$ | $g_5$ |
| $g_1$ | $g_1$ | e     | $g_3$ | $g_5$ | $g_2$ | $g_4$ |
| $g_2$ | $g_2$ | $g_5$ | e     | $g_4$ | $g_1$ | $g_3$ |
| $g_3$ | $g_3$ | $g_4$ | $g_1$ | e     | $g_5$ | $g_2$ |
| $g_4$ | $g_4$ | $g_3$ | $g_5$ | $g_2$ | e     | $g_1$ |
| $g_5$ | $g_5$ | $g_2$ | $g_4$ | $g_1$ | $g_3$ | e     |

We see a special quality of this loop viz. in this loop xy ≠ yx for any x, y ∈ L \ {e} with x ≠ y. This is left as an exercise for the reader to verify.

**PROBLEMS**

1. Does there exist a loop of order 4?
2. Give an example of a commutative loop of order 5.
3. How many loops of order 5 exist?
4. Can we as in the case of groups say all loops of order 5 are commutative?
5. How many loops of order 4 exist?
6. Can we have a loop (which is not a group) to be generated by a single element?
7. Is it possible to have a loop of order 3? Justify your answer.
8. Give an example of a non-commutative loop of order 7.
9. Give an example of a loop L of order 5 in which xy ≠ yx for any x, y ∈ L\{1}, x ≠ y
10. Find an example of a commutative loop of order 11.

## 2.2 Substructures in loops

In this section we introduce the concepts of substructures like, subloop, normal subloop, commutator subloop, associator subloop, Moufang centre and Nuclei of a loop. We recall the definition of these concepts and illustrate them with examples. We have not proved or recalled any results about them but we have proposed some problems at the end of this section for the reader to solve. A new notion called strictly non-commutative loop studied by us in 1994 is also introduced in this section. Finally the definition of disassociative and power associative loops by R.H.Bruck are given.

**DEFINITION 2.2.1**: *Let L be a loop. A non-empty subset H of L is called a subloop of L if H itself is a loop under the operation of L.*

***Example 2.2.1***: Consider the loop L given in example 2.1.3 we see $H_i = \{e, g_i\}$ for i = 1, 2, 3, … , 7 are subloops of L.

**DEFINITION 2.2.2**: *Let L be a loop. A subloop H of L is said to be a normal subloop of L, if*

1. *xH = Hx.*
2. *(Hx)y = H(xy).*
3. *y(xH) = (yx)H*

*for all x, y ∈ L.*

**DEFINITION 2.2.3**: *A loop L is said to be a simple loop if it does not contain any non-trivial normal subloop.*

***Example 2.2.2***: The loops given in example 2.1.1 and 2.1.3 are simple loops for it is left for the reader to check that these loops do not contain normal subloops, in fact both of them contain subloops which are not normal.

**DEFINITION 2.2.4**: *The commutator subloop of a loop L denoted by L' is the subloop generated by all of its commutators, that is, ⟨{x ∈ L / x = (y, z) for some y, z ∈ L}⟩ where for A ⊆ L, ⟨A⟩ denotes the subloop generated by A.*

**DEFINITION 2.2.5**: *If x, y and z are elements of a loop L an associator (x, y, z) is defined by, (xy)z = (x(yz)) (x, y, z).*

**DEFINITION 2.2.6**: *The associator subloop of a loop L (denoted by A(L)) is the subloop generated by all of its associators, that is ⟨{x ∈ L / x = (a, b, c) for some a, b, c ∈ L}⟩.*

**DEFINITION 2.2.7**: *A loop L is said to be semi alternative if (x, y, z) = (y, z, x) for all x, y, z ∈ L, where (x, y, z) denotes the associator of elements x, y, z ∈ L.*

**DEFINITION 2.2.8**: *Let L be a loop. The left nucleus $N_\lambda = \{a \in L / e = (a, x, y)$ for all x, y ∈ L} is a subloop of L. The middle nucleus $N_\mu = \{a \in L / e = (x, a, y)$* *for all $x, y \in L$} is a subloop of L. The right nucleus $N_\rho$ = {$a \in L$ / $e = (x, y, a)$ for all $x, y \in L$} is a subloop of L.*

*The nucleus $N(L)$ of the loop L is the subloop given by $N(L) = N_\lambda \cap N_\mu \cap N_\rho$.*

**DEFINITION 2.2.9**: *Let L be a loop, the Moufang center C(L) is the set of all elements of the loop L which commute with every element of L, that is, $C(L) = \{x \in L / xy = yx$ for all $y \in L\}$.*

**DEFINITION 2.2.10**: *The centre Z(L) of a loop L is the intersection of the nucleus and the Moufang centre, that is $Z(L) = C(L) \cap N(L)$.*

It has been observed by Pflugfelder [50] that N(L) is a subgroup of L and that Z(L) is an abelian subgroup of N(L). This has been cross citied by Tim Hsu [63] further he defines Normal subloops of a loop L is a different way.

**DEFINITION [63]**: *A normal subloop of a loop L is any subloop of L which is the kernel of some homomorphism from L to a loop.*

Further Pflugfelder [50] has proved the central subgroup Z(L) of a loop L is normal in L.

**DEFINITION [63]**: *Let L be a loop. The centrally derived subloop (or normal commutator- associator subloop) of L is defined to be the smallest normal subloop $L' \subset L$ such that $L / L'$ is an abelian group. Similarly nuclearly derived subloop (or normal associator subloop) of L is defined to be the smallest normal subloop $L_1 \subset L$ such that $L / L_1$ is a group. Bruck proves L' and $L_1$ are well defined.*

**DEFINITION [63]**: *The Frattini subloop $\phi(L)$ of a loop L is defined to be the set of all non-generators of L, that is the set of all $x \in L$ such that for any subset S of L, $L = \langle x, S \rangle$ implies $L = \langle S \rangle$. Bruck has proved as stated by Tim Hsu $\phi(L) \subset L$ and $L / \phi(L)$ is isomorphic to a subgroup of the direct product of groups of prime order.*

It was observed by Pflugfelder [50] that the Moufang centre C(L) is a loop.

**DEFINITION [63]**: *A p-loop L is said to be small Frattini if $\phi(L)$ has order dividing p. A small Frattini loop L is said to be central small Frattini if $\phi(L) \leq Z(L)$. The interesting result proved by Tim Hsu is.*

**THEOREM [63]**: *Every small Frattini Moufang loop is central small Frattini.*

The proof is left to the reader.

**DEFINITION [40]**: *Let L be a loop. The commutant of L is the set (L) = {a ∈ L / ax = xa ∀ x ∈ L}. The centre of L is the set of all a ∈ C(L) such that a • xy = ax • y = x • ay = xa • y and xy • a = x • ya for all x, y ∈ L. The centre is a normal subloop. The commutant is also known as Moufang Centre in literature.*

**DEFINITION [39]**: *A left loop (B, •) is a set B together with a binary operation '•' such that (i) for each a ∈ B, the mapping x → a • x is a bijection and (ii) there exists a two sided identity 1 ∈ B satisfying 1 • x = x • 1 = x for every x ∈ B. A right loop is defined similarly. A loop is both a right loop and a left loop.*

**DEFINITION [11]**: *A loop L is said to have the weak Lagrange property if, for each subloop K of L, |K| divides |L|. It has the strong Lagrange property if every subloop K of L has the weak Lagrange property.*

A loop may have the weak Lagrange property(For more about these notions refer Orin Chein et al [11]).

**DEFINITION [41]**: *Let L be a loop. The flexible law FLEX: $x • yx = xy • x$ for all $x, y \in L$. If a loop L satisfies left alternative laws that is $y • yx = yy • x$ then LALT. L satisfies right alternative laws RALT: $x • yy = xy • y$.*

There is also the inverse property IP (Michael K. Kinyon et al [41] have proved in IP loops RALT and LALT are equivalent).

*In a loop L, the left and right translations by $x \in L$ are defined by $yL(x) = xy$ and $yR(x) = yx$ respectively. The multiplication group of L is the permutation group on L, $Mlt(L) = \langle R(x), L(x): x \in L \rangle$ generated by all left and right translations. The inner mapping group is the subgroup $Mlt_1(L)$ fixing 1. If L is a group, then $Mlt_1(L)$ is the group of inner automorphism of L. In an IP loop, the AAIP implies that we can conjugate by J to get $L(x)^J = R(x^{-1})$, $R(x)^J = L(x^{-1})$ where $\theta^J = J^{-1} \theta J = J \theta J$ for a permutation $\theta$. If $\theta$ is an inner mapping so is $\theta^J$.*

**DEFINITION [41]**: *ARIF loop is an IP loop L with the property $\theta^J = \theta$ for all $\theta \in Mlt_1(L)$. Equivalently, inner mappings preserve inverses that is $(x^{-1}) \theta = (x\theta)^{-1}$ for all $\theta \in Mlt_1(L)$ and for all $x \in L$.*

**DEFINITION 2.2.11**: *A map $\phi$ from a loop L to another loop $L_1$ is called a loop homomorphism if $\phi(ab) = \phi(a) \phi(b)$ for all $a, b \in L$.*

**DEFINITION [65]**: *Let L be a loop L is said to be a strictly non-commutative loop if xy ≠ yx for any x, y ∈ L (x ≠ y, x ≠ e, y ≠ e where e is the identity element of L).*

**DEFINITION 2.2.12**: *A loop L is said to be power-associative in the sense that every element of L generates an abelian group.*

**DEFINITION 2.2.13**: *A loop L is diassociative loop if every pair of elements of L generates a subgroup.*

*Example 2.2.3*: Let L be a loop given by the following table:

| • | e | $a_1$ | $a_2$ | $a_3$ | $a_4$ | $a_5$ |
|---|---|---|---|---|---|---|
| e | e | $a_1$ | $a_2$ | $a_3$ | $a_4$ | $a_5$ |
| $a_1$ | $a_1$ | e | $a_3$ | $a_5$ | $a_2$ | $a_4$ |
| $a_2$ | $a_2$ | $a_5$ | e | $a_4$ | $a_1$ | $a_3$ |
| $a_3$ | $a_3$ | $a_4$ | $a_1$ | e | $a_5$ | $a_2$ |
| $a_4$ | $a_4$ | $a_3$ | $a_5$ | $a_2$ | e | $a_1$ |
| $a_5$ | $a_5$ | $a_2$ | $a_4$ | $a_1$ | $a_3$ | e |

The nucleus of this loop is just {e}. The left nucleus of L, $N_\lambda(L) = \{e\}$. The Moufang centre of the loop L is C(L) = {e}. Thus for this L we see the center is just {e}.

The reader is requested to prove the above facts, which has been already verified by the author for this loop L.

**PROBLEMS:**

1. Give an example of a loop L of order 11, which has a non-trivial centre.
2. Find a loop L in which $N_\lambda(L) = N_\mu(L) = N_\rho(L) \neq \{e\}$, {e} the identity element of L.
3. Find a loop L in which C(L) ≠ e or C(L) ≠ L.
4. Give an example of a loop L in which Z(L) ≠ {e} and Z(L) ≠ L.
5. Find a loop L in which the commutator of L is different from {e}.
6. Give an example of a non-simple loop of order 13.
7. Find a loop L in which all subloops are normal.
8. Can there exist a loop L in which the associator subloop and the commutator subloop are equal but not equal to L?
9. For the loops given in examples 2.1.2 and 2.1.3 construct a loop homomorphism.
10. Give an example of a strict non-commutative loop.
11. Prove or disprove in a strict non-commutative loop the Moufang center, centre, $N_\lambda$, $N_\mu$ and $N_\rho$ are all equal and is equal to {e}.
12. Does there exist an example of a loop which has no subloops?

## 2.3 Special identities in loops

In this section we recall several special identities in loops introduced by Bruck, Bol, Moufang, Hamiltonian, etc and illustrate them with example whenever possible. As all these notions are to be given a Smarandache analogue in the next chapter we have tried our level best to give all the identities.

**DEFINITION 2.3.1**: *A loop L is said to be a Moufang loop if it satisfies any one of the following identities:*

1. *(xy) (zx) = (x(yz))x*
2. *((xy)z)y = x(y(zy))*
3. *x(y(xz)) = ((xy)x)z*

*for all x, y, z ∈ L.*

**DEFINITION 2.3.2**: *Let L be a loop, L is called a Bruck loop if x(yx)z = x(y(xz)) and $(xy)^{-1} = x^{-1}y^{-1}$ for all x, y, z ∈ L.*

**DEFINITION 2.3.3**: *A loop (L, •) is called a Bol loop if ((xy)z)y = x((yz)y) for all x, y, z ∈ L.*

**DEFINITION 2.3.4**: *A loop L is said to be right alternative if (xy)y = x(yy) for all x, y ∈ L and L is left alternative if (xx)y = x(xy) for all x, y ∈ L. L is said to be an alternative loop if it is both a right and left alternative loop.*

**DEFINITION 2.3.5**: *A loop (L, •) is called a weak inverse property loop (WIP-loop) if (xy)z = e imply x(yz) = e for all x, y, z ∈ L.*

**DEFINITION 2.3.6**: *A loop L is said to be semi alternative if (x, y, z) = (y, z, x) for all x, y, z ∈ L, where (x, y, z) denotes the associator of elements x, y, z ∈ L.*

**THEOREM** (Moufang's theorem): *Every Moufang loop G is diassociative more generally, if a, b, c are elements in a Moufang loop G such that (ab)c = a(bc) then a, b, c generate an associative loop.*

The proof is left for the reader; for assistance refer Bruck R.H. [6].

**PROBLEMS:**

1. Can a loop of order 5 be a Moufang loop?
2. Give an example of a strictly non-commutative loop of order 9.
3. Can a loop of order 11 be not simple?

4. Does there exist a loop of order p, p a prime that is simple?
5. Give an example of a loop, which is not a Bruck loop.
6. Does there exist an example of a loop L that is Bruck, Moufang and Bol?
7. Give an example of a power associative loop of order 14.

## 2.4 Special types of loops

In this section we introduce several special types of loops like unique product loop, two unique product loop, Hamiltonian loop, diassociative loop, strongly semi right commutative loop, inner commutative loop etc. The loops can be of any order finite or infinite. Further we recall these definitions and some simple interesting properties are given for the deeper understanding of these concepts. Several proofs are left for the reader to prove. The concept of unique product groups and two unique product groups were introduced in 1941 by Higman [28]. He studied this relative to zero divisors in group rings and in fact proved if G is a two unique product group or a unique product group then the group ring FG has no zero divisors where F is a field of characteristic zero. In 1980 Strojnowski.A [62] proved that in case of groups the notion of unique product and two unique product coincide. We introduce the definition of unique product (u.p) and two unique product (t.u.p) to loops.

**DEFINITION 2.4.1**: *Let L be a loop, L is said to be a two unique product loop (t.u.p) if given any two non-empty finite subsets A and B of L with |A| + |B| > 2 there exist at least two distinct elements x and y of L that have unique representation in the from x = ab and y = cd with a, c $\in$ A and b, d $\in$ B.*

*A loop L is called a unique product (u.p) loop if, given A and B two non-empty finite subsets of L, then there always exists at least one x $\in$ L which has a unique representation in the from x = ab, with a $\in$ A and b $\in$ B.*

It is left as an open problem to prove whether the two concepts u.p and t.u.p are one and the same in case of loops

**DEFINITION 2.4.2**: *An A-loop is a loop in which every inner mapping is an automorphism.*

It has been proved by Michael K. Kinyon et al that every diassociative A-loop is a Moufang loop. For proof the reader is requested to refer [41, 42]. As the main aim of this book is the introduction of Smarandache loops and study their Smarandache properties we have only recalled the definition from several authors. A few interesting results are given as exercise to the readers.

**DEFINITION 2.4.3**: *A loop L is said to be simple ordered by a binary relation (<) provided that for all a, b, c in L*

> (i) exactly one of the following holds: $a < b$, $a = b$, $b < a$.
> (ii) if $a < b$ and $b < c$ then $a < c$.
> (iii) if $a < b$ then $ac < bc$ and $ca < cb$.

*The relation $a > b$ is interpreted as usual, to mean $b < a$. If 1 is the identity element of L, an element a is called positive if $a > 1$, negative if $a < 1$. This notion will find its importance in case of loop rings.*

**DEFINITION 2.4.4**: *A loop L is called Hamiltonian if every subloop is normal.*

**DEFINITION 2.4.5**: *A loop is power associative (diassociative) if every element generates (every two element generate) a subgroup.*

**DEFINITION 2.4.6**: *A power associative loop is a p-loop (p a prime) if every element has p-power order.*

In view of these definitions we have the following theorems.

**THEOREM 2.4.1**: *A power associative Hamiltonian loop in which every element has finite order is a direct product of Hamiltonian p-loops.*

The proof of the theorem is left for the reader. Refer [6] for more information.

**THEOREM 2.4.2**: *A diassociative Hamiltonian loop G is either an abelian group or a direct product $G = A \otimes T \otimes H$ where A is an abelian group whose elements have finite odd order. T is an abelian group of exponent 2 and H is a non-commutative loop with the following properties*

> (i) The centre $Z = Z(H)$ has order two element 1, where $e \neq 1$, $e^2 = 1$.
> (ii) If x is a non-central element of H; $x^2 = e$.
> (iii) If x, y are in H and $(x, y) \neq e$ and x, y generate a quaternion group.
> (iv) If x, y, z are in H and $(x, y, z) \neq 1$ then $(x, y, z) = e$.

The proof is left for the reader and is requested to refer [6].

Norton [44,45] showed conversely if A, T, H are as specified in the above theorem then $A \otimes T \otimes H$ is a diassociative Hamiltonian loop. Further more if H satisfies the additional hypothesis that any three elements x, y, z for which $(x, y, z) = 1$ are contained in a subgroup then H is either a quaternion group or a Cayley group. We define CA-loop, inner commutative loops etc. and prove the following results.

**DEFINITION 2.4.7**: *Let L be a loop. An element x ∈ L is called a CA- element in L if (ax)b = (xb)a and a(xb) = b(ax) for all a, b ∈ L. A loop L having a CA-element is called a CA-loop.*

**DEFINITION 2.4.8**: *Let L be a loop. L is said to be inner commutative if every proper subloop of L is commutative but L is not commutative. We say L is strictly inner commutative if every proper subloop of L is commutative but they are not cyclic groups.*

**DEFINITION 2.4.9**: *Let L be a loop, L is said to be semi right commutative if for every pair of elements x, y ∈ L; we can always find an element c ∈ L such that ab = c(ba) or (cb)a.*

**DEFINITION 2.4.10**: *Let L be a loop if for every triple x, y, z in L at least one of the following equality is true:*

(i)     *xy = z(yx) or (zy)x*
(ii)    *yz = x(zy) or (xz)y*
(iii)   *zx = y(xz) or (yx)z*

*then we call the loop L to be a strongly semi right commutative.*

**THEOREM 2.4.3**: *Every strongly semi right commutative loop is commutative.*

*Proof*: Clearly for every triple {1, a, b}, 1, a, b ∈ L; 1 the unit element of L we have ab = ba. Hence the claim.

**THEOREM 2.4.4**: *Every commutative loop L in general need not be strongly semi right commutative.*

*Proof*: By an example. Consider the loop L = {e, a, b, c, d, g} given by the following table; e is the identity element of L.

| • | e | a | b | c | d | g |
|---|---|---|---|---|---|---|
| e | e | a | b | c | d | g |
| a | a | e | d | b | g | c |
| b | b | d | e | g | c | a |
| c | c | b | g | e | a | d |
| d | d | g | c | a | e | b |
| g | g | c | a | d | b | e |

Clearly L is commutative. For the triple a, b, c ∈ L, ab ≠ c (ba) or (cb)a for ab = d and c(ba) = cd = a; d ≠ a, hence ab ≠ c(ba).

Now (cb)a = ga = c; ab = d ≠ (cb)a = c. So L is not strongly semi right commutative.

**THEOREM 2.4.5**: *Every commutative loop L is semi right commutative.*

*Proof*: Clear from the fact for every pair x, y ∈ L we can choose the identity element e of L so that xy = e(yx).

**THEOREM 2.4.6**: *Every semi right commutative loop is not strongly semi right commutative.*

*Proof*: The loop L given as an example in Theorem 2.4.4 is commutative, but L is not strongly semi right commutative.

We have the following relation

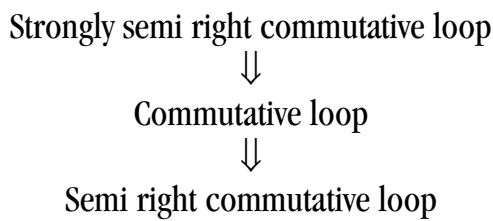

Strongly semi right commutative loop
⇓
Commutative loop
⇓
Semi right commutative loop

**DEFINITION 2.4.11**: *Let L be a loop. If for a, b ∈ L with ab = ba in L if we have (ax)b = (bx)a (or b(xa)) for all x ∈ L then we say the pair a, b is pseudo commutative. [we can have a(xb) = (bx)a or b(xa) that is we can take a(xb) instead of (ax)b also]. If in L every commutative pair is pseudo commutative then we say the loop L is a pseudo commutative loop. If for every distinct pair of elements a, b ∈ L we have axb = bxa for all x ∈ L take the associative brackets in a way which suits the equality; then the loop is called a strongly pseudo commutative loop. Using this concept of pseudo commutativity we now define the concept of pseudo commutator of L.*

**DEFINITION 2.4.12**: *Let L be a loop. The pseudo commutator of L demoted by P(L) = 〈{ p ∈ L / a(xb) = p [(bx)a]}〉 where 〈〉 denotes the set generated by p. Similarly we define the strongly pseudo commutator of L denoted by SP(L) = 〈{p ∈ L / (ax)b = (pb) (ax)}〉.*

**DEFINITION 2.4.13**: *Let L be a loop. An associative triple a, b, c ∈ L is said to be pseudo associative if (ab) (xc) = (ax) (bc) for all x ∈ L. If (ab) (xc) = (ax) (bc) for some x ∈ L we say the triple is pseudo associative relative to those x in L. If in particular associative triple is pseudo associative then we say the loop is a pseudo associative loop.*

*If in the non-associative loop L if for a non-associative triple a, b, c in L we have (ax) (bc) = (ab) (xc) then for all $x \in L$, we say the triple is strongly pseudo associative. A loop is strongly pseudo associative if every triple is strongly pseudo associative.*

**DEFINITION 2.4.14**: *Let L be a loop. PA(L) = $\langle \{t \in L$ / (ab) (tc) = (at) (bc) where a(bc) = (ab)c for a, b, c $\in L\} \rangle$ denotes the pseudo associator of L generated by $t \in L$, satisfying the condition given in PA(L).*

Similarly we can define strongly pseudo associator SPA(L) of the loop L. Now we proceed onto another type of loop called Jordan loop.

**DEFINITION [78]**: *Let L be a loop. We say L is a Jordan loop if ab = ba. $a^2$ (ba) = ($a^2$b) a for all a, b $\in L$.*

*Example 2.4.1*: Let L = {e, $g_1$, ... , $g_7$} given the following table:

| • | e | $g_1$ | $g_2$ | $g_3$ | $g_4$ | $g_5$ | $g_6$ | $g_7$ |
|---|---|---|---|---|---|---|---|---|
| e | e | $g_1$ | $g_2$ | $g_3$ | $g_4$ | $g_5$ | $g_6$ | $g_7$ |
| $g_1$ | $g_1$ | e | $g_5$ | $g_2$ | $g_6$ | $g_3$ | $g_7$ | $g_4$ |
| $g_2$ | $g_2$ | $g_5$ | e | $g_6$ | $g_3$ | $g_7$ | $g_4$ | $g_1$ |
| $g_3$ | $g_3$ | $g_2$ | $g_6$ | e | $g_7$ | $g_4$ | $g_1$ | $g_3$ |
| $g_4$ | $g_4$ | $g_6$ | $g_3$ | $g_7$ | e | $g_1$ | $g_5$ | $g_2$ |
| $g_5$ | $g_5$ | $g_3$ | $g_7$ | $g_4$ | $g_1$ | e | $g_2$ | $g_6$ |
| $g_6$ | $g_6$ | $g_7$ | $g_4$ | $g_1$ | $g_5$ | $g_2$ | e | $g_3$ |
| $g_7$ | $g_7$ | $g_4$ | $g_1$ | $g_3$ | $g_2$ | $g_6$ | $g_3$ | e |

This is easily verified to be a Jordan loop of order 8.

**DEFINITION [65]**: *A loop L is said to be strictly non-right alternative if (xy)y $\neq$ x(yy) for any distinct pair x, y in L with x $\neq$ e and y $\neq$ e. Similarly we define strictly non-left alternative loop. We say a loop is strictly non-alternative if L is simultaneously a strictly non-left alternative and strictly non-right alternative loop.*

The following examples give a strictly non-left alternative and a strictly non-right alternative loops each of order 6.

*Example 2.4.2* : Let $L_5(2)$ be a loop in $L_5$ given by the following table:

|   | e | 1 | 2 | 3 | 4 | 5 |
|---|---|---|---|---|---|---|
| e | e | 1 | 2 | 3 | 4 | 5 |
| 1 | 1 | e | 3 | 5 | 2 | 4 |
| 2 | 2 | 5 | e | 4 | 1 | 3 |
| 3 | 3 | 4 | 1 | e | 5 | 2 |
| 4 | 4 | 3 | 5 | 2 | e | 1 |
| 5 | 5 | 2 | 4 | 1 | 3 | e |

This loop is strictly non-left alternative.

***Example 2.4.3***: Let $L_5(4)$ be a loop of order 6 given by the following table:

| • | e | 1 | 2 | 3 | 4 | 5 |
|---|---|---|---|---|---|---|
| e | e | 1 | 2 | 3 | 4 | 5 |
| 1 | 1 | e | 5 | 4 | 3 | 2 |
| 2 | 2 | 3 | e | 1 | 5 | 4 |
| 3 | 3 | 5 | 4 | e | 2 | 1 |
| 4 | 4 | 2 | 1 | 5 | e | 3 |
| 5 | 5 | 4 | 3 | 2 | 1 | e |

This loop is strictly non-right alternative.

**PROBLEMS:**

1. Give an example of a Jordan loop of order 7.
2. Can a strongly pseudo associative loop of order 8 exist?
3. Is the loop given by the following table:

| • | e | $a_1$ | $a_2$ | $a_3$ | $a_4$ | $a_5$ |
|---|---|---|---|---|---|---|
| e | e | $a_1$ | $a_2$ | $a_3$ | $a_4$ | $a_5$ |
| $a_1$ | $a_1$ | e | $a_5$ | $a_4$ | $a_3$ | $a_2$ |
| $a_2$ | $a_2$ | $a_3$ | e | $a_1$ | $a_5$ | $a_4$ |
| $a_3$ | $a_3$ | $a_5$ | $a_4$ | e | $a_2$ | $a_1$ |
| $a_4$ | $a_4$ | $a_2$ | $a_1$ | $a_5$ | e | $a_3$ |
| $a_5$ | $a_5$ | $a_4$ | $a_3$ | $a_2$ | $a_1$ | e |

   a. A Jordan loop?
   b. strongly pseudo commutative? diassociative? power associative?
   c. CA-loop? Strictly inner commutative loop?

4. Can a loop of order p , p a prime be a power associative loop? Justify your answer.

5. Is a loop of order 19 diassociative? If so give an example.
6. Give an example of a strongly semi right commutative loop.
7. What is the smallest order of a loop L so that L is diassociative?
8. Can an ordered loop of order 12 be a Hamiltonian loop?
9. Give an example of a loop which is not a Hamiltonian loop.
10. Find a loop of odd order say 15 which is Hamiltonian and not a A-loop. (Hint: If no solution exists for these problems can we give nice characterization theorem about inter relation of these concepts).
11. Find an example of a loop L in which A(L) = SPA(L).
12. Does there exists a loop L such that A(L) = PA(L)?

## 2.5 Representation and isotopes of loops

This section is mainly devoted to give the definition of right regular representation of loops. We also recall the definition of isotopes and the concept of G-loops. As we assume the reader to have a good knowledge of not only algebra but a very strong foundation about loops we just recall the definition; as our motivation is the introduction and study of Smarandache loops. Here the right regular representation or in short representation of loops are given. For a detailed notion about these concepts the reader is requested to refer Albert and Burn [1, 8, 9, 10].

***Result:*** [1]: A set $\pi$ of permutations on a set L is the representation of a loop (L, •) if and only if

    a. $I \in \pi$
    b. $\pi$ is transitive on L and
    c. for $\alpha, \beta \in \pi$, if $\alpha\beta^{-1}$ fixes one element of L, then $\alpha = \beta$.

Let (L, •) be a finite loop. For $\alpha \in L$, define a right multiplication $R_\alpha$ as a permutation of the loop (L, •) as follows:

$R_\alpha : x \to x \bullet \alpha$ we will call the set $\{R_\alpha / \alpha \in L\}$ the right regular representation of (L, •) or briefly the representation of L.

***Result:*** [8, 9, 10]: If $\pi$ is a representation of a loop L, then L is a Bol loop if and only if $\alpha, \beta \in \pi$ implies $\alpha\beta\alpha \in \pi$.

***Result:*** [8, 9, 10]: If $\pi$ is the representation of a Bol loop and if $\alpha \in \pi$ then $\alpha^n \in \pi$ for $n \in Z$.

***Result:*** [8, 9, 10]: If $\pi$ is the representation of a finite Bol loop for $\alpha \in \pi$, $\alpha$ is a product of disjoint cycles of equal length.

*Example 2.5.1*: Let L be a loop given the following table:

| • | e | $g_1$ | $g_2$ | $g_3$ | $g_4$ | $g_5$ | $g_6$ | $g_7$ |
|---|---|---|---|---|---|---|---|---|
| e | e | $g_1$ | $g_2$ | $g_3$ | $g_4$ | $g_5$ | $g_6$ | $g_7$ |
| $g_1$ | $g_1$ | e | $g_5$ | $g_2$ | $g_6$ | $g_3$ | $g_7$ | $g_4$ |
| $g_2$ | $g_2$ | $g_5$ | e | $g_6$ | $g_3$ | $g_7$ | $g_4$ | $g_1$ |
| $g_3$ | $g_3$ | $g_2$ | $g_6$ | e | $g_7$ | $g_4$ | $g_1$ | $g_5$ |
| $g_4$ | $g_4$ | $g_6$ | $g_3$ | $g_7$ | e | $g_1$ | $g_5$ | $g_2$ |
| $g_5$ | $g_5$ | $g_3$ | $g_7$ | $g_4$ | $g_1$ | e | $g_2$ | $g_6$ |
| $g_6$ | $g_6$ | $g_7$ | $g_4$ | $g_1$ | $g_5$ | $g_2$ | e | $g_3$ |
| $g_7$ | $g_7$ | $g_4$ | $g_1$ | $g_3$ | $g_2$ | $g_6$ | $g_3$ | e |

The right regular representation of the loop L is given by

$$I$$
$$(eg_1)\ (g_2g_5g_3)\ (g_4g_6g_7)$$
$$(eg_2)\ (g_1g_5g_7)\ (g_3g_6g_4)$$
$$(eg_3)\ (g_1g_2g_6)\ (g_4g_7g_5)$$
$$(eg_4)\ (g_1g_6g_5)\ (g_2g_3g_7)$$
$$(eg_5)\ (g_1g_3g_4)\ (g_2g_7g_6)$$
$$(eg_6)\ (g_1g_7g_3)\ (g_2g_4g_5)$$
$$(eg_7)\ (g_1g_4g_2)\ (g_3g_5g_6)$$

where I is the identity permutation on the loop L.

Now we go on to give the definition of isotopes for more information please refer Bruck [6].

**DEFINITION 2.5.1**: *Let (L, •) be a loop. The principal isotope (L, *) of (L, •) with respect to any predetermined a, b ∈ L is defined by x * y = XY, for all x, y ∈ L, where Xa = x and bY = y for some X, Y ∈ L.*

**DEFINITION 2.5.2**: *Let L be a loop, L is said to be a G-loop if it is isomorphic to all of its principal isotopes.*

**PROBLEMS:**

1. Does there exist a loop L of order 7 which is a G-loop?
2. Can we have loops L of odd order n, n finite such that L is a G-loop?
3. Give an example of a loop L of order 19 which in not a G-loop.
4. Which class of loops are G-loops? Moufang? Bruck? Bol?
5. Is an alternative loop of order 10 a G-loop?

## 2.6 On a new class of loops and its properties.

The main objective of this section is the introduction of a new class of loops with a natural simple operation. As to introduce loops several functions or maps are defined satisfying some desired conditions we felt that it would be nice if we can get a natural class of loops built using integers.

Here we define the new class of loops of any even order, they are distinctly different from the loops constructed by other researchers. Here we enumerate several of the properties enjoyed by these loops.

**DEFINITION [65]:** *Let $L_n(m) = \{e, 1, 2, \ldots, n\}$ be a set where $n > 3$, $n$ is odd and $m$ is a positive integer such that $(m, n) = 1$ and $(m-1, n) = 1$ with $m < n$.*

*Define on $L_n(m)$ a binary operation '•' as follows:*

1. *$e \bullet i = i \bullet e = i$ for all $i \in L_n(m)$*
2. *$i^2 = i \bullet i = e$ for all $i \in L_n(m)$*
3. *$i \bullet j = t$ where $t = (mj - (m-1)i) \pmod{n}$*

*for all $i, j \in L_n(m)$; $i \neq j$, $i \neq e$ and $j \neq e$, then $L_n(m)$ is a loop under the binary operation '•'.*

***Example 2.6.1***: Consider the loop $L_5(2) = \{e, 1, 2, 3, 4, 5\}$. The composition table for $L_5(2)$ is given below:

| • | e | 1 | 2 | 3 | 4 | 5 |
|---|---|---|---|---|---|---|
| e | e | 1 | 2 | 3 | 4 | 5 |
| 1 | 1 | e | 3 | 5 | 2 | 4 |
| 2 | 2 | 5 | e | 4 | 1 | 3 |
| 3 | 3 | 4 | 1 | e | 5 | 2 |
| 4 | 4 | 3 | 5 | 2 | e | 1 |
| 5 | 5 | 2 | 4 | 1 | 3 | e |

This loop is of order 6 which is both non-associative and non-commutative.

*Physical interpretation of the operation in the loop $L_n(m)$:*

We give a physical interpretation of this class of loops as follows: Let $L_n(m) = \{e, 1, 2, \ldots, n\}$ be a loop in this identity element of the loop are equidistantly placed on a circle with e as its centre.

We assume the elements to move always in the clockwise direction.

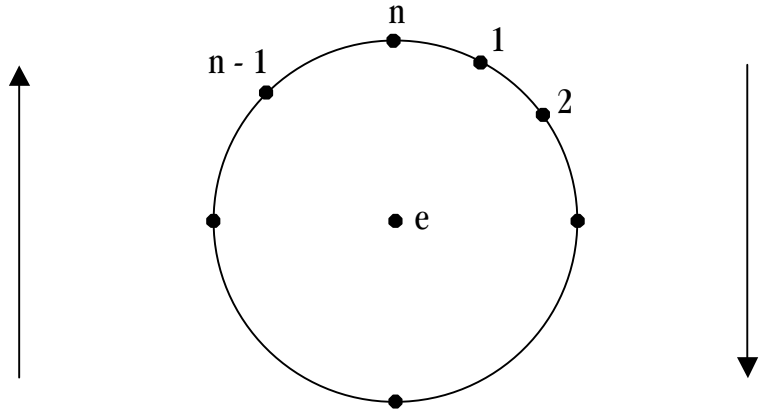

Let i, j ∈ $L_n(m)$ (i ≠ j, i ≠ e, j ≠ e). If j is the $r^{th}$ element from i counting in the clockwise direction the i • j will be the $t^{th}$ element from j in the clockwise direction where t = (m −1)r. We see that in general i • j need not be equal to j • i. When i = j we define $i^2$ = e and i • e = e • i = i for all i ∈ $L_n(m)$ and e acts as the identity in $L_n(m)$.

***Example 2.6.2***: Now the loop $L_7(4)$ is given by the following table:

| • | e | 1 | 2 | 3 | 4 | 5 | 6 | 7 |
|---|---|---|---|---|---|---|---|---|
| e | e | 1 | 2 | 3 | 4 | 5 | 6 | 7 |
| 1 | 1 | e | 5 | 2 | 6 | 3 | 7 | 4 |
| 2 | 2 | 5 | e | 6 | 3 | 7 | 4 | 1 |
| 3 | 3 | 2 | 6 | e | 7 | 4 | 1 | 5 |
| 4 | 4 | 6 | 3 | 7 | e | 1 | 5 | 2 |
| 5 | 5 | 3 | 7 | 4 | 1 | e | 2 | 6 |
| 6 | 6 | 7 | 4 | 1 | 5 | 2 | e | 3 |
| 7 | 7 | 4 | 1 | 5 | 2 | 6 | 3 | e |

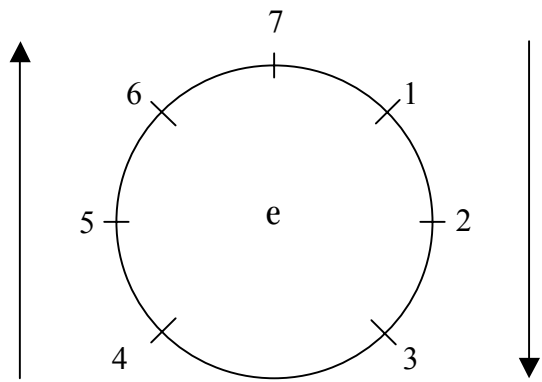

Let 2, 4 ∈ $L_7(4)$. Now 4 is the $2^{nd}$ element from 2 in the clockwise direction. So 2.4 will be (4 −1)2 that is the $6^{th}$ element from 4 in the clockwise direction which is 3.

Hence 2.4 = 3.

**Notation**: Let $L_n$ denote the class of loops. $L_n(m)$ for fixed n and various m's satisfying the conditions $m < n$, $(m, n) = 1$ and $(m - 1, n) = 1$, that is $L_n = \{L_n(m) \mid n > 3, n \text{ odd}, m < n, (m, n) = 1 \text{ and } (m-1, n) = 1\}$.

***Example 2.6.3***: Let n = 5. The class $L_5$ contains three loops; viz. $L_5(2)$, $L_5(3)$ and $L_5(4)$ given by the following tables:

$L_5(2)$

| • | e | 1 | 2 | 3 | 4 | 5 |
|---|---|---|---|---|---|---|
| e | e | 1 | 2 | 3 | 4 | 5 |
| 1 | 1 | e | 3 | 5 | 2 | 4 |
| 2 | 2 | 5 | e | 4 | 1 | 3 |
| 3 | 3 | 4 | 1 | e | 5 | 2 |
| 4 | 4 | 3 | 5 | 2 | e | 1 |
| 5 | 5 | 2 | 4 | 1 | 3 | e |

$L_5(3)$

| • | e | 1 | 2 | 3 | 4 | 5 |
|---|---|---|---|---|---|---|
| e | e | 1 | 2 | 3 | 4 | 5 |
| 1 | 1 | e | 4 | 2 | 5 | 3 |
| 2 | 2 | 4 | e | 5 | 3 | 1 |
| 3 | 3 | 2 | 5 | e | 1 | 4 |
| 4 | 4 | 5 | 3 | 1 | e | 2 |
| 5 | 5 | 3 | 1 | 4 | 2 | e |

$L_5(4)$

| • | e | 1 | 2 | 3 | 4 | 5 |
|---|---|---|---|---|---|---|
| e | e | 1 | 2 | 3 | 4 | 5 |
| 1 | 1 | e | 5 | 4 | 3 | 2 |
| 2 | 2 | 3 | e | 1 | 5 | 4 |
| 3 | 3 | 5 | 4 | e | 2 | 1 |
| 4 | 4 | 2 | 1 | 5 | e | 3 |
| 5 | 5 | 4 | 3 | 2 | 1 | e |

**THEOREM [56]**: *Let $L_n$ be the class of loops for any $n > 3$, if $n = p_1^{\alpha_1} p_2^{\alpha_2} \ldots p_k^{\alpha_k}$ ($\alpha_i > 1$, for $i = 1, 2, \ldots, k$), then $|L_n| = \prod_{i=1}^{k}(p_i - 2) p_i^{\alpha_i - 1}$ where $|L_n|$ denotes the number of loops in $L_n$.*

The proof is left for the reader as an exercise.

**THEOREM [56]**: *$L_n$ contains one and only one commutative loop. This happens when $m = (n + 1) / 2$. Clearly for this m, we have $(m, n) = 1$ and $(m - 1, n) = 1$.*

It can be easily verified by using simple number theoretic techniques.

**THEOREM [56]**: *Let $L_n$ be the class of loops. If $n = p_1^{\alpha_1} p_2^{\alpha_2} \ldots p_k^{\alpha_k}$, then $L_n$ contains exactly $F_n$ loops which are strictly non-commutative where $F_n = \prod_{i=1}^{k}(p_i - 3) p_i^{\alpha_i - 1}$.*

The proof is left for the reader as an exercise.

<u>Note</u>: If $n = p$ where p is a prime greater than or equal to 5 then in $L_n$ a loop is either commutative or strictly non-commutative. Further it is interesting to note if $n = 3t$ then the class $L_n$ does not contain any strictly non-commutative loop.

**THEOREM [65]**: *The class of loops $L_n$ contains exactly one left alternative loop and one right alternative loop but does not contain any alternative loop.*

*Proof*: We see $L_n(2)$ is the only right alternative loop that is when $m = 2$ (Left for the reader to prove using simple number theoretic techniques). When $m = n - 1$ that is $L_n(n - 1)$ is the only left alternative loop in the class of loops $L_n$.

From this it is impossible to find a loop in $L_n$, which is simultaneously right alternative and left alternative. Further it is clear from earlier result both the right alternative loop and the left alternative loop is not commutative.

**THEOREM [56]**: *Let $L_n$ be the class of loops:*

1. *$L_n$ does not contain any Moufang loop*
2. *$L_n$ does not contain any Bol loop*
3. *$L_n$ does not contain any Bruck loop.*

The reader is requested to prove these results using number theoretic techniques.

**THEOREM [65]**: *Let $L_n(m) \in L_n$. Then $L_n(m)$ is a weak inverse property (WIP) loop if and only if $(m^2 - m + 1) \equiv 0 \pmod{n}$.*

*Proof*: It is easily checked that for a loop to be a WIP loop we have "if $(xy)z = e$ then $x(yz) = e$ where $x, y, z \in L$." Both way conditions can be derived using the defining operation on the loop $L_n(m)$.

***Example 2.6.4***: L be the loop $L_7(3) = \{e, 1, 2, 3, 4, 5, 6, 7\}$ be in $L_7$ given by the following table:

| • | e | 1 | 2 | 3 | 4 | 5 | 6 | 7 |
|---|---|---|---|---|---|---|---|---|
| e | e | 1 | 2 | 3 | 4 | 5 | 6 | 7 |
| 1 | 1 | e | 4 | 7 | 3 | 6 | 2 | 5 |
| 2 | 2 | 6 | e | 5 | 1 | 4 | 7 | 3 |
| 3 | 3 | 4 | 7 | e | 6 | 2 | 5 | 1 |
| 4 | 4 | 2 | 5 | 1 | e | 7 | 3 | 6 |
| 5 | 5 | 7 | 3 | 6 | 2 | e | 1 | 4 |
| 6 | 6 | 5 | 1 | 4 | 7 | 3 | e | 2 |
| 7 | 7 | 3 | 6 | 2 | 5 | 1 | 4 | e |

It is easily verified $L_7(3)$ is a WIP loop. One way is easy for $(m^2 - m + 1) \equiv 0 \pmod{7}$ that is $9 - 3 + 1 = 9 + 4 + 1 \equiv 0 \pmod{7}$. It is interesting to note that no loop in the class $L_n$ contain any associative loop.

**THEOREM [56]**: *Let $L_n$ be the class of loops. The number of strictly non-right (left) alternative loops is $P_n$ where $P_n = \prod_{i=1}^{k} (p_i - 3) p_i^{\alpha_i - 1}$ and $n = \prod_{i=1}^{k} p_i^{\alpha_i}$.*

The proof is left for the reader to verify.

Now we proceed on to study the associator and the commutator of the loops in $L_n$.

**THEOREM [56]**: *Let $L_n(m) \in L_n$. The associator $A(L_n(m)) = L_n(m)$.*

<u>Hint</u>: Recall.

$A(L_n(m)) = \langle \{ t \in L_n(m) \,/\, t = (x, y, z) \text{ for some } x, y, z \in L_n(m) \} \rangle$. Construct t using the definition of $L_n(m)$ it is easily verified $A(L_n(m)) = L_n(m)$.

Similarly if we take $L_n(m)$ any non-commutative loop in $L_n(m)$, we see $L'_n(m) = L_n(m)$ where $L'_n(m) = \langle \{ t \in L_n(m) \,/\, t = (x, y) \text{ for some } x, y \in L_n(m) \} \rangle$. To prove $L_n'(m) = L_n(m)$ it is sufficient to show that for every $t \in L_n(m)$ there exists $i, j \in L_n(m)$ such that $(i, j) = t$. Thus we see in case of the new class of loops $L_n$. $A(L_n(m)) = L'_n(m) = L_n(m)$ for meaningful $L_n(m)$. One of the very interesting properties about this new class of loops is their subloops. We see these subloops enjoy some special properties.

Let $L_n(m) \in L_n$, for every t / n there are t subloops of order k + 1 where k = n/t. Let $L_n(m) = \{e, 1, 2, \ldots, n\}$, for $i \leq t$ consider the subset. $H_i(t) = \{e, i, i + t, i + 2t, \ldots, i + (k - 1) t\}$ of $L_n(m)$. Clearly $e \in H_i$; it is left for the reader to prove $H_i(t)$ is a subloop of $L_n(m)$ of order k + 1. As i can vary from 1 to t, there exists t such subloops. Let $H_i(t)$ be subloops for $i \neq j$ we have $H_i(t) \cap H_j(t) = \{e\}$. Also we have $L_n(m) = \bigcup_{i=1}^{t} H_i(t)$ for every t dividing n. Finally subloops $H_i(t)$ and $H_j(t)$ are isomorphic for every t dividing n.

Suppose H is any subloop of $L_n(m)$ of order t + 1 then t / n in view of this we have $L_n(m) \in L_n$ contains a subloop of order k + 1 if and only if k/n. Now we want to find out for what values of n; the Lagrange's theorem for groups is satisfied by every subloop of $L_n(m)$.

Let $L_n(m) \in L_n$. The Lagrange's theorem for groups is satisfied by every subloop of $L_n(m)$ if and only if n is an odd prime. This is easily verified when n is a prime we have $|L_n(m)| = p + 1$ as the order of the subloop is p + 1 as no t/p so only 1/p so there exists a subloop of order 2 and as p / p hence there exists a subloop of order p +1. Thus trivially the Lagrange's theorem for groups is satisfied only by loops $L_n(m)$ of order p + 1 where p is a prime ( n = p ). Further we have for any loop $L_n(m)$ of $L_n$ there exist only 2-sylow subloops. Clearly $|L_n(m)| = n + 1$ where n is odd (n > 3). As 2/n + 1 so $L_n(m)$ has 2-sylow subloops. Also if $L_n(m) \in L_n$; n is an odd prime say n = p then we see that there exist no element of order p in $L_n(m)$. In view of this property we state as "If $L_n(m) \in L_n$, n an odd prime then no element in $L_n(m)$ satisfies Cauchy Theorem for groups. It is left for the reader to show that no loop in the class of loops $L_n$ has a normal subloop; this leads us to frame a nice property that all loops in the class $L_n$ are simple.

***Example 2.6.5***: The subloops of the loop $L_5(2)$ are $\{e\}$, $\{e, 1\}$, $\{e, 2\}$, $\{e, 3\}$, $\{e, 4\}$, $\{e, 5\}$ and $L_5(2)$.

The lattice of subloops of $L_5(2)$

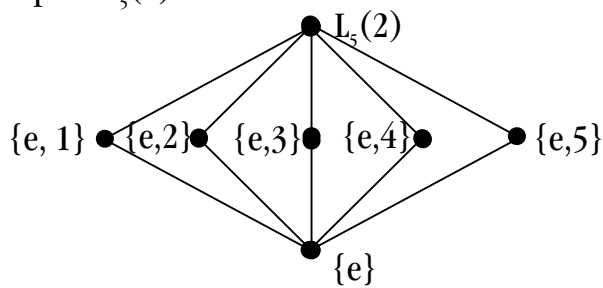

Here it is important to note that all these proper subloops are subgroups. Thus the lattice of subloops (subgroups) forms a modular lattice. We see $L_n(m) \in L_n$ are all S-loops. Thus we define in case of S-loops the Smarandache subgroups, Smarandache subloops and Smarandache normal subloops. We call the related lattice of substructures as S-lattices. In view of this we prove the following theorem.

**THEOREM 2.6.1**: *Let $L_n(m) \in L_n$ ( n a prime.) The subloops (or subgroups) of $L_n(m)$ form a non-distributive modular lattice of order n + 2 where each element* *is complement of the other and each chain connecting {e} and $L_n(m)$ is of length 3.*

*Proof:* We know the loops $L_n(m) \in L_n$ (n-prime) have no subloops other than the subgroups. $L_n(m)$, {e}, $A_i$ = {e, i}, i = 1, 2, ..., n. The following is the lattice of subloops:

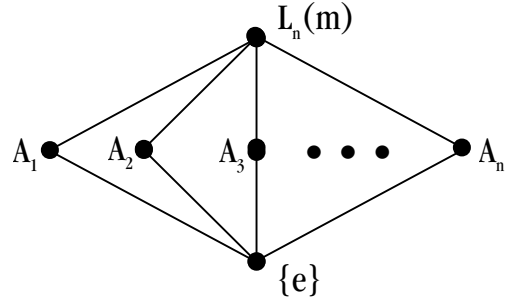

The order of this loop is n + 2.

Finally before we proceed to define about normalizers of the loops we see by the very definition of the loops $L_n(m)$ in the class of loop $L_n$ that all loops are power associative as every element in the loop $L_n(m)$ is such that i • i = e for all i ∈ $L_n(m)$ that is every element generates a cyclic group of order two; now it is left for the reader to prove that no loop in this class $L_n$ is diassociative.

We in case of loops $L_n(m) \in L_n$ define two normalizers: First normalizer and Second normalizer.

**DEFINITION [56]**: *Let $L_n(m) \in L_n$ and $H_i(t)$ be a subloop of $L_n(m)$. The first normalizer of $H_i(t)$ is given by $N_1(H_i(t)) = \{a \in L_n(m) / a(H_i(t)) = (H_i(t))a\}$ the second normalizer $N_2(H_i(t)) = \{x \in L_n(m) / x(H_i(t))x = H_i(t)\}$.*

*Example 2.6.6*: Take $L_{15}(2) \in L_{15}$. Consider the subloop $H_1(3)$ = {e, 1, 4, 7, 10, 13}, it is easily verified $N_1(H_1(3)) = L_{15}(2)$.

It is important question to analyse, whether the two normalizers are equal or in general distinct. The answer is the two normalizers in general are different for a subloop in $L_n(m)$.

We prove this by an example.

*Example 2.6.7*: Let $L_{45}(8) \in L_{45}$ be a loop of order 46. Consider $H_1(15)$ be a subloop of $L_{45}(8)$. $N_1(H_1(15))$ can be found to equal to $L_{45}(8)$ that is the whole loop. But $N_2(H_1(15) = H_1(5)$. So the two normalizers are not always equal.

Now the natural concept would be the class equation for groups. The class equation for groups is not satisfied by any loop in $L_n(m)$ except only by the commutative loop $L_n(m) \in L_n$. Now we discuss about the isotopes of the loops $L_n(m)$. Let $(L, \bullet)$ be a loop. The principal isotope $(L_1, *)$ of $(L, \bullet)$ with respect to any predetermined a, b ∈ L is defined by x * y = XY, for all x, y ∈ L, where Xa = x and bY = y for some X, Y ∈ L. We give the following examples of loops and their principal isotopes.

***Example 2.6.8***: Let $(L_5(2), \bullet)$ be the loop and its principal isotope is $(L_5(2), *)$, is given below:

Table of $(L_5(2), \bullet)$:

| $\bullet$ | e | 1 | 2 | 3 | 4 | 5 |
|---|---|---|---|---|---|---|
| e | e | 1 | 2 | 3 | 4 | 5 |
| 1 | 1 | e | 3 | 5 | 2 | 4 |
| 2 | 2 | 5 | e | 4 | 1 | 3 |
| 3 | 3 | 4 | 1 | e | 5 | 2 |
| 4 | 4 | 3 | 5 | 2 | e | 1 |
| 5 | 5 | 2 | 4 | 1 | 3 | e |

Table of $(L_5(2), *)$:

| * | e | 1 | 2 | 3 | 4 | 5 |
|---|---|---|---|---|---|---|
| e | 3 | 2 | 5 | e | 1 | 4 |
| 1 | 5 | 3 | 4 | 1 | e | 2 |
| 2 | 4 | e | 3 | 2 | 5 | 1 |
| 3 | e | 1 | 2 | 3 | 4 | 5 |
| 4 | 2 | 5 | 1 | 4 | 3 | e |
| 5 | 1 | 4 | e | 5 | 2 | 3 |

The main observation is that the principal isotopes do not in general preserve basic properties of the loops. This is illustrated by the following:

The principal isotope of a commutative loop can also be strictly non-commutative.

***Example 2.6.9***: $L_5(3) \in L_5$ is a commutative loop. Take a = 4 and b = e.

Composition table of $(L_5(3), \bullet)$, which is a commutative loop :

| $\bullet$ | e | 1 | 2 | 3 | 4 | 5 |
|---|---|---|---|---|---|---|
| e | e | 1 | 2 | 3 | 4 | 5 |
| 1 | 1 | e | 4 | 2 | 5 | 3 |
| 2 | 2 | 4 | e | 5 | 3 | 1 |
| 3 | 3 | 2 | 5 | e | 1 | 4 |
| 4 | 4 | 5 | 3 | 1 | e | 2 |
| 5 | 5 | 3 | 1 | 4 | 2 | e |

The principle isotope $(L_5(3), *)$ of $(L_5(3), \bullet)$ with respect to $a = 4$ and $b = e$ is strictly non-commutative is given by the following table:

Composition table of $(L_5(3), *)$

| * | e | 1 | 2 | 3 | 4 | 5 |
|---|---|---|---|---|---|---|
| e | 4 | 5 | 3 | 1 | e | 2 |
| 1 | 3 | 2 | 5 | e | 1 | 4 |
| 2 | 5 | 3 | 1 | 4 | 2 | e |
| 3 | 2 | 4 | e | 5 | 3 | 1 |
| 4 | e | 1 | 2 | 3 | 4 | 5 |
| 5 | 1 | e | 4 | 2 | 5 | 3 |

It is left for the reader to prove no loop in class $L_n$ is a G-loop.

**PROBLEMS:**

1. Find the commutative loop in $L_{27}$.
2. Let $L_{15}(2) \in L_{15}$. Find for the subloop $H_2(3)$ the first normalizer and the second normalizer.
3. Find the strictly non-commutative loop in $L_{129}$.
4. Find the WIP loop in $L_{211}$.
5. Find left semi alternative loop in $L_{301}$.
6. Find the right semi alternative loop in $L_{203}$.
7. Does $L_{203}(2)$ have normal subloops? Justify.
8. Find two subloops in $L_{33}(13)$ which are not subgroups of order 2.
9. Find the first normalizer of $H_1(15)$, which is a subloop in $L_{45}(8)$.
10. Find the associator of $L_{13}(2)$, that is $A(L_{13}(2))$.
11. Find the commutator subloop of $L_{17}(2)$.
12. Find the Moufang centre of $L_{17}(5)$.
13. Find the centre of $L_{19}(3)$, that is $Z(L_{19}(3))$.
14. Find $N(L_{23}(4))$ where $N(L)$ denotes the nucleus of $L_{23}(4)$.
15. Find elements in the loop $L_{25}(7)$ which satisfy Cauchy's theorem. (Hint: $x \in L_{25}(7)$ $x^t = e$ then $t / 26$).
16. Show the loop $L_{27}(5)$ is not a Bruck loop.
17. Prove the loop $L_9(5)$ is not a Moufang loop.
18. Can $L_9(5)$ be a Bol loop? Justify your answer.
19. Find all subloops of $L_9(5)$.

## 2.7 The new class of loops and its application to proper edge colouring of the graph $K_{2n}$

In this section we study the right regular representation of the new class of loops introduced in section 2.6. We formulate the necessary condition for the existence of k-cycle in the representation of $L_n(m)$ and using this we find out the cycle decomposition of $L_n(m)$. Further, we prove that the number of different representations of loops of order 2n (n ≥ 3), which are right alternative and in which square of each element is identity, is equal to the different proper edge colouring of the graph $K_{2n}$, using exactly (2n − 1) colours. As this application happens to be an important one we prove each and every other concepts needed first.

**THEOREM [65]**: *Let $L_n(m) \in L_n$. Then for any $a \in L_n(m)$ (a ≠ e) the transposition (a, e) belongs to the permutation $R_a$.*

*Proof*: Since in $L_n(m)$ we have a • a = e, $R_a$ : e → e • a = a and $R_a$ : a → a • a = e. Hence the transposition (a, e) belongs to the permutation $R_a$.

**THEOREM [65]**: *Let $L_n(m) \in L_n$ and $a \in L_n(m)$ (a ≠ e). If the permutation $R_a$ contains a k-cycle then, $((m − 1)^k + (− 1)^{k−1}) (a − x) \equiv 0 \pmod{n}$ where x is any element in the k-cycle.*

*Proof*: Consider the permutation $R_a$ (a ∈ $L_n(m)$ and a ≠ e). Let the permutation $R_a$ contain a k-cycle say $(x_1, x_2, \ldots, x_k)$ (as we are taking the cycle decomposition of elements other than a and e, $x_i \neq e$ and $x_i \neq a$, for i = 1, 2, …, k)

Now $x_1 • a = x_2$ that is, $x_2 \equiv (ma − (m − 1)x_1) \pmod{n}$
$x_2 • a = x_3$ so, $x_3 \equiv (ma − (m − 1)x_2) \pmod{n}$ or
$x_3 \equiv (ma − (m − 1)(ma − (m − 1) x_1)) \pmod{n}$ or $x_3 \equiv (ma −(m−1)a + (m- 1)^2 x_1) \pmod{n}$.

Finally $x_k • a = x_1$ so we get ma $(1 − (m − 1) + (m − 1)^2 + \ldots + (-1)^{k-1} (m-1)^{k-1}) + (-1)^k (m − 1)^k x_1 \equiv x_1 \pmod{n}$. Simplifying this we get $((m − 1)^k + (-1)^{k−1})(a − x_1) \equiv 0 \pmod{n}$ Since any $x_i$'s for 1 ≤ i ≤ k can be taken as the first entry in the cycle $(x_1, x_2, \ldots, x_k)$ we have $((m − 1)^k + (−1)^{k−1}) (a − x) \equiv 0 \pmod{n}$ for any x ∈ $\{x_1, x_2, \ldots, x_k\}$.

Hence the result.

The converse of the above theorem is not true in general. For consider the loop $L_7(4)$ given by the following table:

*Example 2.7.1*:

| • | e | 1 | 2 | 3 | 4 | 5 | 6 | 7 |
|---|---|---|---|---|---|---|---|---|
| e | e | 1 | 2 | 3 | 4 | 5 | 6 | 7 |
| 1 | 1 | e | 5 | 2 | 6 | 3 | 7 | 4 |
| 2 | 2 | 5 | e | 6 | 3 | 7 | 4 | 1 |
| 3 | 3 | 2 | 6 | e | 7 | 4 | 1 | 5 |
| 4 | 4 | 6 | 3 | 7 | e | 1 | 5 | 2 |
| 5 | 5 | 3 | 7 | 4 | 1 | e | 2 | 6 |
| 6 | 6 | 7 | 4 | 1 | 5 | 2 | e | 3 |
| 7 | 7 | 4 | 1 | 5 | 2 | 6 | 3 | e |

The right regular representation is given by

I

(e 1) (2 5 3) (4 6 7)
(e 2) (1 5 7) (3 6 4)
(e 3) (1 2 6) (4 7 5)
(e 4) (1 6 5) (2 3 7)
(e 5) (1 3 4) (2 7 6)
(e 6) (1 7 3) (2 4 5)
(e 7) (1 4 2) (3 5 6)

where I is the identity permutation of the loop $L_7(4)$.

Here m = 4 so $(4-1)^6 + (-1)^{6-1} = 728 \equiv 0 \pmod{7}$ so $(4-1)^6 + (-1)^{6-1}(a-x) \equiv 0 \pmod{7}$ for any $a, x \in L_7(4)$, but in the representation of $L_7(4)$, no permutation has a cycle of length 6.

As an immediate consequence of this we have the following:

**COROLLARY [65]**: *Let $L_n(m) \in L_n$. If m, k (k < n) such that $((m-1)^k + (-1)^{k-1}, n) = 1$, then any permutation $R_a \in \pi$ ($\pi$ the representation of $L_n(m)$) does not have a cycle of length k.*

*Proof*: Suppose the permutation $R_a \in \pi$ has a cycle of length k, then $((m-1)^k + (-1)^{k-1})(a-x) \equiv 0 \pmod{n}$ for any x in the k-cycle. But this cannot happen as $((m-1)^k + (-1)^{k-1}, n) = 1$ and $(|x-a|, n) < n$. So the permutation $R_a$ does not contain a k-cycle.

**COROLLARY [65]**: *Let $L_n(m) \in L_n$. If there exists $k \in N$ (N-natural numbers, k < n), such that $((m-1)^k + (-1)^{k-1}) \equiv 0 \pmod{n}$ then there exists no cycle of length greater than k in any permutation $\alpha \in \pi$ where $\pi$ is the representation of $L_n(m)$.*

*Proof*: $((m - 1)^k + (-1)^{k-1}) \equiv 0 \pmod{n}$ implies $((m - 1)^k + (-1)^{k-1})(a - x) \equiv 0$ mod n for any $a, x \in L_n(m)$. So all the elements of $L_n(m)$ will decompose into the cycles of length at most k, hence no cycle of length greater than k exists.

**THEOREM [65]**: *Let $L_n(m) \in L_n$. If n is an odd prime and k is the least positive integer such that $(m - 1)^k \equiv (-1)^k \pmod{n}$, then any permutation $R_a \in \pi (a \neq e)$ in the representation of $L_n(m)$ is a product of 2 cycles and t, k-cycles where $t = (n - 1) / k$.*

*Proof*: We know by earlier result and by the construction of $L_n(m)$ there exists a 2 cycle namely (a, e). Now we will show that the other $p - 1$ elements will decompose only as k-cycles. Since it is given that $((m - 1)^k + (-1)^{k-1}) \equiv 0 \pmod{n}$ by the earlier corollaries there exists no cycle of length greater than k. Now it only remains to show that there exists no cycle of length less than k. If not, let there exist a cycle of length s ($s \leq k$) in the permutation $R_a$ ($a \in L_n(m)$ and $a \neq e$). Then we have by earlier theorem $((m - 1)^s + (-1)^{s-1})(a - x) \equiv 0 \pmod{n}$. Now as $x \neq a$ and n is a prime number we have $(m - 1)^s \equiv (-1)^s \pmod{n}$, which is a contradiction as k is the least such number. So no cycle of length less than k-exists. Thus the remaining $(n - 1)$ elements will decompose into $(n - 1) / k$ cycles of length k.

Using the theorem and corollaries already proved we find out the decomposition of $L_n(m)$ when n is a composite number.

**THEOREM [65]**: *Let $L_n(m) \in L_n$ n a composite number. If $((m - 1)^k + (-1)^{k-1}, n) = d_k$ ($d_k < n$) and k is the least integer for this $d_k$, then in any permutation $R_a \in \pi (a \neq e$ and $\pi$ the representation of $L_n(m))$, there exists $[(d_k - 1) / k]$ cycles of length k, where [x] denotes the largest integer not exceeding x for any $x \in R$ (R: the set of real numbers) and the remaining elements will decompose into cycles of length t where t is the least integer such that $((m - 1)^t + (-1)^{t-1}), n) = n$.*

We will illustrate this result by an example.

***Example 2.7.2***: Let $L_{45}(8) \in L_{45}$. We will find out the cycle decomposition in the representation of the loop $L_{45}(8)$. Here $m = 8$ and $k = 2$ is the least integer such that $(((m - 1)^2 + (-1)^{2-1}), 45) = 3$, so there exists $[(3 - 1) / 2] = 1$ cycle of length 2. $k = 4$ is the least integer such that $(((m - 1)^4 + (-1)^{4-1}), 45) = 15$, so that there exist $[(15-1)/4] = 3$ cycles of length 4.

Further $k = 6$ is the least number such that $(((m - 1)^6 + (-1)^{6-1}), 45) = 9$ hence there exists $[(9 - 1) / 8] =$ cycles of length 6. Now $k = 12$ is the least number such that $(((m - 1)^{12} + (-1)^{12-1}), 45) = 45$. So the remaining 24 element will decompose into cycles of length 12, hence there exist 2 cycles of length 12.

**COROLLARY [65]**: *Let $L_n(m) \in L_n$. Then each permutation $R_a \in \pi$ ($a \neq e$ and $\pi$ is the representation of $L(m)$) belongs to the same cycle class.*

*Proof*: As in theorems and results we see the cycle decomposition does not depend on the choice of $a \in L_n(m)$, that is it is the same for all a, $a \neq e$. So each permutation $R_a \in \pi$ belongs to the same cycle class. We know that if the cycle class of permutation is $(k_1, k_2, \ldots, k_n)$, then its order is the l.c.m. of the non-zero $k_i$'s.

Thus we have the following result:

**COROLLARY [65]**: *Order of each permutation $R_a$ in $\pi$ is the same for all $a \in L_n(m)$ ($a \neq e$).*

**LEMMA [65]**: *Let L be a right alternative loop of even order say 2s ($s \geq 3$) in which the square of each element is the identity. Then its representation contains identity permutation on L and (2s – 1) other permutations which are products of disjoint transpositions. Further no two permutations have any transpositions in common.*

*Proof*: Let L be a right alternative loop of even order in which square of each element is identity and $\pi$ be its representation. Clearly $I \in \pi$ (I is the identity permutation on L). Take $a \in L$ ($a \neq e$). Let $R_a: x \to y \Rightarrow xa = y$. Now $ya = (xa) a = x(aa)$ (as L is right alternative) = x (as the square of each element is identity) so if $R_a: x \to y$ then $R_a: y \to x$ that is $R_a$ is the product of disjoint transpositions. If the transposition (x, y) belongs to both $R_a$ and $R_b$ for some a, b $\in$ L, then we have $R_a: x \to y$ and $R_b: x \to y$ that is xa = y and xb = y so xa = ab which gives a = b.

Thus no two different permutations have any transposition in common.

**LEMMA [65]**: *Let L = {1, 2, ..., n + 1} n odd and $S_{n+1}$ be the set of all permutations on the set L. If $\pi \in S_{n+1}$ is such that*

> i. *$\pi$ contains the identity permutation on L.*
> ii. *$\pi$ contains n non-identity permutations satisfying the following two conditions*
>> a. *Each permutation is the product of disjoint transpositions*
>> b. *No two permutations have any transposition in common.*
>> *Then $\pi$ is a representation of a right alternative loop of order n + 1 in which square of each element is identity.*

*Proof*: Let $\pi \in S_{n+1}$ be a set of permutations on L satisfying the above conditions. First we will show that $\pi$ is a representation of some loop on the set L. To prove the 3 results of Albert [1].

    a. Clearly $\pi$ has identity.

    b. To show that $\pi$ is transitive on L. Now these n non-identity distinct permutations have $(n(n+1))/2$ different transpositions. But this is the number of possible transpositions over the set of $(n+1)$ elements. So for any a, b $\in$ L, there exists a $\alpha \in \pi$ such that (a, b) $\in \alpha$, $\alpha: a \to b$. Hence $\pi$ is transitive on L.

    c. We will show that $\alpha\beta^{-1}$ ($\alpha \neq \beta$) where $\alpha, \beta \in \pi$, does not fix any element of L. For if x $\in$ L then the transposition (x, y) $\in \beta$ for some y $\in$ L. For this y $\in$ L, the transposition (y, z) $\in \alpha$ for some z $\in$ L. Hence z $\neq$ x otherwise both $\alpha$ and $\beta$ will have (x, y) which is not possible as $\alpha \neq \beta$. Hence $\pi$ is a representation for some loop on the set L. Now we will show that

        i.    $x^2 = e$ (where e is the identity of the loop L) for all x $\in$ L and
        ii.   L is a right alternative loop.

To prove $x^2 = e$ (x $\neq$ e). Let x $\in$ L, consider $R_x \in \pi$ if $R_x: x \to a$ then $R_x: a \to x$ for some a $\in$ L, that is if x • x = a and a • x = x. Thus (x • x) x = x implies x • x = e or $x^2$ = e. To prove L is right alternative let x, y $\in$ L consider $R_y \in \pi$ if $R_y : x \to a$ then $R_y: a \to x$ that is x • y = a and a • y = a. Therefore (xy) y = x as (a = xy) or (xy)y = x • e = x(yy). Hence L is right alternative. Hence $\pi$ represents a loop on L, which is right alternative and in which square of every element is identity.

**THEOREM [65]**: *Let $\pi \in S_{n+1}$, satisfy the conditions (i) and (ii) of the above lemma. If f(n) is the number of distinct possible choices of $\pi$, then any right alternative loop of order n + 1 (n-odd) in which square of each element is identity has the representation out of these f(n) representations.*

*Proof*: Clearly from the above results, the theorem is true.

The major application to edge colouring of the graph $K_{2n}$ is given below.

**THEOREM [65]**: *The number of different representations of loops of order 2n (n > 3) which are right alternative and in which square of each element is* *identity, is equal to different proper edge colourings of the graph $k_{2n}$ using exactly $2n - 1$ colours.*

*Proof*: Let $L = \{1, 2, \ldots, 2n\}$ be a set of $2n$ elements and $\pi \subset S_{2n}$ be the set of permutations on L as in the above theorem. Now construction of $\pi$ can be put in other words as follows. Let the number of objects be $2n$ say $1, 2, \ldots, 2n$. From these $2n$ numbers we can form $(2n)(2n-1)/2$ (this is $n(2n-1)$) distinct non-ordered pairs in $(2n - 1)$ lines such that each line contains n pairs in which all elements from 1 to $2n$ occurs. Obviously, no two lines will have any pair in common. Such partition of pairs together with identity will give $\pi$.

Now let us look at this problem graph theoretically. Consider the graph $K_{2n}$ which has $2n$ vertices and $2n(2n-1)/2$ edges that is $n(2n - 1)$ edges. So we can treat each non-ordered pair as an edge. Now if we give different $(2n - 1)$ colours to the pairs of different $(2n - 1)$ lines, then no two edges which meet at a vertex will have the same colour, which is a proper edge colouring of the graph $K_{2n}$ using $(2n - 1)$ colours. Hence the number of different representations of our loops and the number of different proper edge colourings of the graph $K_{2n}$ using $(2n - 1)$ colours are the same.

*Example 2.7.3*: There exists only six different representations of loops of order 6, which are right alternative and in which square of each element is identity

I
(e 1) (2 3) (4 5)
(e 2) (1 5) (3 4)
(e 3) (1 4) (2 5)
(e 4) (1 2) (3 5)
(e 5) (1 3) (2 4)

I
(e 1) (2 3) (4 5)
(e 2) (1 4) (3 5)
(e 3) (1 5) (2 4)
(e 4) (1 3) (2 5)
(e 5) (1 2) (3 4)

I
(e 1) (2 5) (3 4)
(e 2) (1 3) (4 5)
(e 3) (1 5) (2 4)
(e 4) (1 2) (3 5)
(e 5) (1 4) (2 3)

I

(e 1) (2 4) (3 5)
(e 2) (1 5) (3 4)
(e 3) (1 2) (4 5)
(e 4) (1 3) (2 5)
(e 5) (1 4) (2 3)

I
(e 1) (2 4) (3 5)
(e 2) (1 3) (4 5)
(e 3) (1 4) (2 5)
(e 4) (1 5) (2 3)
(e 5) (1 2) (3 4)

I
(e 1) (2 5) (3 4)
(e 2) (1 4) (3 5)
(e 3) (1 2) (4 5)
(e 4) (1 5) (2 3)
(e 5) (1 3) (2 4)

Now we suggest some problems as exercise for the reader.

**PROBLEMS**:

1. Prove using representation theory $L_5(2)$ is not a Bol loop.
2. Can you prove in $L_n(m) \in L_n$ if $m \neq 2$ then $L_n(m)$ cannot be a Bol loop?
3. Show that the loop $L_{13}(7)$ is such that the number of different representation of the loop $L_{13}(7)$ and the number of different proper edge colourings of the graph $k_{14}$ using 13 colours are the same.
4. Can any loop in $L_{27}$ be a Bol loop?

**Chapter three**

# SMARANDACHE LOOPS

In this chapter we introduce the notion of Smarandache loops and illustrate them by examples. This chapter has ten sections. The first section defines Smarandache loop. In section 2 we introduce substructures in loops like S-subloops, S-subgroups etc. Section 3 is completely spent to introduce new classical Smarandache loops. We introduce commutator subloops and Smarandache commutative loops in section 4. Section 5 is used to define Smarandache associativity and Smarandache associator subloops. Since study of loops, which satisfy identity, is a major study in loops in section 6 we introduce Smarandache identities in loops. The next section studies the special structure like S-left nuclei, S-middle nuclei, S-Moufang centre etc.

The concept of Smarandache mixed direct product is introduced in section 8 which enables us to define Smarandache loop level II using this only we are in a position to extend the classical theorems, Smarandache Lagrange's criteria and Smarandache Sylow criteria. Further in section 9 we introduce Smarandache cosets in loops. The last section is used to define the Smarandache analogue of the very recently introduced loops by Michael K. Kinyon [40] in 2002, Smarandache ARIF loops, Smarandache RIF loops, Smarandache Steiner loops. Only these concepts has helped us to get odd order Moufang loops which are S-loops.

Each section of this chapter gives a brief introduction of what is carried out in that particular section.

## 3.1 Definition of Smarandache Loops with examples

This section is devoted to the introduction of Smarandache loops and illustrate them with examples. Smarandache loops were introduced very recently (2002) [81]. Here we also give which of the well-known classes of loops are Smarandache loops. We would call this Smarandache loop of level I but we don't add the adjective level I but just call it as Smarandache loops as Smarandache loops of level II will be distinguished by putting II beside the definition.

**DEFINITION 3.1.1**: *The Smarandache loop (S-loop) is defined to be a loop L such that a proper subset A of L is a subgroup (with respect to the same induced operation) that is $\phi \neq A \subset L$.*

**Example 3.1.1**: Let L be the loop given the following table:

| • | e | $a_1$ | $a_2$ | $a_3$ | $a_4$ | $a_5$ | $a_6$ | $a_7$ |
|---|---|---|---|---|---|---|---|---|
| e | e | $a_1$ | $a_2$ | $a_3$ | $a_4$ | $a_5$ | $a_6$ | $a_7$ |
| $a_1$ | $a_1$ | e | $a_4$ | $a_7$ | $a_3$ | $a_6$ | $a_2$ | $a_5$ |
| $a_2$ | $a_2$ | $a_6$ | e | $a_5$ | $a_1$ | $a_4$ | $a_7$ | $a_3$ |
| $a_3$ | $a_3$ | $a_4$ | $a_7$ | e | $a_6$ | $a_2$ | $a_5$ | $a_1$ |
| $a_4$ | $a_4$ | $a_2$ | $a_5$ | $a_1$ | e | $a_7$ | $a_3$ | $a_6$ |
| $a_5$ | $a_5$ | $a_7$ | $a_3$ | $a_6$ | $a_2$ | e | $a_1$ | $a_4$ |
| $a_6$ | $a_6$ | $a_5$ | $a_1$ | $a_4$ | $a_7$ | $a_3$ | e | $a_2$ |
| $a_7$ | $a_7$ | $a_3$ | $a_6$ | $a_2$ | $a_5$ | $a_1$ | $a_4$ | e |

This is a S-loop as $H_i = \{e, a_i\}$; i = 1, 2, 3, …, 7 are subsets of L which are subgroups.

We have the following classes of loops to be S-loops.

**THEOREM 3.1.1**: *The natural class of loops $L_n(m) \in L_n$ (n odd, n > 3, (m, n) = 1, (m – 1, n) = 1 for varying m) are S-loops.*

*Proof*: We see by the very construction of loops $L_n(m)$ in $L_n$ each i ∈ $L_n(m)$ is such that i • i = e where e is the identity element of $L_n(m)$. Thus all proper subsets of the form {e, i} ⊂ $L_n(m)$ for varying i are groups. Thus the class of loops $L_n$ are S-loops.

It is very exciting to note that the class of loops $L_n$ are not Moufang or Bol or Bruck. To be more precise that none of the loops in the class $L_n$ are Moufang or Bol or Bruck. But every loop in the class $L_n$ are S-loops.

**THEOREM 3.1.2**: *Every power associative loop is a S-loop.*

*Proof*: By the very definition a power associative loop L, we see every element in L generates a subgroup in L. Hence the claim.

**THEOREM 3.1.3**: *Every diassociative loop is a S-loop.*

*Proof*: Since a loop L is diassociative if every pair of elements of L generate a subgroup in L. Thus every diassociative loop is a S-loop.

**THEOREM 3.1.4**: *Let L be a Moufang loop which is centrally nilpotent of class 2. Then L is a S-loop.*

*Proof*: We know if L is a Moufang loop which is centrally nilpotent of class 2, that is, a Moufang loop L such that the quotient of L by its centre Z(L) is an abelian group; and let $L_p$ denote the set of all elements of L whose order is a power of p. Recall that the nuclearly derived subloop, or normal associator subloop of L, which we denote by $L^*$ is the smallest normal subloop of L such that $L/L^*$ is associative (i.e. a group). Recall also that the torsion subloop (subloop of finite order elements) of L is isomorphic to the (restricted) direct product of the subloops $L_p$ where p runs over all primes.

(Theorem 6.2. of R. H. Bruck [6] and [7] some theorems on Moufang loops and Theorem 3.9 of Tim Hsu [63] or in the finite case, Corollary 1 of Glauberman and Wright [22]). Using the main theorem A of Tim Hsu [63] we see that for p > 3, $L_p$ is a group. So L is S-loop when L is a Moufang loop which is centrally nilpotent of class 2.

In view of the Theorem 2.11 of Tim Hsu [63] we have the following:

**THEOREM 3.1.5**: *Let L be a Moufang loop. Then L is a S-loop.*

*Proof*: Now according to the theorem 2.11 of Tim Hsu [63] we see if L is a Moufang loop then L is diassociative that is for (x,y) ∈ L, ⟨x, y⟩ is associative. By theorem 3.1.3. of ours a Moufang loop is a S-loop.

Now we know Bol loops are power-associative (H.O. Pflugfelder [49, 50] and D. A. Robinson [55]) leading to the following theorem.

**THEOREM 3.1.6**: *Every Bol loop is a S-loop.*

*Proof*: Bol loops are power-associative from the above references and by our theorem 3.1.2. All Bol loops are S-loops.

Now Michael Kinyon [42] (2000) has proved every diassociative A-loop is a Moufang loop in view of this we have.

**THEOREM 3.1.7**: *Let L be a A-loop. If L is diassociative or L is a Moufang loop then L is a S-loop.*

*Proof*: In view of Corollary 1 of Kinyon [42] we see for an A-loop L. The two concepts; L is diassociative is equivalent to L is a Moufang loop so we have L to be S-loop.

In view of this we propose an open problem in chapter 5.

From the results of Kinyon [41] 2001 we see that every ARIF loop is a diassociative loop in this regard we have the following.

**THEOREM 3.1.8**: *Every ARIF loop is a S-loop.*

*Proof*: Kinyon [41] has proved in Theorem 1.5 that every ARIF loop is diassociative from our theorem 3.1.3 every diassociative loop is a S-loop. So all ARIF loops are S-loops.

Further using the Lemma 2.4 of Kinyon [41] we have another interesting result.

**THEOREM 3.1.9**: *Every RIF-loop is a S-loop.*

*Proof*: According to Kinyon's [41] Lemma 2.4. every RIF loop is an ARIF loop. But by our Theorem 3.1.8 every ARIF loop is a S-loop. Hence every RIF loop is a S-loop.

Fenyves [20] (1969) has introduced the concept of C-loops.

By Corollary 2.6 of Kinyon [41], every flexible C-loop is a ARIF-loop, we have a theorem.

**THEOREM 3.1.10**: *Every flexible C-loop is a S-loop.*

*Proof*: We see every flexible C-loop is a ARIF loop by Corollary 2.6. of Kinyon [41], Hence by Theorem 3.1.8 every flexible C-loop is a S-loop.

This theorem immediately does not imply that all non-flexible C-loops are not S-loops; for we have the example 4.2 given in Kinyon [41] is a C-loop which is not flexible given by the following table:

| •  | 0  | 1  | 2  | 3  | 4  | 5  | 6  | 7  | 8  | 9  | 10 | 11 |
|----|----|----|----|----|----|----|----|----|----|----|----|----|
| 0  | 0  | 1  | 2  | 3  | 4  | 5  | 6  | 7  | 8  | 9  | 10 | 11 |
| 1  | 1  | 2  | 0  | 4  | 5  | 3  | 7  | 8  | 6  | 10 | 11 | 9  |
| 2  | 2  | 0  | 1  | 5  | 3  | 4  | 8  | 6  | 7  | 11 | 9  | 10 |
| 3  | 3  | 4  | 5  | 0  | 1  | 2  | 10 | 11 | 3  | 8  | 6  | 7  |
| 4  | 4  | 5  | 3  | 1  | 2  | 0  | 11 | 9  | 10 | 6  | 7  | 8  |
| 5  | 5  | 3  | 4  | 2  | 0  | 1  | 9  | 10 | 11 | 7  | 8  | 8  |
| 6  | 6  | 7  | 8  | 11 | 9  | 10 | 0  | 1  | 2  | 4  | 5  | 3  |
| 7  | 7  | 8  | 6  | 9  | 10 | 11 | 1  | 2  | 0  | 5  | 3  | 4  |
| 8  | 8  | 6  | 7  | 10 | 11 | 9  | 2  | 0  | 1  | 3  | 4  | 5  |
| 9  | 9  | 10 | 11 | 7  | 8  | 6  | 5  | 3  | 4  | 0  | 1  | 2  |
| 10 | 10 | 11 | 9  | 8  | 6  | 7  | 3  | 4  | 5  | 1  | 2  | 0  |
| 11 | 11 | 9  | 10 | 6  | 7  | 8  | 4  | 5  | 3  | 2  | 0  | 1  |

This C-loop is a S-loop as the set {0, 1, 2, •} is a group given by the following table:

| • | 0 | 1 | 2 |
|---|---|---|---|
| 0 | 0 | 1 | 2 |
| 1 | 1 | 2 | 0 |
| 2 | 2 | 0 | 1 |

Thus we have a non-flexible C-loop which is a S-loop leading to the following open problem proposed in Chapter 5. For more about C-loops please refer (F. Fenyves) [20].

Thus we have seen the definition of S-loops and the classes of loops which are S-loops. Now we proceed on to define substructures in S-loops.

**PROBLEMS:**

1. Give an example of a loop of order 8 which is not a S-loop.
2. Can a Bruck loop be a S-loop?
3. Find a S-loop of order 9 which is not a Moufang loop.
4. What is the order of the smallest S-loop?
5. Does there exist a S-loop of order p, p a prime?

## 3.2 Smarandache substructures in Loops

In this section we introduce the concept of Smarandache subloops (S-subloops) and Smarandache normal subloops (S-normal subloops). The absence of S-normal subloops in a loop L leads us to the definition of Smarandache simple loops (S-simple loops). We prove that the class of loops $L_n$ is a S-simple loop of order n + 1. We obtain some results about these definitions.

**DEFINITION 3.2.1**: *Let L be a loop. A proper subset A of L is said to be a Smarandache subloop (S-subloop) of L if A is a subloop of L and A is itself a S-loop; that is A contains a proper subset B contained in A such that B is a group under the operations of L. We demand A to be a S-subloop which is not a subgroup.*

In view of this we have the following:

**THEOREM 3.2.1**: *Let L be a loop. If L has a S-subloop then L is a S-loop.*

*Proof*: If a loop L has S-subloop then we have a subset $A \subset L$ such that A is a subloop and contains a proper subset B such that B is a group. Hence $B \subset A \subset L$ so L is a S-loop. So a subloop can have a S-subloop only when L is a S-loop.

*Example 3.2.1*: Consider the loop $L_{15}(2) \in L_{15}$. H = {e, 1, 4, 7, 10, 13} is a subloop of the loop $L_{15}(2)$. Clearly H is a S-subloop of L. But it is interesting to note that in general all S-loops need not have every subloop to be a S-subloop or more particularly a S-loop need not have S-subloops at all.

This following example will show it.

*Example 3.2.2*: Let $L_5(2)$ be a loop given by the following table:

| • | e | 1 | 2 | 3 | 4 | 5 |
|---|---|---|---|---|---|---|
| e | e | 1 | 2 | 3 | 4 | 5 |
| 1 | 1 | e | 3 | 5 | 2 | 4 |
| 2 | 2 | 5 | e | 4 | 1 | 3 |
| 3 | 3 | 4 | 1 | e | 5 | 2 |
| 4 | 4 | 3 | 5 | 2 | e | 1 |
| 5 | 5 | 2 | 4 | 1 | 3 | e |

It can be easily verified that this loop is a S-loop. This has in fact no subloops but only 5 subgroups each of order two given by {e, 1}, {e, 2}, {e, 3}, {e, 4} and {e, 5}.

In view of this we have the following definition:

**DEFINITION 3.2.2**: *Let L be a S-loop. If L has no subloops but only subgroups we call L a Smarandache subgroup (S-subgroup) loop.*

We have very many loops of even order given by the following.

**THEOREM 3.2.2**: *Let $L_n(m) \in L_n$ where n is a prime. Then the class of loops $L_n$ is a S-subgroup loop.*

*Proof*: Given n is a prime. So $L_n(m) \in L_n$ has n + 1 elements and further no number t divides n. By the very construction of $L_n(m)$ we see $L_n(m)$ is a S-loop every element generates a cyclic group of order 2. Thus we have a class of loops $L_n$ which are S-subgroup loops for each prime n = p, n > 3.

From this we see the existence of S-subgroup loops.

Now we proceed onto define Smarandache normal subloops.

**DEFINITION 3.2.3**: *Let L be a loop. We say a non-empty subset A of L is a Smarandache normal subloop (S-normal subloop) of L if*

1. *A is itself a normal subloop of L.*
2. *A contains a proper subset B where B is a subgroup under the operations of L. If L has no S-normal subloop we say the loop L is Smarandache simple (S-simple).*

We have the following interesting results:

**THEOREM 3.2.3**: *Let L be a loop. If L has a S-normal subloop then L is a S-loop.*

*Proof*: Obvious by the very definition of S-normal subloop we see L is a S-loop.

Now we see that a loop may have normal subloops but yet that normal subloop may not be a S-normal subloop.

In view of this we have the following:

**THEOREM 3.2.4**: *Let L be a loop. If L has a S-normal subloop then L has a normal subloop, so L is not simple.*

*Proof*: By the very definition of S-normal subloop in a loop L we are guaranteed, that the loop L must have a normal subloop so L is not simple. Hence the claim.

We define analogous to the case of groups where the normal subgroups of a group form a modular lattice in the following:

**DEFINITION 3.2.4**: *Let L be a S-loop. The lattice of all S-subloops of L is denoted by S(L). We call this representation by Smarandache lattice representation of S-subloops. Similarly we define Smarandache lattice representation of S-normal subloops and subgroups of a loop L.*

The study about these structures in S-loops unlike in groups or loops remain at a very dormant state. We prove the following theorems only in case of the new class of loops $L_n$ which are analogous to results in groups.

**THEOREM 3.2.5**: *Let $L_n(m) \in L_n$ be the class of S-loops. The Smarandache lattice representation of S-normal subloops of the loop $L_n(m)$ form a two element chain lattice.*

*Proof*: Every $L_n(m) \in L_n$ has no S-normal subloops. So the only trivial S-normal subloops are e and $L_n(m)$ giving the two element chain

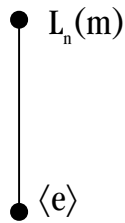

This can be compared with the normal subgroups in the alternating group $A_n$, $n \geq 5$.

**THEOREM 3.2.6**: *Let $L_n(m) \in L_n$, n an odd prime. S(L) the lattice of subgroups of every loop $L_n(m)$ in $L_n$ has a modular lattice with n + 2 elements for the Smarandache lattice representation of S-subgroups.*

*Proof*: Now $L_n(m)$ when n is a prime has only n subgroups each of order two. Say $A_i$ = {i, e} for i = 1, 2, …, n having the following lattice representation:

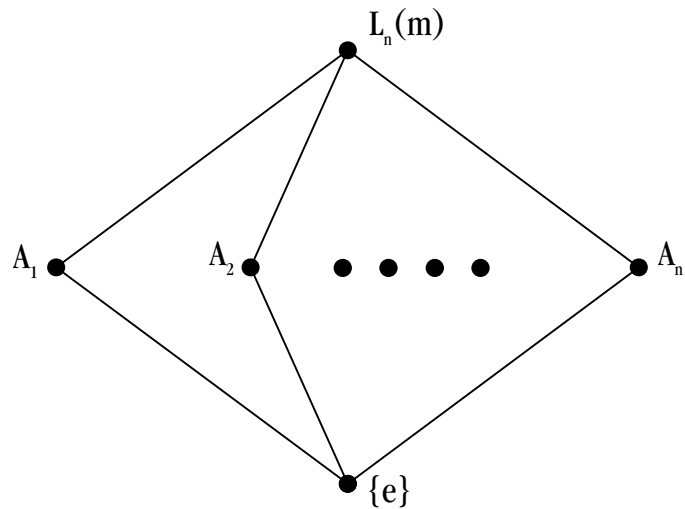

We are not able to state whether in every loop L the normal subloop of L is a S-normal subloop of L. This problem we leave it as a problem to be solved in chapter 5. We have a class of loops which are S-simple.

**THEOREM 3.2.7**: *Let $L_n(m) \in L_n$, n a composite number. The S-subloops of $L_n(m)$ has a modular lattice representation with t-subloops where n/t = p, p and t are primes, so p-subloops as n/p = t (p > 3 and t > 3).*

*Proof*: We know $L_n(m) \in L_n$ where n = pt where both p and t are primes. We have p + t S-subloops which form a S-modular lattice.

We have $A_1$, $A_2$, …, $A_p$ and $B_1$, $B_2$, …, $B_t$ as S-subloops having the following lattice representation:

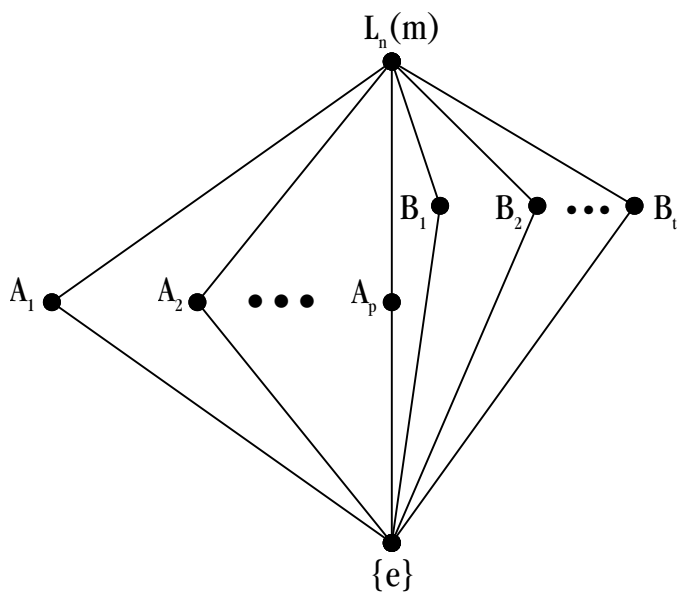

For n not a product of two primes the question of representation of S-subloops remains an open problem in case of loops in $L_n$.

***Example 3.2.3***: Let $L_{15}(8)$ be a loop in $L_{15}$. The subgroups of $L_{15}(8)$ are $B_1$ = {e, 1, 6, 11}, $B_2$ = {e, 2, 7, 12}, $B_3$ = {e, 3, 8, 13}, $B_4$ = {e, 4, 9, 14}, $B_5$ = {e, 5, 10, 15} which are of order 4 and $A_i$ = {e, i} for i = 1, 2, …, 15 each of order 2. Thus we have 22 subgroups which has the following S-lattice representation.

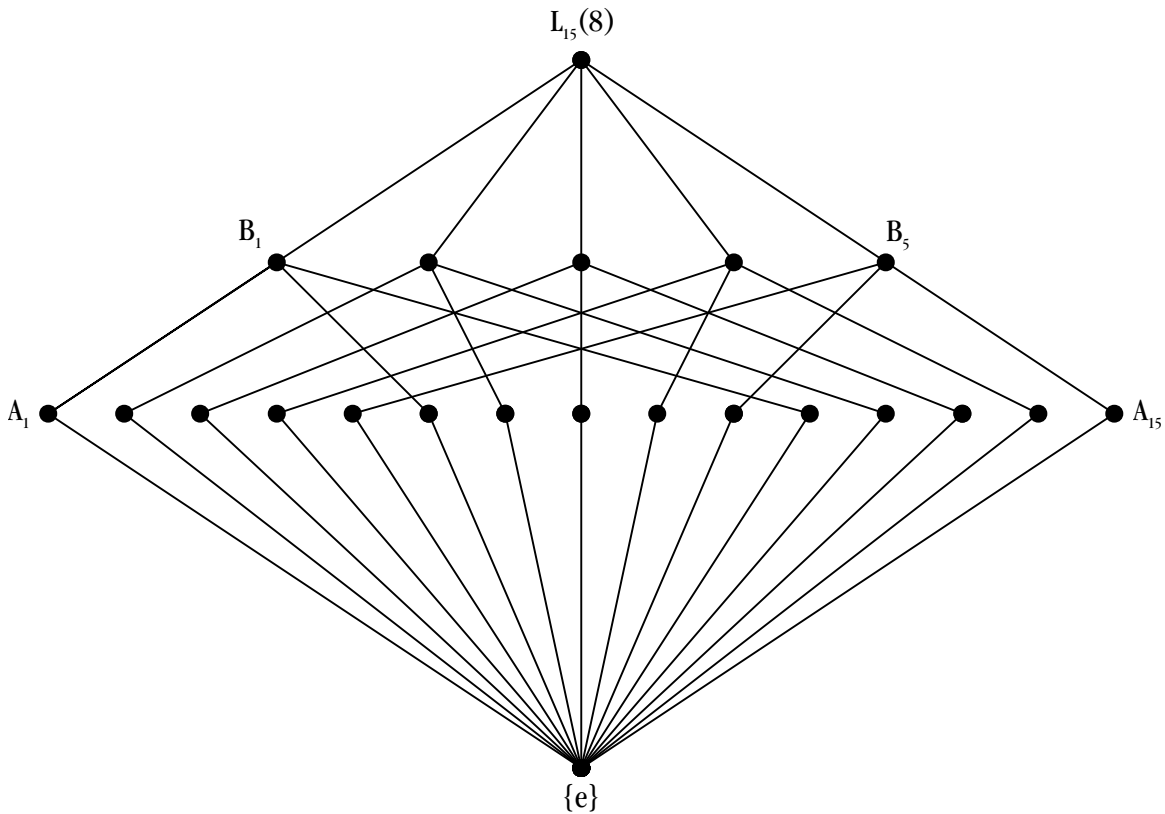

This lattice is both non-modular and non-distributive. This set ($L_{15}(8)$, $A_1'$ = {e, 11, 4, 7, 10, 13}, $A_2'$ = {e, 2, 5, 8, 11, 14}, $A_3'$ = {e, 3, 6, 9, 12, 15}, {e}) forms the diamond lattice which is modular.

**THEOREM 3.2.8**: *Let $L_n(m) \in L_n$. Then $L_n(m)$ are S-simple.*

*Proof*: Each $L_n(m) \in L_n$ is simple, so each $L_n(m)$ has no non-trivial normal subloops so each $L_n(m) \in L_n$ are S-simple. Thus we have a class of loops which are S-simple.

<u>PROBLEMS:</u>

1. Find all the S-subloops of $L_{17}(2)$.
2. Can the loop $L_9(5)$ have S-normal subloop? Justify.
3. Prove $L_{19}(4)$ is not a S-simple loop.

4. Is $L_{23}(4)$ a S-subgroup loop?
5. How many subloops does $L_{19}(4)$ have? Is every subloop in $L_{19}(4)$ a S-subloop?
6. Give an example of S-simple loop of order 29.
7. Can a loop of order 23 have S-normal subloops?
8. Give an example of a loop L in which every normal subloop is a S-normal subloop.
9. Give an example of a loop of odd order in which every subloop is a subgroup.
10. Can a Bruck loop be a S-simple loop?
11. Give an example of a Bruck loop which is not S-simple.
12. Give an example of a non-simple loop which is a S-simple loop.
13. Give an example of a Moufang loop which is S-subgroup loop.
14. Can a Moufang loop be S-simple?
15. Can a Bol loop be S-simple?
16. Can a Bol loop be a S-subgroup loop?
17. Find a loop L of order 17 which has no S-subloops.

## 3.3 Some new classical S-loops

In this section we define some new classical S-loops. When we say new they are loops not contemplated till date, the term classical means that these loops are built or makes use of classical properties in groups like Cauchy, Sylow or Lagranges so we choose to call this section 'Some new classical S-loops'. We define and give examples of them.

**DEFINITION 3.3.1**: *Let L be a finite S-loop. An element $a \in A$, $A \subset L$; A the subgroup of L is said to be Smarandache Cauchy element (S-Cauchy element) of L if $a^r = 1$ ( $r > 1$), 1 the identity element of L and r divides order of L; otherwise a is not a Smarandache Cauchy element of L.*

**DEFINITION 3.3.2**: *Let L be a finite S-loop. If every element in every subgroup is a S-Cauchy element then we say S is a Smarandache Cauchy loop (S-Cauchy loop).*

**Example 3.3.1**: Consider the loop $L_{11}(3)$. Clearly $L_{11}(3)$ is a S-Cauchy loop as every subgroup of $L_{11}(3)$ is a cyclic group of order two and $2/|L_{11}(3)|$ as the order of $L_{11}(3)$ is 12.

In view of this we have a nice theorem.

**THEOREM 3.3.1**: *Every loop in the class of loops $L_n$ are S-Cauchy loops.*

*Proof*: Every loop in the class of loops $L_n$ is of even order. Further by the very construction of the loops in $L_n(m) \in L_n$, every element in $L_n(m)$ is of order two. Since all loops in $L_n$ are S-loops and $2/|L_n(m)|$ we see every loop in $L_n$ are S-Cauchy loops.

**DEFINITION 3.3.3**: *Let L be a finite loop. If the order of every subgroup in L divides the order of L we say L is a Smarandache Lagrange loop (S-Lagrange loop).*

**DEFINITION 3.3.4**: *Let L be a finite loop. If there exists at least one subgroup in L whose order divides order of L then L is a Smarandache weakly Lagrange (S-weakly Lagrange) loop.*

**THEOREM 3.3.2**: *Every S-Lagrange loop is a S-weakly Lagrange loop.*

*Proof*: By the very definition of S-Lagrange loop and S-weakly Lagrange loop the result follows.

In view of this we get a class of S-Lagrange loops of even order p + 1 where p is a prime greater than three.

**THEOREM 3.3.3**: *Let $L_p$ be the class of loops, for any prime p. Then every loop in the class of loops $L_p$ are S-Lagrange loop.*

*Proof*: We know $L_p(m) \in L_p$ for every prime p is a loop of order p + 1. The only subgroups in $L_p(m)$ are of order 2. Hence the claim.

**THEOREM 3.3.4**: *Let $L_n(m) \in L_n$. Every loop in $L_n$ is a S-weakly Lagrange loop.*

*Proof*: To show every $L_n(m)$ in $L_n$ is a S-weakly Lagrange loop we have to show $L_n(m)$ has at least one subgroup whose order divides the order of $L_n(m)$. Clearly for any $i \in L_n(m)$ we have (e, i) is a subgroup of order 2 as all loops in $L_n$ are of even order $2/|L_n(m)|$. Hence every $L_n(m) \in L_n$ are S-weakly Lagrange loops.

But do we have loops in the class $L_n$ such that we have subloops H in $L_n(m)$ where the o(H) does not divide the $o(L_n(m))$. To this end we have such examples.

***Example 3.3.2***: Let $L_{15}(2) \in L_{15}$ be a loop of order 16. Clearly H = {e, 2, 5, 8, 11, 14} is a subloop of $L_{15}(2)$ which is not a subgroup but this is a S-subloop of $L_{15}(2)$. o(H) = 6 and $o(L_{15}(2))$ = 16. We see 6 ∤ 16.

But this loop is a S-Lagrange group as it contains only subgroups of order 2 and 4 and 2 and 4 divide 16.

**DEFINITION 3.3.5**: *Let L be a loop of finite order, if the order of every S-subloop divides the order of L, we say L is a Smarandache pseudo Lagrange loop. (S-pseudo Lagrange loop).*

**DEFINITION 3.3.6**: *Let L be a finite loop. If L has at least one S-subloop K such that o(K)/o(L) then we say the loop L is a Smarandache weakly pseudo Lagrange loop (S-weakly pseudo Lagranges loop).*

**THEOREM 3.3.5**: *Every S-pseudo Lagrange loop is a S-weakly pseudo Lagrange loop.*

*Proof*: Obvious by the definitions 3.3.5 and 3.3.6.

Now we proceed onto define Smarandache p-Sylow subloops.

**DEFINITION 3.3.7**: *Let L be a finite loop. Let p be a prime such that p/o(L). If L has a S-subloop of order p, then we say L is a Smarandache p-Sylow subloop (S-p-Sylow subloop). It is important to note that $p^r$ need not divide the order of L.*

**DEFINITION 3.3.8**: *Let L be a finite loop. Let p be a prime such that p/o(L). If L has a S-subloop A of order m and B is a subgroup of order p, $B \subset A$, and p/m them we say L is a Smarandache p-Sylow subgroup (S-p-Sylow subgroup).*

We note here that m need not divide the order of L. Still L can be a S-p-Sylow subgroup. Suppose L is a loop such that L is S-subgroup loop then how to define Sylow structure in L.

**DEFINITION 3.3.9**: *Let L be a S-subgroup loop of finite order if every subgroup is either of a prime order or a prime power order and divides o(L) then we call L a Smarandache strong p-Sylow loop (S-strong p-Sylow loop).*

It is a notable fact that all subgroups in a finite S-strong p-Sylow loop should have order $p^\alpha$ where p is a prime and $\alpha$ is a positive integer $\alpha \geq 1$.

Such things occur and we define such class of loops also.

**THEOREM 3.3.6**: *Let $L_n(m) \in L_n$ where n is a prime number. Then every loop in $L_n$ is Smarandache strong 2-sylow loop.*

*Proof*: Given $L_n(m) \in L_n$ and n is a prime number. But order of each $L_n(m)$ in $L_n$ is (n + 1), so 2/n + 1. Every element in $L_n(m)$ is of order two. Hence each $L_n(m) \in L_n$ is a Smarandache strong 2-sylow loop.

The only classical theorem, which is left behind, is the Cayley's theorem. Can we have any nice analogue of Cayley's theorem? The answer is yes; but for this we have to introduce the concept of Smarandache loop homomorphism between loops.

**DEFINITION 3.3.10**: *Let L and L' be two Smarandache loops with A and A' its subgroups respectively. A map $\phi$ from L to L' is called Smarandache loop homomorphism (S-loop homomorphism) if $\phi$ restricted to A is mapped onto a subgroup A' of L'; that is $\phi : A \to A'$ is a group homomorphism. The concept of Smarandache loop isomorphism and automorphism are defined in a similar way. It is important to observe the following facts:*

1. *The map $\phi$ from L to L' need not be even be a loop homomorphism.*
2. *Two loops of different order can be isomorphic.*
3. *Further two loops which are isomorphic for a map $\phi$ may not be isomorphic for some other map $\eta$.*
4. *If L and L' have at least subgroups A and A' in L and L' respectively which are isomorphic then certainly L and L' are isomorphic.*

**Example 3.3.3**: Let the two loops $L_5(3)$ and $L_7(3)$ be given by the following tables:

| • | e | 1 | 2 | 3 | 4 | 5 |
|---|---|---|---|---|---|---|
| e | e | 1 | 2 | 3 | 4 | 5 |
| 1 | 1 | e | 4 | 2 | 5 | 3 |
| 2 | 2 | 4 | e | 5 | 3 | 1 |
| 3 | 3 | 2 | 5 | e | 1 | 4 |
| 4 | 4 | 5 | 3 | 1 | e | 2 |
| 5 | 5 | 3 | 1 | 4 | 2 | e |

$L_5(3)$ is a S-loop of order 6 which is commutative.

Now the table for the loop $L_7(3)$ is given below:

| • | e | 1 | 2 | 3 | 4 | 5 | 6 | 7 |
|---|---|---|---|---|---|---|---|---|
| e | e | 1 | 2 | 3 | 4 | 5 | 6 | 7 |
| 1 | 1 | e | 4 | 7 | 3 | 6 | 2 | 5 |
| 2 | 2 | 6 | e | 5 | 1 | 4 | 7 | 3 |
| 3 | 3 | 4 | 7 | e | 6 | 2 | 5 | 1 |
| 4 | 4 | 2 | 5 | 1 | e | 7 | 3 | 6 |
| 5 | 5 | 7 | 3 | 6 | 2 | e | 1 | 4 |
| 6 | 6 | 5 | 1 | 4 | 7 | 3 | e | 2 |
| 7 | 7 | 3 | 6 | 2 | 5 | 1 | 4 | e |

This is S-loop of order 8 which is non-commutative, but both the loops are S-isomorphic. For take A = {4, e} and A' = {7, e} we have a group isomorphism between the groups A and A' hence they are S-isomorphic as loops.

In view of this we have the following nice result about the new class of loops $L_n$ for n > 3, n odd for varying n.

**THEOREM 3.3.7**: *All S-loops in $L_n$ and for odd n, n > 3 and varying appropriate m, have isomorphic subgroups. So all loops in the class $L_n$ are S-isomorphic loops.*

*Proof*: Follows from the fact that every loop $L_n(m) \in L_n$ has a subgroup of order 2. Hence we have got a S-loop isomorphism between any two loops in $L_n$ even for varying n as is evident by the example.

Now a few problems are left open.

Can we have a class of S-loops which have $S_n$ the symmetric group of degree n to be a subgroup for every appropriate n? If this is achieved certainly extension of Cayley's theorem for S-loops will be possible.

To this end we will have the concept of Smarandache mixed direct product and some related results which will be dealt in this book.

**PROBLEMS:**

1. Does there exist a loop L which has $S_n$ as its subgroup?
2. Find a Moufang loop in which every element is a Cauchy element.
3. Can a Bol loop be a S-Cauchy loop?
4. Give an example of a Bol loop which is S-Lagrange loop.
5. Can all Bol loops be S-Lagrange loops?
6. Will a Bruck-loop be a S-loop? Give at least one example of a Bruck loop which is a S-loop.
7. Give an example of a Bruck loop which is not a S-Lagrange loop.
8. Give an example of a S-loop in which every subgroup is a S-5-Sylow subgroup.
9. Give an example of a S-loop which has S-3-Sylow subloop.
10. Can a loop L have all its subloops to be S-p-Sylow subloops (the order of L is a composite number)?
11. Can a Moufang loop and a Bol loop be S-isomorphic loops? Justify.
12. Give an example of a loop of order 15 which has S-subloops.
13. Give an example of a loop of order 10 which has no S-subloops.
14. Can every non-simple loop have a S-normal subloop? Illustrate with examples.
15. Find an example of a finite loop L which has $S_7$ to be a subgroup.
16. Give an example of a loop L which has the dihedral group of order 8 to be a subgroup.

## 3.4. Smarandache commutative and commutator subloops

In this section we introduce the concept of Smarandache commutative loop, Smarandache strongly commutative loop, Smarandache commutator loop etc. and study them. Further we introduce concepts like Smarandache cyclic loop, Smarandache pseudo commutative loops and characterize them. We illustrate each of these by examples. The main interesting concept about S-loops is that when the S-loop has no S-subloop then the commutator of the loop coincides with the S-commutator of the loop. If the S-loop has many S-subloops we can get more than one S-commutator subloop. This is a unique property enjoyed only by S-loops and not by loops. This study has forced us to propose open research problems in chapter 5.

**DEFINITION 3.4.1**: *Let L be a loop. We say L is a Smarandache commutative loop (S-commutative loop) if L has a proper subset A such that A is a commutative group.*

*Example 3.4.1*: Consider the loop $L_7(3)$ given by the following table:

| • | e | 1 | 2 | 3 | 4 | 5 | 6 | 7 |
|---|---|---|---|---|---|---|---|---|
| e | e | 1 | 2 | 3 | 4 | 5 | 6 | 7 |
| 1 | 1 | e | 4 | 7 | 3 | 6 | 2 | 5 |
| 2 | 2 | 6 | e | 5 | 1 | 4 | 7 | 3 |
| 3 | 3 | 4 | 7 | e | 6 | 2 | 5 | 1 |
| 4 | 4 | 2 | 5 | 1 | e | 7 | 3 | 6 |
| 5 | 5 | 7 | 3 | 6 | 2 | e | 1 | 4 |
| 6 | 6 | 5 | 1 | 4 | 7 | 3 | e | 2 |
| 7 | 7 | 3 | 6 | 2 | 5 | 1 | 4 | e |

Clearly this loop is non-commutative but it is a S-commutative loop. Thus we see the following if a loop is commutative then it is clearly S-commutative leading to the nice result.

**THEOREM 3.4.1**: *Let L be a S-commutative loop then L need not in general be commutative.*

*Proof*: By an example. The loop given in example 3.4.1 is non-commutative but is S-commutative.

In view of this we give a definition of Smarandache strongly commutative loops.

**DEFINITION 3.4.2**: *Let L be a loop. L is said to be a Smarandache strongly commutative (S-strongly commutative) loop if every proper subset which is a group is a commutative group.*

**THEOREM 3.4.2**: *Let L be a S-strongly commutative loop then L is a S-commutative loop.*

*Proof*: Obvious by the definitions 3.4.1 and 3.4.2.

**THEOREM 3.4.3**: *Every S-strongly commutative loop in general need not be a commutative loop.*

*Proof*: The proof is given by an example. Consider the loop $L_7(3)$ given in example 3.4.1. Clearly $L_7(3)$ is non-commutative; but $L_7(3)$ is S-strongly commutative as every proper subset which is a group in $L_7(3)$ is a commutative group. Hence the claim.

**THEOREM 3.4.4**: *Let L be a power-associative loop. L is then a S-commutative loop.*

*Proof*: Given L is a power associative loop so any element generates a group so it must be a cyclic group, hence abelian. So a power associative loop is a S-commutative loop.

**THEOREM 3.4.5**: *Every loop in the class $L_n$ where n is a prime; n > 3 is a S-strongly commutative loop.*

*Proof:* We know every loop in the class $L_n$ when n is a prime is of order n + 1. We have proved the only subloops of $L_n(m)$ are subgroups of order 2 and it has no other subloops. Further every element generates a subgroup of order 2. Hence all loops $L_n(m)$ in $L_n$, when n is a prime is a S-strongly commutative loop for every prime p, p > 3.

**THEOREM 3.4.6**: *Let L be a loop, if Z(L) is a non-trivial centre then L is a S-commutative loop.*

*Proof*: From the result of Pflugfelder [50] in his book on Quasigroups and loops: Introduction shows Z(L) is an abelian subgroup. Hence L is a S-commutative if Z(L) is non-trivial.

**DEFINITION 3.4.3**: *Let L be a loop. L is said to be a Smarandache cyclic loop (S-cyclic loop) if L contains at least a proper subset A which is a cyclic group.*

**DEFINITION 3.4.4**: *Let L be a loop. If every proper subset A of L which is a subgroup is a cyclic group then we say the loop L is a Smarandache strong cyclic loop (S-strong cyclic loop).*

**THEOREM 3.4.7**: *Every S-strong cyclic loop is a S-cyclic loop.*

*Proof*: Given L is a S-loop by the very definitions 3.4.3 and 3.4.4 we see every S-strong cyclic loop is a S-cyclic loop.

It is pertinent to mention here that we cannot in case of loops define the concept of cyclic loop but we are gifted for in case of S-loops we can define S-strongly cyclic loops.

We have in fact classes of loops which are S-strongly cyclic.

**THEOREM 3.4.8**: *Every power associative loop which has no other non-cyclic subgroups in them are S-strong cyclic loops.*

*Proof*: We know if L is a loop which is power associative and if no subgroup which is non-cyclic then clearly L is a S-strongly cyclic loop.

Here we give an example of a S-strongly cyclic group.

**Example 3.4.2**: Let $L_5(4)$ be a loop given by the following table:

| • | e | 1 | 2 | 3 | 4 | 5 |
|---|---|---|---|---|---|---|
| e | e | 1 | 2 | 3 | 4 | 5 |
| 1 | 1 | e | 5 | 4 | 3 | 2 |
| 2 | 2 | 3 | e | 1 | 5 | 4 |
| 3 | 3 | 5 | 4 | e | 2 | 1 |
| 4 | 4 | 2 | 1 | 5 | e | 3 |
| 5 | 5 | 4 | 3 | 2 | 1 | e |

Clearly this is a S-strongly cyclic group which is also power associative.

Now we proceed onto show $L_n$ is a S-strongly cyclic loop when n is a prime.

**THEOREM 3.4.9**: *Let $L_n(m) \in L_n$, n > 3 and n a prime. Then every loop in $L_n(m)$ is cyclic so $L_n(m)$ is a S-strongly cyclic loop.*

*Proof*: Obvious from the fact every element in $L_n(m)$ generates a cyclic group of order two and every distinct pair generates only a subloop or the whole loop.

It is still interesting to note that even strictly non-commutative loop can be a S-commutative loop and S-strongly commutative loop. Also these strictly non-commutative loops can also be S-cyclic loops or S-strongly cyclic loops. We show these by examples.

**Example 3.4.3**: Let $L_5(2)$ be the loop given by the following table:

| • | e | 1 | 2 | 3 | 4 | 5 |
|---|---|---|---|---|---|---|
| e | e | 1 | 2 | 3 | 4 | 5 |
| 1 | 1 | e | 3 | 5 | 2 | 4 |
| 2 | 2 | 5 | e | 4 | 1 | 3 |
| 3 | 3 | 4 | 1 | e | 5 | 2 |
| 4 | 4 | 3 | 5 | 2 | e | 1 |
| 5 | 5 | 2 | 4 | 1 | 3 | e |

This loop is non-commutative and non-associative. So $L_5(2)$ is a strictly non-commutative loop but $L_5(2)$ is a S-strongly commutative loop and S-strongly cyclic loop.

In view of this we have the following:

**THEOREM 3.4.10**: *Let $L_n$ be a class of loops for any $n > 3$, if $n = p_1^{\alpha_1} p_2^{\alpha_2} \ldots p_k^{\alpha_k}$ ($\alpha_i \geq 1$ for $i = 1, 2, \ldots, k$) then it contains exactly $F_n$ loops which are strictly non-commutative and they are:*

1. *S-strongly commutative loops and*
2. *S-strongly cyclic loops where $F_n = \prod_{i=1}^{k} (p_i - 3) p_i^{\alpha_i - 1}$*

*Proof*: It can be shown that every element generates a cyclic group of order two and any two distinct elements different from identity does not generate a subgroup but only a subloop or a S-subloop. Hence all these strictly non-commutative loops are S-strongly cyclic and S-strongly commutative loops.

It is left as an open research problem in Chapter 5 to find S-strongly commutative loops which are not S-strongly cyclic loops.

**DEFINITION 3.4.5**: *Let L be a loop. If A (A proper subset of L) is a S-subloop of L is such that A is a pseudo commutative loop then we say L is a Smarandache pseudo commutative loop (S-pseudo commutative loop) i.e. for every $a, b \in A$ we have an $x \in B$ such that $a(xb) = b(xa)$ (or $(bx)a$), B is a subgroup in A. Thus we see for a loop to be a S-pseudo commutative we need not have the loop L to be a pseudo-commutative loop. If L is itself a pseudo commutative loop then trivially L is a S-pseudo commutative loop.*

**DEFINITION 3.4.6**: *Let L be a loop. If a proper subset $A \subset L$ where A is a S-subloop is such that for every distinct pair $a, b \in A$ we have $a \bullet x \bullet b = b \bullet x \bullet a$ for all $x \in A$ then we say the loop L is a S-strongly pseudo commutative loop.*

The following can be easily derived which is left as an exercise to the reader.

1. Prove all commutative loops are

    a. S-pseudo commutative
    b. S-strongly pseudo commutative

2. Prove a S-strongly pseudo commutative loop is a S-pseudo commutative loop.

*Example 3.4.4*: Let L be a loop given by the table:

|   | e | 1 | 2 | 3 | 4 | 5 |
|---|---|---|---|---|---|---|
| e | e | 1 | 2 | 3 | 4 | 5 |
| 1 | 1 | e | 3 | 5 | 2 | 4 |
| 2 | 2 | 5 | e | 4 | 1 | 3 |
| 3 | 3 | 4 | 1 | e | 5 | 2 |
| 4 | 4 | 3 | 5 | 2 | e | 1 |
| 5 | 5 | 2 | 4 | 1 | 3 | e |

This loop is non-commutative. It is left as an exercise for the reader to verify whether L is pseudo commutative. But in any case the loop has no S-subloops. So we are not in a position to define S-pseudo commutative in this case.

Similarly we can have examples of loops which are commutative but which do not contain S-subloops. This is an example of a commutative loop which has no subloops, then how to define S-pseudo commutativity in these cases?

*Example 3.4.5*: Let $L_{19}(10)$ be a loop which is a commutative loop of order 20. This loop has no S-subloops but $L_{19}(10)$ is a S-loop. Then for this loop we are not in a position to define S-pseudo commutativity.

**DEFINITION 3.4.7**: *Let L be a loop. The Smarandache commutator subloop (S-commutator subloop) of L is denoted by $L^s$ is the subloop generated by all the commutators in A, A a S-subloop of L. ($A \subset L$).*

$L^s = \langle \{x \in A \,/\, x = (y, z) \text{ for some } y, z \in A\} \rangle$ with the usual notation and if L has no S-subloops but L is a S-loop we replace A by L.

Thus we have the notion of S-commutator subloop only when the loop L has a proper S-subloop or L is a S-loop. We see when L has no S-subloop and if L is a S-loop then we have the concept of S-commutator subloop to be one and the same for the loop L and for the S-loop.

**THEOREM 3.4.11**: *Let L be a S-loop which has no S-subloops then $L' = L^s$.*

*Proof*: Obvious from the very definition of S-commutator in loop; we see if L has no S-subloops but is a S-loop we have L' = L$^s$.

**THEOREM 3.4.12**: *Let $L_n(m) \in L_n$ where $L_n(m)$ is assumed to be a non-commutative loop and n a prime. Then we have $L'_n(m) = L_n(m) = L_n^s(m)$ for all $L_n(m) \in L_n$.*

*Proof*: Clear from the fact that each $L_n(m)$ when n is a prime is a S-loop which has no proper S-subloop or even subloops. In view of this we have $L'_n(m) = L_n(m) = L_n^s(m)$ as in case of non-commutative loop we have $L'_n(m) = \langle \{t \in L_n(m) / t = (x, y)$ for some $x, y \in L_n(m)\}\rangle$. Hence the claim.

***Example 3.4.6:*** Let L be a loop given by $L_{21}(11)$. Find subloops of order 3 and 12 which are S-subloops that are S-commutator subloops. This is given as an exercise for the reader as it needs only simple number theoretic methods.

In view of the definition of S-commutator subloops in a loop L, unlike in a loop we can have many S-commutator subloops if L is a loop which has more than one S-subloop.

**THEOREM 3.4.13**: *Let L be a loop having more than one S-subloop. Then the S-commutator subloops of L can have more than one subloop.*

*Proof*: If L has more than one S-subloop then we see if in particular the S-subloops are disjoint i.e. $S_1 \cap S_2 = \{e\}$ where $S_1$ and $S_2$ are two S-subloops of L then certainly L has more than one S-commutator subloop.

**PROBLEMS:**

1. Find all cyclic groups in $L_{19}(3)$.
2. Is $L_{23}(4)$ a S-strongly cyclic loop? Justify your answer.
3. Can $L_{13}(4)$ be a S-strongly commutative loop which is not a S-strongly cyclic loop?
4. Is $L_{15}(8)$ a S-cyclic loop? Justify.
5. Find all the S-subloops of $L_{25}(7)$. Find $L_{25}^s(7)$ for each of these S-subloops. Is $L_{25}^s(7)$ the same for all S-subloops?
6. Give an example of a S-loop L, of odd non-prime order in which L$^s$ = L' = L.
7. Give an example of a S-loop of odd non-prime order which has more than one S-commutator subloop.
8. Find an example of a loop in which every distinct S-subloop has the same S-commutator subloop.

9. Find a S-loop which is S-pseudo commutative of order 15.
10. Find a S-loop which is S-strongly commutative of order 27.

## 3.5. Smarandache associative and associator subloops

In this section we introduce the notion of Smarandache associativity in loops, Smarandache pseudo associativity and Smarandache associator subloops. Here also like Smarandache commutator subloops we can for a given loop have many Smarandache associator subloops which is a distinctly an important property enjoyed only by Smarandache associators. Further these concepts leads us to several open unsolved problems which are enlisted in Chapter 5 of this book.

Further it is exciting to note that the study when the Smarandache associator loop will be coincident with the associator loop and strongly pseudo associator loop. It is not even imaginable to obtain a relation even when L happens to a S-loop having no S-subloops the relation between the sets A(L), PA(L), SPA(L), $L^s$, PA($L^s$) and SPA($L^s$).

**DEFINITION 3.5.1**: *Let L be a loop. We say L is a Smarandache associative loop (S-associative loop) if L has a S-subloop A such that A contains a triple x, y, z (all three elements distinct and different from e, the identity element of A) such that x • (y • z) = (x • y) • z. This triple is called the S-associative triple.*

**DEFINITION 3.5.2**: *Let L be a loop. We say L is a Smarandache strongly associative loop (S-strongly associative loop) if in L every S-subloop has an S-associative triple.*

**THEOREM 3.5.1**: *Let L be a loop. If L is a S-strongly associative loop then L is a S-associative loop.*

*Proof*: Obvious from the very definitions 3.5.1 and 3.5.2.

We have yet another nice theorem which guarantees or gives a condition on a loop L to be a S-associative loop.

**THEOREM 3.5.2**: *Let L be a loop. Let A be a S-subloop of L. If A contains a proper subset B which is a subgroup and |B| > 3 then clearly L is a S-associative loop.*

*Proof*: Follows from the fact if L is a loop such that it has A to be a S-subloop containing a subgroup B of order greater than or equal to 4 then we can take the triple in B\{e} and it will be an S-associative triple of L. Hence the claim.

**COROLLARY 3.5.1**: *If L is a S-associative loop, then L is a S-loop.*

*Proof*: By the very definition of S-associative loop we see L must contain a S-subloop so that L becomes a S-loop.

**THEOREM 3.5.3**: *All S-loops are not S-associative loops but all S-associative loops are S-loops.*

*Proof*: By an example. Consider the S-loop L given by L = $L_{13}(2)$. Clearly $L_{13}(2)$ is a S-loop but it has no S-subloop which has a S-associative triple. Hence the claim. But if L is an S-associative loop we know L has S-subloop so L is a S-loop.

**THEOREM 3.5.4**: *Let $L_n(m) \in L_n$ be the class of loops where n is a prime. None of the loops in this class of loops $L_n$ is a S-associative loop.*

*Proof*: Obvious from the very fact that each $L_n(m) \in L_n$ is a S-loop but every loop $L_n(m)$ has no S-subloop for $L_n(m)$ has subgroups of order two and no subloops of any order other than $L_n(m)$ itself. Hence the claim.

**DEFINITION 3.5.3**: *Let L be a loop if L has a S-subloop A such that for x, y $\in$ A we have (xy)x = x(yx) then we say the loop L is a Smarandache pairwise associative loop (S-pairwise associative loop).*

**DEFINITION 3.5.4**: *Let L be a loop if every S-subloop A of L is a S-pairwise associative loop we say L is a Smarandache strongly pairwise associative loop (S-strongly pairwise associative loop).*

**THEOREM 3.5.5**: *Every S-strongly pairwise associative loop is a S-pairwise associative loop.*

*Proof*: By the very definition of S-pairwise associative loop and S-strongly pairwise associative loop we have the theorem to be true.

**THEOREM 3.5.6**: *If L is a S-strongly pairwise associative loop or a S-pairwise associative loop, L is a S-loop.*

*Proof*: True by the very definition of the two concepts S-pairwise associative and S-strongly pairwise associative.

**THEOREM 3.5.7**: *Every loop in class of loops $L_n$ are S-strongly pairwise associative.*

*Proof*: All loops in the class $L_n$ are S-loops. Further it can be easily verified using simple number theoretic techniques that for every x, y $\in L_n(m) \in L_n$ we have (xy)x = x(yx). Hence the claim.

**DEFINITION 3.5.5**: *Let L be a loop. The Smarandache associator subloop (S-associator subloop) of L is denoted by $L^A$ is the subloop generated by all the associators in A, where A is a S-subloop of L ($A \subset L$) i.e. $L^A = \langle \{x \in A / x = (a,b,c)$ for some $a, b, c \in L\} \rangle$ with the usual notation $\langle B \rangle$ denotes the loop generated by the set B. If L has no S-subloop but L is a S-loop we replace A by L itself. Thus we have the notion of S-associator subloop only when the loop L has a S-subloop or L itself is a S-loop. We see when L has no S-subloops but L a S-loop then the concept of S-associator subloop and the associator subloop coincide.*

**THEOREM 3.5.8**: *Let L be a S-loop which has no S-subloops then we have $A(L) = L^A$ i.e. the associator subloop of L coincides with the S-associator subloop.*

*Proof*: Obvious by the very definition of S-associator subloop and the associator subloop of a S-loop.

**THEOREM 3.5.9**: *Let $L_n(m) \in L_n$ where $L_n(m)$ is a S-loop and has no S-subloops. Then $A(L_n(m)) = L_n^A(m) = L_n(m)$.*

*Proof*: We know for all loops $L_n(m)$ in the class of loops $L_n$ we have $A(L_n(m)) = L_n(m)$, using number theoretic techniques. So for the loop $L_n(m)$ in $L_n$ which are always S-loops we have, if $L_n(m)$ has no S-subloops then $A(L_n(m)) = L_n^A(m) = L_n(m)$. Hence the claim.

The natural question which arises is that what happens in case $L_n(m)$ has S-subloops. Further we see if L is a loop having many S-subloops we will have several S-associator subloops associated with them? Will those S-associator subloops be distinct or identical even for different S-subloops? These problems are proposed as open problems in the chapter 5 of this book.

*Example 3.5.1*: Let $L = L_{45}(8)$ be a loop of order 46 in $L_{45}$. Now $L_{45}(8)$ is a S-loop. Further $A(L_{45}(8)) = L_{45}(8)$. Now consider the S-subloop $A = \{e, 1, 16, 31\}$. A is a S-subloop for $H = \{1,e\}$ is a subgroup of A. We see in fact A itself is a subgroup, now the $L_{45}^A(8) = \{1, 16, 31, e\}$. Clearly $A(L_{45}(8)) \neq L_{45}^A(8)$. From this we have a nice consequence.

**THEOREM 3.5.10**: *Let L be a loop if L is a S-loop having a S-subloop A then in general $A(L) \neq L^A$.*

*Proof*: By an example. Consider the loop $L_{45}(8)$ given in example 3.5.1. Clearly for the S-subloop, $A = \{e, 1, 16, 31\}$ we see $L_{45}^A(8) = \{e, 1, 16, 31\}$ but $A(L_{45}(8)) \neq L_{45}(8)$. Hence the claim.

Now we proceed onto define Smarandache pseudo associativity in loops.

**DEFINITION 3.5.6**: *Let L be a loop. Suppose A is a S-subloop in L, if we have an associative triple a, b, c in A such that (ax)(bc) = (ab)(xc) for some x ∈ A. We say L is a Smarandache pseudo associative loop (S-pseudo associative loop) if every associative triple in A is a pseudo associative triple of A for some x in A. If L is a S-loop having no S-subloops then we replace A by L.*

**DEFINITION 3.5.7**: *Let L be a loop and A be a subloop. If for any associative triple a, b, c ∈ A we have (ax)(bc) = (ab)(xc) for some x ∈ L then we say the triple is strongly pseudo associative. If in particular every associative triple in A is strongly pseudo associative then we say the loop L is Smarandache strongly pseudo associative loop (S-strongly pseudo associative loop).*

If in particular L is a S-loop which has no S-subloop we replace A by L in the definition.

**DEFINITION 3.5.8**: *Let L be a Smarandache pseudo associative loop. Then $PA(L^s) = \langle \{t \in A / (ab)(tc) = (at)(bc) \text{ where } (ab)c = a(bc) \text{ for } a, b, c \in A\} \rangle$ denotes the Smarandache pseudo associator subloop of L (S-pseudo associator subloop of L).*

**DEFINITION 3.5.9**: *Let L be a loop. A be a S-subloop of L. The set $SPA(L^s) = \langle \{t \in L / (ab)(tc) = (at)(bc) \text{ where } (ab)c = a(bc) \text{ and } a, b, c \text{ are in } A\} \rangle$ denotes the Smarandache strongly pseudo associator (S-strongly pseudo associator) subloop of L.*

Clearly $PA(L^s) \subset SPA(L^s)$ this is left for the reader to prove.

**PROBLEMS:**

1. Find an S-associative triple in $L_{47}(3)$.
2. Find $A(L_7(3))$. Is $A(L_7(3)) = L_7^A(3)$ for any S-subloop of $L_7(3)$? Justify your answer.
3. Find two S-subloops A and B in a loop L such that $A^A \neq B^A$ (i.e. S-associator subloop of A ≠ S-associator subloop of B).
4. Find for all S-subloops of $L_{45}(17)$ and their S-associative subloops.
5. Does $L_{45}(17)$ have a S-associative triple? If so give it explicitly.
6. Find a loop of odd order in which the S-associator subloop is the same as associator subloop.
7. Give an example of a loop of order 15 which has a S-subloop and find its S-associator subloop.

8. For the loop $L_7(4)$ find $PA(L_7^S(4))$ and $SPA (L_7^S(4))$.

9. For the loop $L_{21}(11)$ find a S-subloop A and obtain $PA(L_{21}^S(11))$ and $SPA (L_{21}^S(11))$.

10. Find the S-associator subloop of $L_{21}(11)$ and compare it with $PA(L_{21}^S(11))$ and $SPA(L_{21}^S(11))$.

## 3.6. Smarandache identities in Loops

This section is devoted to the introduction of concepts like Smarandache Bol loops, Smarandache Bruck loops etc. Also we define Smarandache strong power associative loop, Smarandache strong Bruck loop and so on. We study those properties and show in case of the class of loops in $L_n$ when n is prime these loops cannot be S-Bruck or S-Bol or S-Moufang but the class of loops are S-power associative.

We obtain some results about these properties. In fact we have proposed several problems of finding conditions for them to be S-strong Bol loops, S-strong Moufang loops etc. We define a Smarandache Bol triple in a loop L. We extend this idea to the case of Smarandache strong Bol triple. Similarly one can define using Bruck identity and Moufang identity.

**DEFINITION 3.6.1**: *A Smarandache Bol loop (S-Bol loop) is defined to be a loop L such that a proper subset A of L which is a S-subloop is a Bol-loop (with respect to the induced operations). That is $\phi \neq A \subset S$.*

Similarly Smarandache Bruck loops (S-Bruck loop), Smarandache Moufang loops (S-Moufang loop), Smarandache WIP loop (S-WIP loop) and Smarandache right (left) alternative loop (S-right (left) alternative loop) are defined. It is pertinent to mention that in the definition 3.6.1 we insist A should be a subloop of L and not a subgroup of L. For every subgroup is a subloop but a subloop in general is not a subgroup. Further every subgroup will certainly be a Moufang loop, Bol loop, Bruck loop and right (left) alternative loop, WIP loop etc. since in groups the operation is associative. Hence only we make the definition little rigid so that automatically we will not have all Smarandache loops to be Smarandache Bol loops, Smarandache Moufang loops etc.

**THEOREM 3.6.1**: *Every Bol loop is a Smarandache Bol loop but every Smarandache Bol loop is not a Bol loop.*

*Proof*: Clearly every Bol loop is a Smarandache Bol loop as every subloop of a Bol loop is a Bol loop. But a Smarandache Bol loop L is one which has a proper subset A which is a Bol loop. Hence L need not in general be a Bol loop.

**THEOREM 3.6.2**: *No loop in the class of loops $L_n(m)$, n an odd prime is a*

1. *Smarandache Bol loop*
2. *Smarandache Bruck loop*
3. *Smarandache Moufang loop*

*Proof*: We know these loops $L_n(m) \in L_n$ when n is an odd prime greater than 3 do not have any proper subloops and it has only subgroups of order 2. So for $L_n(m) \in L_n$ the class of loops are not S-Bol loop, S-Moufang loop or S-Bruck loop. Hence the theorem.

**THEOREM 3.6.3**: *Every finite ARIF loop of odd order is a S-Moufang loop.*

*Proof*: By Corollary 2.14 to Theorem 2.13 of [41] 2001 we have every finite ARIF-loop of odd order is Moufang; so every subloop of this loop will also be Moufang hence ARIF loop is a S-Moufang loop.

**THEOREM 3.6.4**: *The loop $L_n(m) \in L_n$ is a S-WIP loop if and only if $(m^2 - m + 1) \equiv 0 \pmod{n}$ ( n a non-prime).*

*Proof*: We assume n is a non-prime for if we assume n is a prime. $L_n(m)$ has no subloops only subgroups of order 2.

Recall that $L_n(m)$ is a WIP loop if $(xy)z = e$ implies $x(yz) = e$ where

$x, y, z \in L_n(m)$ …1

Now assume $L_n(m)$ is a WIP loop. x, y, z such that $(x - y, n) = 1$ and $z = xy$. Now $z = xy$ implies.

$z \equiv (my - (m-1)x) \pmod{n}$ …2

Since $L_n(m)$ is a WIP loop and $(xy)z = e$ we must have $x(yz) = e$ or $yz = x$.

That is

$x \equiv (mz - (m-1)y) \pmod{n}$ …3

Putting the value of z from (2) in (3) we get

$(m^2 - m + 1)(x - y) \equiv 0 \pmod{n}$ or $(m^2 - m + 1) \equiv 0 \pmod{n}$ …4

Conversely if $(m^2 - m + 1) \equiv 0 \pmod{n}$ then it is easy to see that (3) holds good whenever (2) holds good, that is $(xy)z = e$ implies $x(yz) = e$. x, y and z are distinct elements of $L_n(m)\setminus\{e\}$. However if one of x or y or z is equal to e or x = y then (1) holds trivially. Hence $L_n(m)$ is a WIP loop. So all subloops will be WIP loops only if n is an odd non-prime greater than 3.

**DEFINITION 3.6.2**: *Let L be a loop we call a triple (x, y, z) where x, y, z ∈ L to be a Smarandache Bol triple (S-Bol triple) if the Bol identity ((xy)z)y = x((yz)y) is true. It is to be noted that x, y, z cannot be interchanged in anyway.*

***Example 3.6.1***: Let $L_{15}(8)$ be a loop. We can verify x = 2, y = 4 and z = 13 is a Smarandache Bol triple. Clearly if x = 13, y = 4 and z = 2 is taken we see the Bol identity is not true. Thus we define a concept called Smarandache strong Bol triple.

**DEFINITION 3.6.3**: *Let L be a loop. A triple (x,y,z) where x, y, z ∈ L is called the Smarandache strong Bol triple if the Bol identity is satisfied by all the 6 permutations (x, y, z).*

We similarly define Smarandache Bruck triple, Smarandache strong Bruck triple, Smarandache Moufang triple and Smarandache strong Moufang triple and propose some problems.

**DEFINITION 3.6.4**: *Let L be a loop which is not a diassociative loop if L has a S-subloop A which is diassociative then we say L is a Smarandache diassociative loop (S-diassociative loop).*

Note: If we assume the loop L to be diassociative then trivially L would be S-diassociative that is why we assume L is not a diassociative loop.

**DEFINITION 3.6.5**: *Let L be a non-power associative loop. L is said to be a Smarandache power associative loop (S-power associative loop) if L contains a S-subloop which is a power associative subloop.*

Study of these types have not been carried out in loop theory, i.e. even when the loop is not a diassociative loop, can it have subloops which are diassociative. So these two definitions not only motivates such study of those loops for which the S-subloops happen to hold to a nature or property not possessed by the loop.

Now we proceed on to define Smarandache strong properties.

**DEFINITION 3.6.6**: *Let L be a loop if every proper S-subloop of L satisfies Bol identity. We say the loop L is a Smarandache strong Bol loop (S-strong Bol loop).*

On similar lines we define Smarandache strong Bruck loop, Smarandache strong Moufang loop, Smarandache strong WIP loops, Smarandache strong diassociative loop and Smarandache strong power-associative loop.

The following theorems are true for all Smarandache strong loops having special identities to be satisfied.

**THEOREM 3.6.5**: *If L is a S-strong Bol loop then it is a S-Bol loop.*

*Proof*: Obvious by the very definitions of S-strong Bol loop and S-Bol loops. Similarly we can prove theorems on S-strong Bruck loop, S-strong Moufang loops etc.

**PROBLEMS:**

1. Can the loop $L_{27}(8)$ be a S-Moufang loop?
2. Is $L_{25}(7)$ a S-diassociative loop?
3. Does $L_{49}(9)$ have a S-subloop which is a S-Bruck loop?
4. Give an example of a loop L of order 15 which is S-Moufang loop but L is not a Moufang loop.
5. Does there exist a loop L of order 9, which is a S-diassociative loop, but L is not a diassociative loop?
6. Give an example of a loop of order 21 which is a power associative loop.
7. Find a loop of order 11 which is a S-Bruck loop but L is not a Bruck loop.
8. Can a loop of order 11 which is not a Moufang loop be a S-strong Moufang loop? Justify.
9. Can a loop of prime order p which is not a Bol loop be is a S-strong Bol loop?
10. Give an example of an odd order loop which is not a S-strong Bol loop but is a S-Bol loop.

## 3.7 Some special structures in S-loops

In this section we introduce the concepts of Smarandache left right and middle Nucleus for loops, which have S-subloop, or for S-loops. We also define Smarandache centre, Moufang centre first and second normalizer. We show in general Smarandache first and second normalizer are different and they coincide only in certain loops in $L_n$.

Finally we prove in case of loop $L_n(m) \in L_n$ when n is a prime the concept of Smarandache and the general definition coincide. We prove some results for these class of loops using number theoretic technique. The section ends with problems which are given as exercise for to reader to solve.

**DEFINITION 3.7.1**: *Let L be a loop the Smarandache left nucleus (S-left nucleus) $S(N_\lambda) = \{a \in A / (a, x, y) = e$ for all $x, y \in A\}$ is a subloop of A, where A is a S-subloop of L. Thus we see in case of S-left nucleus we can have many subloops; it is not unique as in case of loops. If L has no S-subloops but L is a S-loop then $S(N_\lambda) = SN_\lambda = N_\lambda$.*

**DEFINITION 3.7.2**: *Let L be a loop, the Smarandache right nucleus (S-right nucleus) $S(N_\rho) = \{ a \in A / (x, y, a) = e$ for all $x, y \in A\}$ is a subloop of L where A is a S-subloop of L. If L has no S-subloops but L is a S-loop then $S(N_\rho) = SN_\rho = N_\rho$.*

**DEFINITION 3.7.3**: *Let L be a loop, the Smarandache middle nucleus (S-middle nucleus). $S(N_\mu) = \{ a \in A / (x, a, y) = e$ for all $x, y \in A\}$ is a subloop of L where A is a S-subloop of L. If L has no S-subloop but L is a S-loop then $S(N_\mu) = SN_\mu = N_\mu$.*

**DEFINITION 3.7.4**: *The Smarandache nucleus $S(N(L))$ of a loop L is a subloop given by $S(N(L)) = SN_\mu \cap SN_\lambda \cap SN_\rho$. It is important to note that unlike the nucleus we can have several Smarandache nucleus depending on the S-subloops. If L has no S-subloops but L is S-loop then $S(N(L)) = N(L)$.*

**DEFINITION 3.7.5**: *Let L be a loop. The Smarandache Moufang centre (S-Moufang centre) $SC(L)$ is the set of elements in a S-subloop A of L which commute with every element of A, that is $SC(L) = \{x \in A / xy = yx$ for all $y \in A\}$. If L has no S-subloops but L is a S-loop then we have $SC(L) = C(L)$. If L has many S-subloops then we have correspondingly many S-Moufang centres $SC(L)$.*

**DEFINITION 3.7.6**: *Let L be a loop, the Smarandache centre (S-centre) $(SZ(L))$ of a loop. L is the intersection of $SZ(L) = SC(L) \cap S(N(L))$ for a S-subloop A of L.*

**DEFINITION 3.7.7**: *Let L be a loop if A is a S-subloop, then the Smarandache first normalizer (S-first normalizer) of A is given by $SN_1(A) = \{a \in L / a \bullet A = A \bullet a\}$.*

**Example 3.7.1**: Let $L_{15}(2) \in L_{15}$ be a loop. Clearly $A = \{e, 1, 4, 7, 10, 13\}$ is a S-subloop of $L_{15}(2)$. Clearly $SN_1(A) = L_5(2)$. It is left for the reader to verify in this case $N_1(A) = SN_1(A)$; but if A is not a S-subloop we do not define S-first normalizer.

**DEFINITION 3.7.8**: *Let L be a loop. A a S-subloop of L, the Smarandache second normalizer (S-second normalizer) of A is given by $SN_2(A) = \{x \in L / xAx^{-1} = A\}$. The two Smarandache normalizer are not equal on any S-subloop in general.*

**Example 3.7.2**: Consider the loop $L_{45}(8)$ where $H_1(15) = \{e, 1, 16, 31\}$, $SN_1(H_1(15)) = L_{45}(8)$ this is left for the reader to verify: $SN_2(H_1(15)) = \{e, 1, 6, 16, 21, 26, 31, 36, 41\} = H_1(5)$. Thus $SN_1(H_1(15)) \neq SN_2(H_1(15))$ in general.

Now a natural question would be can we have S-subloops A in $L_n(m) \in L_n$ for which $SN_1(A) = SN_2(A)$. The answer to this question is yes. We to answer this problem introduce a notation for S-subloops in $L_n(m)$. When we say $L_n(m)$ has subloops we assume n is an odd non-prime for otherwise it will have no subloops.

**Notation**: Let $L_n(m)$ be a loop, n an odd number which is not a prime for every t/n there exist subloops of order $k + 1$ where $k = n/t$. So if t/n we denote the subloop by $H_i(t) = \{e, i, i + t, \ldots, i + (k - t)t\}$, it is a S-subloop of $L_n(m)$ for every $i \leq t$.

***Example 3.7.3***: Now $L_9(8) \in L_9$ the number which divides 9 is 3. So we have the subloop

$$H_1(t) = \{e, 1, 4, 7\}$$
$$H_2(t) = \{e, 2, 5, 8\}$$
$$H_3(t) = \{e, 3, 6, 9\}$$

The tables for them are as follows

Table for $H_1(t)$

|   | e | 1 | 4 | 7 |
|---|---|---|---|---|
| e | e | 1 | 4 | 7 |
| 1 | 1 | e | 7 | 4 |
| 4 | 4 | 7 | e | 1 |
| 7 | 7 | 4 | 1 | e |

Table for $H_2(t)$

|   | e | 2 | 5 | 8 |
|---|---|---|---|---|
| e | e | 2 | 5 | 8 |
| 2 | 2 | e | 8 | 5 |
| 5 | 5 | 8 | e | 2 |
| 8 | 8 | 5 | 2 | e |

Table for $H_3(t)$

|   | e | 3 | 6 | 9 |
|---|---|---|---|---|
| e | e | 3 | 6 | 9 |
| 3 | 3 | e | 9 | 6 |
| 6 | 6 | 9 | e | 3 |
| 9 | 9 | 6 | 3 | e |

Thus all subloops in $L_9(8)$ are abelian groups. So this is also a Smarandache Hamiltonian loop as every subloop is an abelian group of order 4.

Now we obtain a condition for the normalizers to be coincident in the loops $L_n(m) \in L_n$ (n not a prime).

**THEOREM 3.7.1**: *Let $L_n(m) \in L_n$ and $H_i(t)$ be its S-subloop, then $SN_1(H_i(t)) = SN_2(H_i(t))$ if and only if $(m^2 - m + 1, t) = (2m - 1, t)$.*

*Proof*: Let $L_n(m) \in L_n$ and $H_i(t)$ be a S-subloop of $L_n(m)$ first we show the first S-normalizer $SN_1(H_i(t)) = H_i(k)$ where $k = t/d$ and $d = (2m - 1, t)$ we use only simple number theoretic arguments and the definition. $SN_1(H_i(t)) = \{j \in L_n(m) / jH_i(t) = H_i(t)j\}$ is the S-first normalizer of $H_i(t)$. It is left for the reader to verify $jH_i(t) = H_i(t)j$ if and only if $(2m - 1)(i - j) \equiv 0 \pmod{t}$ for $j \notin H_i(t)$. More over for $j \in H_i(t)$ we have $jH_i(t) = H_i(t)j$.

So the possible values of j in $L_n(m)$ for which $jH_i(t) = H_i(t)j$ are given by e, i, i + k, i + 2k, ..., i + ((n/k) - i) k. Thus $SN_1(H_i(t)) = H_i(k)$ using this fact we see $k = t/d$, $d = (2m - 1, t)$. Now as we have worked for Smarandache first normalizer we can work for Smarandache second normalizer and show that $SN_2(H_i(t)) = H_i(k)$ here $k = t/d$ and $d = (m^2 - m + 1, t)$. This working is left as an exercise to the reader. Using these two facts we can say $SN_1(H_i(t)) = SN_2(H_i(t))$ if and only if $(m^2 - m + 1, t) = ((2m - 1), t)$. Hence the claim.

**THEOREM 3.7.2**: *Let L be a S-loop which has no S-subloops. Then*

1. $S(N_\lambda) = N_\lambda$
2. $S(N_\rho) = N_\rho$
3. $S(N_\mu) = N_\mu$
4. $SN(L) = N(L)$
5. $SC(L) = C(L)$
6. $SZ(L) = Z(L)$

*Proof*: The proof is obvious by the definitions of all these concepts and their corresponding Smarandache definitions.

**THEOREM 3.7.3**: *Let $L_n(m) \in L_n$ where n is a prime then $SN(L_n(m)) = \{e\}$.*

*Proof*: Now we know when n is a prime. $L_n(m) \in L_n$ has no S-subloops but are S-loops. So we have $N(L_n(m)) = \{e\}$ as we know $SN(L_n(m)) = N(L_n(m))$. Hence we have $SN_\rho(L_n(m)) = N_\rho(L_n(m))$, $SN_\mu(L_n(m)) = N_\mu(L_n(m))$ $SN_\lambda(L_n(m)) = N_\lambda(L_n(m))$. Now $N(L) = N_\mu \cap N_\lambda \cap N_\rho$.

It is enough if we prove $N_\lambda(L_n(m)) = \{e\}$ then it would imply $N(L_n(m)) = \{e\}$. So we shall prove $N_\lambda(L_n(m)) = \{e\}$. If $x \neq e \in N_\lambda(L_n(m))$, choose elements y, z $\in L_n(m)$ such that $z \neq x$, $xy \neq z$ and $yz \neq x$. Now $x \in N_\lambda(L_n(m))$ implies $(x, y, z) = e$ which is possible only when x = z. However e ∈ $N_\lambda(L_n(m))$ as (e, x, y) = e for all x, y, ∈ $L_n(m)$. Hence the theorem.

**THEOREM 3.7.4**: *Let $L_n(m) \in L_n$, n a prime, then S-Moufang centre of $L_n(m)$ is either {e} or $L_n(m)$.*

*Proof*: Using the fact $SC(L_n(m)) = C(L_n(m))$ and e ∈ $L_n(m)$ we see if e ≠ i ∈ $C(L_n(m))$ then i • j = j • i for all j ∈ $L_n(m)$. i • j = j • i implies (2m – 1) (i – j) ≡ 0(mod n) choose j such that (|j – i|, n) = 1 which implies (2m – 1) ≡ zero(mod n). Hence $L_n(m)$ is commutative in which case $C(L_n(m)) = SC(L_n(m)) = L_n(m)$. If $L_n(m)$ is non-commutative $C(L_n(m)) = \{e\} = SC(L_n(m))$.

**THEOREM 3.7.5**: *Let $L_n(m) \in L_n$ where n is a prime. Then $NZ(L_n(m)) = Z(L_n(m)) = \{e\}$.*

*Proof*: If n is a prime we have $NZ(L_n(m)) = Z(L_n(m))$. So it is enough if we prove $Z(L_n(m)) = \{e\}$.

We know SZ(L) = S(NCL) ∩ SC(L) = N(L) ∩ C(L) (if n is a prime) = {e} by our earlier results.

**PROBLEMS:**

1. Find all S-Moufang centres for the loop $L_{51}(11)$.
2. Find all S-centres for the loop $L_{55}(13)$.
3. Find $Z(L_{27}(5))$ and $NZ(Z_{25}(5))$. Which of them is the larger subloop?
4. Find $SN_\mu$, $SN_\lambda$ and $SN_\rho$ for all S-subloops in $L_{57}(17)$.
5. For problem 4 find $SN(L_{57}(17))$.
6. Find $SN_1$ and $SN_2$ for all subloops in $SN(L_{57}(29))$.
7. If P and Q be subloops in a loop L such that P ⊂ Q what can you say about
    i. $SN_\rho(P)$ and $SN_\rho(Q)$.
    ii. SZ(P) and SZ(Q)
    iii. $SN_1(P)$ and $SN_1(Q)$ and
    iv. $SN_2(P)$ and $SN_2(Q)$.
8. For all S-subloops in $L_{35}(9)$ find $SN_1, SN_2$, SZ, SC and SN.
9. For what S-subloops A in L. $SN_1$ = SC = $SN_2$ = SZ = SN?
10. Verify problem 9 in case of S-subloops in $L_{15}(8)$.

## 3.8. Smarandache mixed direct product of loops

In this section we define a new notion called Smarandache mixed direct product of loops and prove these loops got as Smarandache mixed direct product are S-loops. We define S-loop II. Using the definition of Smarandache mixed product we are able to get the Smarandache Cauchy theorem for S-loops. Further in this section we extend two classical concepts, call them as Smarandache Lagrange criteria and Smarandache Sylow criteria. These are substantiated by examples. In fact for each of these concepts we have a class of loops, which satisfies it. Finally we give several problems to the reader as exercise to solve. As solving of these problems alone can make the reader get a deeper understanding of S-loops.

**DEFINITION 3.8.1**: *Let $L = L_1 \times S_n$ be the direct product of a loop and a group. We call this the Smarandache mixed direct product of loops (S-mixed direct product of loops). We insist that L is a loop and one of $L_1$ or $S_n$ is group. We do not take all groups in the S-mixed direct product or all loops in this product. Further we can extend this to any number of loops and groups that is if $L = L_1 \times L_2 \times \ldots \times L_n$ where some of the $L_i$'s are loops and some of them are groups. Then L is called the S-mixed direct product of n-terms.*

**THEOREM 3.8.1**: *Let $L = L_1 \times L_2$ be the direct product of a loop and a group be a S-mixed direct product of loops. Then L is a loop.*

*Proof*: By the very definition of S-mixed product we see one of $L_1$ or $L_2$ is a group. So L = $L_1 \times L_2$ under component wise operation can at most be a loop and never a group as one of $L_1$ or $L_2$ is a loop; so the operation is non-associative that is $L_1 \times L_2 = \{(m_1, m_2) / m_1 \in L_1$ and $m_2 \in L_2\}$.

**COROLLARY 3.8.1**: *Let $L = L_1 \times \ldots \times L_n = \{(m_1, m_2, \ldots, m_n) / m_i \in L_i; i = 1, 2, \ldots, n\}$ where some of the $L_i$'s are loops and some of the $L_i$'s are groups. L is a S-mixed direct product of loops and L is a loop under component wise operation. Now because of this S-mixed direct product alone we are in a position to get Cayley's theorem for Smarandache loops.*

**THEOREM 3.8.2**: *Let $L = L_1 \times \ldots \times L_n$ be a S- mixed direct product of loops, then L is a Smarandache loop.*

*Proof*: Given $L = L_1 \times \ldots \times L_n$ is a S-mixed direct product of loops so we have one of the $L_i$ is a group. Thus in the loop L we have $\{(e_1, e_2, \ldots e_{i-1}, m_i, e_{i+1}, \ldots, e_n)$ such that $m_i \in L_i$ and each $e_j$ is the identity element of $L_j$ for $j = 1, 2, \ldots, n\} = H_i$. It is easily verified $H_i$ is a proper subset of L, which is a group. So, L is a S-loop.

*Example 3.8.1*: Let $L = L_5(2) \times S_3 = \{(m, p)/ m \in L_5(2)$ and $p \in S_3\}$. L is a loop in fact a S-loop. It is easily checked L is a non-commutative loop but L is a S-commutative loop but L is not a S-strongly commutative loop.

In fact L is not a Hamiltonian loop and L may not have S- subloops and may or may not have S-normal subloops. This leads to the following.

**THEOREM 3.8.3**: *Let $L = L_1 \times L_2$, be a S-mixed direct product of loops with $L_1$ a loop and $L_2$ a group. The loop L is not S-simple. (This is as in the case of direct product of groups where we have in case of direct product of groups have normal subgroups. Similarly S-mixed direct product of loops contain S-normal subloop that is L is not simple).*

*Proof*: Let $L = L_1 \times L_2$, clearly $H = \{(m_1, e_2) \mid m_1 \in L_1$ and $e_2 \in L_2$, the identity element of $L_2\}$ is a S-normal subloop of L. Hence the claim.

**DEFINITION 3.8.2**: *Let L be a loop. H a subgroup of L. H is said to be a normal subgroup of L if $H = aHa^{-1}$ for all $a \in L$.*

**THEOREM 3.8.4**: *Let $L = L_1 \times L_2$ be a Smarandache mixed direct product of loop. The loop L has a normal subgroup.*

*Proof*: Given $L = L_1 \times L_2$ is a S-loop, where $L_1$ is a loop and $L_2$ is a group. Let $N = \{(e_1, g_2) \mid e_1$ is the identity of the loop $L_1$ and $g_2 \in L_2\}$. Clearly N is a normal subgroup in L. This concept forces us and paves way for the definition of Smarandache loops of level II.

**DEFINITION 3.8.3**: *Let L be a loop, if L has a non-empty subset A of L which is a normal subgroup in L then we say L is a Smarandache loop of level II denoted by S-loop II.*

**THEOREM 3.8.5**: *Let L be a S-loop of level II then L is a S-loop of level I that is L is a S-loop (by default of notation when we say S-loop it is implied that it is a S-loop I).*

*Proof*: Clearly if L is a S-loop II. L has a proper subset A such that A is normal subgroup under the operation of L. So A is a subgroup under the operations of L. Hence the claim.

**THEOREM 3.8.6**: *Every S-loop need not in general be a S-loop II.*

*Proof*: By an example, consider the loop $L_5(2)$ given by the following table

| • | e | 1 | 2 | 3 | 4 | 5 |
|---|---|---|---|---|---|---|
| e | e | 1 | 2 | 3 | 4 | 5 |
| 1 | 1 | e | 3 | 5 | 2 | 4 |
| 2 | 2 | 5 | e | 4 | 1 | 3 |
| 3 | 3 | 4 | 1 | e | 5 | 2 |
| 4 | 4 | 3 | 5 | 2 | e | 1 |
| 5 | 5 | 2 | 4 | 1 | 3 | e |

Clearly $A_1 = \{e, 1\}$, $A_2 = \{e, 2\}$, $A_3 = \{e, 3\}$, $A_4 = \{e, 4\}$ and $A_5 = \{e, 5\}$ are the only subgroups of $L_5(2)$. But it is easily verified that none of them are normal subgroups of $L_5(2)$. So $L_5(2)$ is not a S-loop II but is a S-loop.

***Example 3.8.2***: Let $L = L_5(2) \times S_3$, L is a S-loop II. For $H = \{(e_1, p_1) \mid e_1 \in L_5(2)$ is the identity element of L and $p_i \in S_3\}$ is a normal subgroup of L. Also take

$$P = \left\{ \begin{pmatrix} 1 & 2 & 3 \\ 1 & 2 & 3 \end{pmatrix}, \begin{pmatrix} 1 & 2 & 3 \\ 2 & 3 & 1 \end{pmatrix}, \begin{pmatrix} 1 & 2 & 3 \\ 3 & 1 & 2 \end{pmatrix} \right\}$$

the normal subgroup contained in $S_3$. Then $K = \{(e_1, p_i) \mid e_1$ is the identity element of $L_5(2)$ and $p_i \in P\}$ is a normal subgroup of L.

**THEOREM 3.8.7**: *All S-mixed direct product loops L are S-loop II.*

*Proof*: Obvious from the fact L is a direct product of both loops and groups if $L_K$ is a group in the direct product then in $L = L_1 \times \ldots \times L_n$ we have $H = \{(e_1, e_2, \ldots, g_k, \ldots, e_n) \mid g_k \in L_k, e_i \in L_i$ for $i \neq k$ and $i = 1, 2, \ldots, n$ are the identity elements in $L_i\}$. H is a normal subgroup of L. So L is a S-loop II.

We will discuss about S-loop II separately later in this book. Now we are initiating to get some analogue of the classical Cayley's theorem for groups to S-loops. Cayley's theorem states "Every group is isomorphic to a subgroup of $S_n$ for some appropriate n". (where $S_n$ is the group of permutation of n elements i.e. symmetric group of degree n need not be finite unless specified).

Now an analogue of Cayley's theorem for S-loops, that is, the Smarandache Cayley theorem for S-loops.

**THEOREM 3.8.8** (S-Cayley theorem for S-loops): *Every S-loop is S-isomorphic to a S-loop got from the S-mixed direct product of loops.*

*Proof*: Given L is a S-loop; so L contains a proper subset A where A is a group. Now consider the S-mixed direct product loop $L' = L_1 \times \ldots \times L_k$ where one of $L_i$ is itself the group A or a $L_i$ is a group containing A as a subgroup. The latter is also possible by the classical Cayley's theorem for groups. So we can always define a S-isomorphism from L to L'. Hence the theorem. This theorem will be called Smarandache Cayley theorem for S-loops.

Thus our S-mixed direct product of loops and S-isomorphism of loops help us to get a Cayley's theorem for S-loops.

The next natural question would be can we have the classical Lagranges' theorem for S-loops. The answer is yes only for the S-loops built using S-mixed direct product loops. Thus we call this concept the Smarandache Lagrange theorem.

**DEFINITION 3.8.4**: *Let L be a finite S-loop if the order of every subgroup divides the order of L then we say S satisfies the Smarandache Lagrange's criteria for S-loops. (If L is not a S-loop there is no meaning for this condition. Secondly it is not a must all S-loops must satisfy the Smarandache Lagranges criteria but we have a class of loops which satisfies the S-Lagranges criteria.)*

*Example 3.8.3* : Let $L_5(2)$ be a S-loop of order 6. It is easily verified every subgroup in $L_5(2)$ is of order 2 and 2/6. Hence $L_5(2)$ satisfies the S-Lagranges criteria. A natural question would be does the order of every subgroup divide the order of the loop L where L is a finite loop.

The answer to this question is now by an example.

*Example 3.8.4* : Consider the loop $L_{45}(8)$. This is a loop of order 46. Clearly H = {e, 1, 16, 31} is a subgroup of $L_{45}(8)$. Clearly $4 \nmid 46$. Hence the claim. But we have class of S-loops which satisfies the S-Lagranges criteria.

**THEOREM 3.8.9** : *Every S-loop in the class of loops $L_n$ where n is an odd prime greater than 3 satisfies S-Lagranges criteria.*

*Proof*: Given $L_n(m) \in L_n$ where n is an odd prime. So each $L_n(m)$ can have at most n-subgroups of order 2 and no subgroups of any other order. Further $|L_n(m)| = n + 1$ = even, we see S-Lagranges criteria is satisfied by the class of loops.

Thus for every prime n > 3, we see we have classes of loops satisfying S-Lagrange criteria. Do we have any other loop which satisfy S-Lagranges criteria the answer is yes.

**THEOREM 3.8.10** : *Let $L = L_1 \times \ldots \times L_n$ where each $L_i$ is finite and one of the $L_i$ is a group. Let $L_n$ be the group. Assume no $L_i$ is a S-loop that is no loop $L_i$ has a proper subset which is a group. Then $L = L_1 \times \ldots \times L_n$ satisfies S-Lagrange criteria.*

*Proof*: Obvious by the very definition of L. It is left for the reader to verify S-Lagrange criteria. Clearly $|L| = |L_1| \times |L_2| \times \ldots \times |L_n|$ and by our classical Lagranges theorem for groups each subgroup of $L_r$ divides the order of $L_r$. As only $L_r$ has subgroups and $L_r$ is a finite group the result is true.

On similar lines we define Smarandache Sylow criteria.

**DEFINITION 3.8.5**: *Let L be a finite loop. If for every prime p/|L| the loop L has a subgroup of order p or a subgroup of a power of p then we say the loop L satisfies Smarandache Sylow (S-Sylow) criteria.*

**Example 3.8.5**: Let $L = L_5(2) \times S_3$ be a S-mixed direct product loop. $|L| = 36$. $2/36$ and $3/36$. L has subgroups of order two and order 3 so L satisfies S-Sylow criteria.

**Example 3.8.6**: Let $L = L_5(2)$, $|L_5(2)| = 6$ does not satisfy S-Sylow criteria for $L_5(2)$ does not have any subgroup of order 3.

It is interesting to note that $L = L_5(2) \times S_3$ satisfies S-Sylow criteria but $L_5(2)$ does not satisfy S-Sylow criteria. This is a very unique property solely enjoyed by the S-mixed direct product alone. For in case of all algebric structures with their product this is impossible as even in case of lattice if one lattice is non-distributive, another distributive the direct product gives only a non-distributive lattice. But S-mixed direct product has enabled us to overcome this problem.

Another natural question would be does all loops of $L_n(m) \in L_n$ never satisfy S-Sylow criteria. The answer is not so by an example.

**Example 3.8.7**: Let $L_7(2)$ be a loop. Clearly $L_7(2)$ has only subgroups of order 2 only in fact 7 groups of order 2. $|L_7(2)| = 8 = 2^3$. $2/8$. Hence all loops in $L_7$ satisfy S-Sylow criteria. Thus this class of loops $L_7$ satisfy S-Sylow criteria.

One more problem is can we say only $L_n$ when n is a prime can give loops which satisfy S-Sylow criteria. The answer is no. For by another example.

**Example 3.8.8**: Let $L_{11}(2) \in L_{11}$ is a loop of order 12, $2/12$ and $3/12$. $L_{11}(2)$ has no subgroups of order 3 but $L_{11}(2) \in L_{11}$ where 11 is a prime. $L_{11}(2)$ does not satisfy S-Sylow criteria.

**Example 3.8.9**: Let $L_{15}(8) \in L_{15}$ be a loop of order 16. Only $2/16$ and $L_{15}(8)$ has subgroups of order 2. Hence $L_{15}(8)$ satisfies S-Sylow criteria.

**THEOREM 3.8.11**: *Let $L_n(m) \in L_n$ if n is equal $2^r - 1$. Then $|L_n(m)| = 2^r$ so that every loop in the class $L_{2^r-1}$ satisfies S-Sylow criteria.*

*Proof*: Obvious form the very definition of S-Sylow criteria and the fact $|L_n(m)| = 2^r$ so only 2 is the prime which divides $|L_n(m)|$. Hence the claim.

In view of this we have the following. Let $L_n(m) \in L_n$ be a loop in $L_n$ if n is not of the form $2^r - 1$ then can we say $L_n(m)$ does not satisfy S-Sylow criteria? We have several examples in support of this.

***Example 3.8.10***: Let L = $L_9(2) \times S_3$. L is S-loop which does not satisfy S-Sylow criteria.

**THEOREM 3.8.12**: *Let L = $L_n(m) \times S_{n+1}$ be the S-mixed direct product loop. Then L satisfies S-Sylow criteria.*

*Proof*: Since $|S_{n+1}| = (n+1)!$ so all prime factors of $|L_n(m)| = n + 1$ are also factors of $(n + 1)!$ So obviously L will satisfy S-Sylow criteria.

This is the least condition we can have so that the S-mixed direct product of the loop $L_n(m) \in L_n$ satisfies S-Sylow criteria.

**PROBLEMS**:

1. Give an example of a loop of order 15 which does not satisfy S-Lagrange criteria.
2. Find an example of a S-loop of order 21 which has both 7-Sylow subgroup and 3-Sylow subgroup.
3. Does there exist a loop of order 21 which does not satisfy both S-Sylow criteria and S-Lagranges criteria?
4. Will it be possible for all S-mixed direct product loops to satisfy S-Lagrange criteria? Justify your answer.
5. Let L = $L_{13}(2) \times G$ (where G = $\langle g / g^7 = 1 \rangle$ is the S-mixed direct product of loops. Does L satisfy S-Lagranges criteria? Will L satisfy S-Sylow criteria?
6. Give an example of a S-loop II of order 26.
7. Give an example of a S-loop of order 25 which is not a S-loop II.
8. Can we say all prime order loops are not S-loops?
9. Is every S-mixed direct product of loop a S-loop II?

### 3.9. Smarandache cosets in loops

In this section we study and introduce the Smarandache coset representation in loops when the loop has a subgroup i.e. only when L is a S-loop. This section defines Smarandache right (left) coset of H in L. Several interesting examples are given like cosets in a group, Smarandache cosets does not in general partition the loop or have some nice representation. All these are explained by self-illustrative examples.

This problem leads to the definition of Smarandache right (left) coset equivalence sets in a loop L related to a subgroup A of L. In chapter 5 several problems are given for the interested reader to develop an interest and carry out research on the Smarandache cosets of loops.

**DEFINITION 3.9.1**: *Let L be a S-loop. We define Smarandache right cosets (S-right cosets) in L as follows:*

*Let $A \subset L$ be the subgroup of L and for $m \in L$ we have Am = {am/a $\in$ A}. Am is called a S-right coset of A in L.*

*Similarly we can for any subgroup A of L define Smarandache left coset (S-left coset) for some $m \in L$ as mA = {ma / a $\in$ A}. If the loop L is commutative then only we have mA = Am. Even if L is S-commutative or S-strongly commutative still we need not have Am = mA for all $m \in L$.*

**Example 3.9.1**: Let $L_5(2)$ be the non-commutative loop given by the following table:

| • | e | 1 | 2 | 3 | 4 | 5 |
|---|---|---|---|---|---|---|
| e | e | 1 | 2 | 3 | 4 | 5 |
| 1 | 1 | e | 3 | 5 | 2 | 4 |
| 2 | 2 | 5 | e | 4 | 1 | 3 |
| 3 | 3 | 4 | 1 | e | 5 | 2 |
| 4 | 4 | 3 | 5 | 2 | e | 1 |
| 5 | 5 | 2 | 4 | 1 | 3 | e |

Let A = {e, 1} be the subgroup of the S-loop $L_5(2)$. The S-right coset of A is

$$A \bullet 1 = \{e, 1\} \quad A \bullet 2 = \{2, 3\}$$
$$A \bullet 3 = \{3, 5\} \quad A \bullet 4 = \{4, 2\}$$
$$A \bullet 5 = \{5, 4\}$$

Thus the S-right cosets do not partition $L_5(2)$.

Consider the S-left cosets of A

$$1 \bullet A = \{e, 1\} \quad 2 \bullet A = \{5, 2\}$$
$$3 \bullet A = \{3, 4\} \quad 4 \bullet A = \{4, 3\}$$
$$5 \bullet A = \{5, 2\}$$

But the S-left cosets has partitioned $L_5(2)$. Thus we cannot make any comment about the partition.

**Example 3.9.2**: Let $L_7(4)$ be the commutative loop given by the following table

| • | e | 1 | 2 | 3 | 4 | 5 | 6 | 7 |
|---|---|---|---|---|---|---|---|---|
| e | e | 1 | 2 | 3 | 4 | 5 | 6 | 7 |
| 1 | 1 | e | 5 | 2 | 6 | 3 | 7 | 4 |
| 2 | 2 | 5 | e | 6 | 3 | 7 | 4 | 1 |
| 3 | 3 | 2 | 6 | e | 7 | 4 | 1 | 5 |
| 4 | 4 | 6 | 3 | 7 | e | 1 | 5 | 2 |
| 5 | 5 | 3 | 7 | 4 | 1 | e | 2 | 6 |
| 6 | 6 | 7 | 4 | 1 | 5 | 2 | e | 3 |
| 7 | 7 | 4 | 1 | 5 | 2 | 6 | 3 | e |

Let A = {e, 5}, consider S-right coset of A.

$$A \bullet 1 = \{1, 3\} \quad A \bullet 2 = \{2, 7\}$$
$$A \bullet 3 = \{3, 4\} \quad A \bullet 4 = \{4, 1\}$$
$$A \bullet 5 = \{5, 1\} \quad A \bullet 6 = \{6, 2\}$$
$$A \bullet 7 = \{7, 6\}$$

Clearly A does not partition $L_5(2)$.

Take A = {e, 5} consider S-left coset of A

$$1 \bullet A = \{1, 3\} \quad 2 \bullet A = \{2, 7\}$$
$$3 \bullet A = \{3, 4\} \quad 4 \bullet A = \{4, 1\}$$
$$5 \bullet A = \{5, 1\} \quad 6 \bullet A = \{6, 2\}$$
$$7 \bullet A = \{6, 7\}$$

Clearly S-cosets do not get partitioned in this case. But mA = Am for all m ∈ $L_7(4)$ as $L_7(4)$ is a commutative loop.

Take B = {e, 4}

$$B \bullet 1 = \{1, 6\} \quad B \bullet 2 = \{2, 3\}$$
$$B \bullet 3 = \{3, 7\} \quad B \bullet 5 = \{5, 1\}$$
$$B \bullet 6 = \{6, 5\} \quad B \bullet 7 = \{7, 2\}$$

This also does not partition $L_7(4)$. Clearly Bm = mB. Just we have seen for the new class of loops; the S-right coset and S-left coset when the loop is non-commutative and the loop is commutative.

In both these loops we do not have subgroups of other order. But one of the important observations which is made by these problems are:

We have a set of elements in $L_5(2)$ such that the set with a specified subgroup A gives coset decomposition in a disjoint way which is a partition of $L_5(2)$. For example in the loop $L_5(2)$ we have for the subgroup A, the elements {2, 5} and {3, 4} in $L_5(2)$ are such that we get $L_5(2) = A \cup A \bullet 2 \cup A \bullet 5$ i.e. {1. e} $\cup$ {2, 3} $\cup$ {4. 5} = $L_5(2)$ which is a disjoint union. Further using the set of elements {3, 4} we see $L_5(2) = A \cup A \bullet 3 \cup A \bullet 4 = $ {1, e} $\cup$ {3, 5} $\cup$ {4, 2}.

So unlike in a group we are able to divide the loop into equivalence classes not using all elements but at the same time we cannot say the two sets partition in the same way.

This is possible only when the loop is a non-commutative one for the loop $L_7(4)$ which is a commutative loop we are not able to get any such relation. So one more relevant question would be should n of the loop $L_n(m)$ be a prime number? We will first illustrate examples before we try to answer these questions.

***Example 3.9.3***: Consider the non-commutative loop $L_7(3)$ given by the following table:

| $\bullet$ | e | 1 | 2 | 3 | 4 | 5 | 6 | 7 |
|---|---|---|---|---|---|---|---|---|
| e | e | 1 | 2 | 3 | 4 | 5 | 6 | 7 |
| 1 | 1 | e | 4 | 7 | 5 | 6 | 2 | 5 |
| 2 | 2 | 6 | e | 5 | 1 | 4 | 7 | 3 |
| 3 | 3 | 4 | 7 | e | 6 | 2 | 5 | 1 |
| 4 | 4 | 2 | 5 | 1 | e | 7 | 3 | 6 |
| 5 | 5 | 7 | 3 | 6 | 2 | e | 1 | 4 |
| 6 | 6 | 5 | 1 | 4 | 7 | 3 | e | 2 |
| 7 | 7 | 3 | 6 | 2 | 5 | 1 | 4 | e |

Let $A_1$ = {1, e}, $2A_1$ = {6, 2}, $3A_1$ = {3, 4}, $4A_1$ = {4, 2}, $5A_1$ = {5, 7}, $6A_1$ = {6, 5} and $7A_1$ = {7, 3}.

S-right coset decomposition by $A_1$ is $A_1$ = {1, e}, $A_1 2$ = {2, 4}, $A_1 3$ = {3,7}, $A_1 4$ = {4, 3}, $A_1 5$ = {5, 6}, $A_1 6$ = {6, 2} and $A_1 7$ = {7, 5}. The set {2, 3, 5} and {4, 6, 7} are such that $L_7(3) = A_1 \cup$ {2, 4} $\cup$ {3, 7} $\cup$ {5,6} = $A_1 \cup A_1 2 \cup A_1 3 \cup A_1 5$. Similarly $L_7(3) = A_1 \cup A_1 4 \cup A_1 6 \cup A_1 7 = $ {e, 1} $\cup$ {4, 3} $\cup$ {6, 2} $\cup$ {7, 5}.

S-left coset decomposition by $A_1$ is $A_1$= {1, e}, $2A_1$ = {6, 3}, $3A_1$ = {3, 4}, $4A_1$ = {4, 2}, $5A_1$ ={5, 7}, $6A_1$ = {6, 5} and $7A_1$ = {7, 3}. Here also the set {2, 3, 5} and {4, 6, 7} are such that $L_7(3) = A_1 \cup 2A_1 \cup 3A_1 \cup 5A_1 = $ {e, 1} $\cup$ {2, 6} $\cup$ {3, 4} $\cup$ {5, 7}.

$L_7(3) = A_1 \cup 4A_1 \cup 6A_1 \cup 7A_1 = \{e, 1\} \cup \{4, 2\} \cup \{6, 5\} \cup \{3, 7\}$. The chief thing to be noticed is that for the same subgroup the sets associated with the S-left coset and S-right coset are the same.

Consider $A_2 = \{e, 4\}$. S-right coset decomposition by $A_2$. $A_2 = \{e, 4\}$, $A_2 \bullet 1 = \{1, 2\}$, $A_2 \bullet 2 = \{2, 5\}$, $A_2 \bullet 3 = \{3, 1\}$, $A_2 \bullet 5 = \{5, 7\}$, $A_2 \bullet 6 = \{6, 3\}$ and $A_2 7 = \{7, 6\}$.

$\{2, 3, 7\}$ and $\{1, 5, 6\}$ are such that we have $L_7(3) = A_2 \cup A_2 \bullet 2 \cup A_2 \bullet 3 \cup A_2 \bullet 7 = \{e, 4\} \cup \{2, 5\} \cup \{3, 1\} \cup \{6, 7\}$ and $L_3(7) = \{e, 4\} \cup \{1, 2\} \cup \{5, 7\} \cup \{3, 6\} = A_2 \cup A_2 1 \cup A_2 5 \cup A_2 6$.

Thus we see depending on the subgroup the sets also differ but we have two sets or two decompositions by $A_2$ or by any subgroup in $L_7(3)$. Now we will see the problem when n = odd non-prime.

*Example 3.9.4*: Let $L_9(8)$ be the loop given by the following table:

| • | e | 1 | 2 | 3 | 4 | 5 | 6 | 7 | 8 | 9 |
|---|---|---|---|---|---|---|---|---|---|---|
| e | e | 1 | 2 | 3 | 4 | 5 | 6 | 7 | 8 | 9 |
| 1 | 1 | e | 9 | 8 | 7 | 6 | 5 | 4 | 3 | 2 |
| 2 | 2 | 3 | e | 1 | 9 | 8 | 7 | 6 | 5 | 4 |
| 3 | 3 | 5 | 4 | e | 2 | 1 | 9 | 8 | 7 | 6 |
| 4 | 4 | 7 | 6 | 5 | e | 3 | 2 | 1 | 9 | 8 |
| 5 | 5 | 9 | 8 | 7 | 6 | e | 4 | 3 | 2 | 1 |
| 6 | 6 | 2 | 1 | 9 | 8 | 7 | e | 5 | 4 | 3 |
| 7 | 7 | 4 | 3 | 2 | 1 | 9 | 8 | e | 6 | 5 |
| 8 | 8 | 6 | 5 | 4 | 3 | 2 | 1 | 9 | e | 7 |
| 9 | 9 | 8 | 7 | 6 | 5 | 4 | 3 | 2 | 1 | e |

Take $A = \{e, 7\}$ as the subgroup of $L_9(8)$.

The S-right coset decomposition by A is $A \bullet 1 = \{1, 4\}$, $A \bullet 2 = \{2, 3\}$, $A \bullet 3 = \{3, 2\}$, $A \bullet 4 = \{4, 1\}$, $A_5 = \{5, 9\}$, $A_6 = \{6, 8\}$, $A \bullet 8 = \{8, 6\}$ and $A \bullet 9 = \{5, 9\}$. Thus in this case we see a nice coset decomposition into partition and into classes. $\{1, 2, 6, 5\}$ and $\{4, 3, 8, 9\}$. Here the nicety is

$$A \bullet 1 = A \bullet 4$$
$$A \bullet 2 = A \bullet 3$$
$$A \bullet 5 = A \bullet 9$$
$$A \bullet 6 = A \bullet 8$$

Thus it is verified $L_9(8)$ is got as a disjoint union of cosets. Now take the subgroup of order 4. B = {e, 1, 4, 7}. Find the S-right coset decomposition by B.

$$B \bullet 2 = \{2, 9, 6, 3\}$$
$$B \bullet 8 - \{8, 3, 9, 6\}$$
$$B \bullet 5 = \{5, 6, 3, 9\}$$
$$B \bullet 3 = \{3, 8, 5, 2\}$$
$$B \bullet 6 = \{6, 5, 2, 8\}$$
$$B \bullet 9 = \{9, 2, 8, 5\}$$

So we in case of subgroups of order four this does not work out by this illustration.

Now we get a loop still of higher order say $L_{15}(14)$ and study the coset representation of their subgroups. We will consider the coset representation of the subgroup A = {e, 4}

| | |
|---|---|
| A • 1 = {1, 7} | A • 2 = {2, 6} |
| A • 3 = {3, 5} | A • 5 = {5, 3} |
| A • 6 = {6, 2} | A • 7 = {7, 1} |
| A • 8 = {8, 15} | A • 9 = {9, 14} |
| A • 10 = {10, 13} | A • 11 = {11, 12} |
| A • 12 = {12, 11} | A • 13 = {13, 10} |
| A • 14 = {14, 9} | A • 15 = {15, 8} |

Thus the two sets {1, 2, 3, 8, 9, 10, 11} and {5, 6, 7, 12, 13, 14, 15} are such that $L_{15}(14)$ has a unique decomposition as union of right cosets by these two sets.

Consider the left coset representation by A = {e, 4}

| | |
|---|---|
| 1 • A = {1, 13} | 2 • A = {2, 15} |
| 3 • A = {3, 2} | 5 • A = {5, 6} |
| 6 • A = {6, 8} | 7 • A = {7, 10} |
| 8 • A = {8, 12} | 9 • A = {9, 14} |
| 10 • A = {10, 1} | 11 • A = {11, 3} |
| 12 • A = {12, 5} | 13 • A = {13, 7} |
| 14 • A = {14, 9} | 15 • A = {15, 11} |

{1, 3, 5, 7, 8, 9, 15} and {2, 6, 10, 11, 12, 13, 14} are such that those cosets representations given by the loop $L_{15}(4)$. But we see the sets for the S-right coset is different from the sets for the S-left coset. These observations helps us to define the following:

**DEFINITION 3.9.2**: *Let L be a loop, A a subgroup of L such that $B = \{x_1 \ldots x_t\}$ be the subset of L with $L = \bigcup_{i=1}^{t} Ax_i$ ($Ax_i \cap Ax_j = \phi$, $i \neq j$) and let $C = \{y_1, \ldots y_t\}$ a subset of L such that $L = \bigcup_{i=1}^{t} Ay_i$ ($Ay_i \cap Ay_j = \phi$, $i \neq j$) further $C \cap B = \phi$.*

*Then we say the set C and B are Smarandache right coset equivalent subsets (S-right coset equivalent subsets) of L related to the subgroup A with S-right coset representation.*

*Similarly we define Smarandache left coset equivalent subsets (S-left coset equivalent subsets) of L related to the subgroup A with S-left coset representation.*

It is pertinent to note in general the S-right coset equivalent subsets of the group A need not be the same as S-left coset equivalent subsets of A. Further it is very essential to note that we are not in general guaranteed of such decomposition. This is also well illustrated by examples. We do think if L happens to be a commutative loop we may not in general have such S-right coset equivalent subsets i.e. representation of L as $L = \bigcup_{i=1}^{t} Ax_i$ may not be possible with $Ax_i \cap Ax_j = \phi$ if $i \neq j$.

We define the decomposition whenever it exists as the Smarandache right (left) coset representation of the loop L relative to the subgroup A of L. We have proposed several research problems for the reader in this direction which is given in Chapter 5.

**PROBLEMS:**

1. Find for the loop $L_5(4)$ the S-left and S-right coset decomposition relative to the subgroup $A = \{1, 4\}$.
2. For the loop $L_{11}(3)$ find the S-left coset representation of $A = \{e, 5\}$. Find also the S-right coset representation related to $A = \{e, 5\}$. Are the equivalent set of S-right coset same as S-left coset.
3. Find for the same loop $L_{11}(3)$ when $A = \{e, 8\}$ the S-right coset and S-left coset. Compare it with problem 2.
4. Does a loop of odd order have a subgroup A such that it has S-right or S-left coset representation? Illustrate with examples.
5. Find for the subgroup $A = \{e, 8\}$ in $L_{45}(8)$ the S-left (right) coset representation.
6. For the subgroup $\{e, 7\}$ of $L_{45}(8)$ find S-left (right) coset representation.
7. Compare problems 6 and 5; also find for $B = \{e, 38\}$ the S-left (right) coset representation.

## 3.10. Some special type of Smarandache loops

In this section we introduce 3 special types of S-loops some of them very recently studied by M. K. Kinyon [41] viz. we define Smarandache RIF loops, Smarandache ARIF loops and Smarandache Steiner loops. We prove using the result of M. K. Kinyon [41, 42] certain odd order Moufang loops are S-loops.

**THEOREM 3.10.1**: *All small Frattini Moufang loops (SFM-loops) are S-loops.*

*Proof*: Clearly SMF loops are Moufang p-loop L with a central subgroup Z of order p (Tim Hsu [63]). Hence SMF loops are S-loops.

**DEFINITION 3.10.1**: *A Smarandache RIF loop is a loop L which has a S-subloop A, such that A is a RIF-loop.*

*So we see a Smarandache RIF loop need not be a RIF loop. Further when we say A is a S-subloop we do not include the condition that A can be a subgroup of L.*

**DEFINITION 3.10.2**: *Smarandache ARIF loop is a loop (S-ARIF loop) L which has a S-subloop A which is a ARIF loop by itself. Thus we do not need L to be a ARIF loop only if a S-subloop of L happens to be a ARIF loop then we say L is a S-ARIF loop.*

Thus only from the results of M. Kinyon 2001 [41], we see ARIF loops are S-ARIF loop further we get a class of odd order Moufang loops which are S-loops. Without this we would be pondering for odd order S-loops.

Steiner loops arise naturally in combinations, since they correspond uniquely to Steiner triple systems, specifically the Steiner loop L corresponds to the triple systems $\{\{x, y, xy\}: x \neq y$ and $x, y \neq 1\}$ on $L \setminus \{1\}$ but we may have in a loop L which has a proper subloop to satisfy the conditions for it to be a Steiner loop. So we are once again necessitated to define Smarandache – Steiner loops.

**DEFINITION 3.10.3**: *Let L be a loop, L is said to be a Smarandache Steiner loop (S-Steiner loop) if L contains a S-subloop A such that A is a Steiner loop.*

As RIF loops includes Steiner loops we have S-Steiner loops (please refer M.K.Kinyon [41])

**PROBLEMS:**

1. Can any of the loops in the class $L_n$ be S-Steiner loops?
2. Can any of the loops in the class $L_n$ be S-RIF loops?

3. Does the class of loops $L_n$ have S-ARIF loops? If so give example.
4. Give an example of a loop L which is not a Steiner loop but a S-Steiner loop.
5. Explain by an example a loop which is not a ARIF loop but is a S-ARIF loop.
6. Does their exist a S-RIF loop which is not a RIF loop?
7. Find a loop of even order which is Moufang and at the same time S-RIF loop.

# Chapter four

# PROPERTIES ABOUT S-LOOPS

In this chapter we deal mainly with S-loops II, define S-subloop II, Smarandache Moufang loops II and likewise all loops which satisfy an identity. In section 1 we define S-loops of level II to satisfy special identities. In the second section, Smarandache homomorphisms of level II are defined. Also Smarandache Nuclei II etc are defined. Section 3 is fully devoted to the Smarandache representation of S-loop II. The Smarandache isotopes are studied in section 4. The final section is completely utilized to the study of S-hyperloops, S-A-hyperloops. This study can be restricted only to the class of new loops as the hyperstructure have meaning only when they are integers modulo n. Every section ends with a list of problems for the reader to solve.

## 4.1 Smarandache loops of level II

This section is devoted to study properties of S-loop II by defining it as in the case of S-loop I. We do not claim that all properties, which have been defined and studied for S-loop I have been introduced and studied in the case of S-loop II in this chapter. What we mainly aim in this chapter is to sketch some definitions and properties of S-loops II as we assume the reader may try many analogous results. Just for the sake of completeness we recall the definition of Smarandache loops of level II and do not think it is a repetition.

**DEFINITION 4.1.1**: *Let L be a loop, L is said to be a Smarandache loop of level II if L has a normal subgroup. (A $\subset$ L where A is a subgroup of L such that mA = Am for all m $\in$ L is a normal subgroup.)*

Several properties were dealt in the chapter 3. Here we view it in a different way and distinct from the results derived in chapter 3.

**DEFINITION 4.1.2**: *Let L be a loop. We say L has a Smarandache subloop II (S-subloop II) if A $\subset$ L is a S-loop II that is A is a subloop and has a normal subgroup B in it. (We do not demand B to a normal subgroup of L it is sufficient that if B happens to be a normal subgroup in A).*

**THEOREM 4.1.1**: *Let L be a loop if L has a S-subloop II L need not be a S-loop II but it is a S-loop.*

*Proof*: From the very definition of S-subloop II we need L to contain a subloop, A which has a normal subgroup relative to A. So we can say L has a subgroup, hence L is S-loop. But in general it need not be a S-loop II as B ⊂ A, which is a normal subgroup of A, B may not be a normal subgroup in L.

We obtain a necessary and sufficient condition for a loop L to have a S-subloop II and L to be a S-loop II.

**THEOREM 4.1.2**: *Let L be a loop and A be a S-subloop II. L is a S-loop II only if at least one of the normal subgroups in a S-subloop II is a normal subgroup of L.*

*Proof*: The proof is simple and the reader is advised to give the proof of this theorem.

**DEFINITION 4.1.3**: *Let L be a loop. L is said to be a Smarandache Moufang loop of level II (S-Moufang loop II) if L has a S-subloop II, A ⊂ L such that (xy) (zx) = (x (yz)) x for all x, y, z ∈ A. Clearly even if L is a Moufang loop still L may fail to be a Smarandache Moufang loop II for L may not contain any normal subgroups. This is once again a main deviation from the definition of S-Moufang loops.*

*On similar lines we define Smarandache Bol loop of level II, Smarandache Bruck loop of level II, Smarandache di-associative loop of level II and Smarandache power associative loop of level II. We see none of these properties (say Bol) held by a loop guarantees of it to be a Smarandache Bol loop of level II.*

**DEFINITION 4.1.4**: *Let L be a loop, L is said to be a Smarandache weak inverse property loop (S-WIP loop) II if L has a S-subloop II A(A ⊂ L) such that (xy) z = e imply x(yz) = e for all x, y, z ∈ A, e identity element of L.*

**Example 4.1.1**: Let $L = L_5(2) \times S_3$ be a S-loop II and $L_1 = L_7(2) \times S_4$ be a S-loop II. We have a map $\phi: A \to A_1$ where $A = \{e\} \times A_3$ and $A_1 = \{e\} \times A_4$ we have $\phi: A \to A_1$ is a S-homomorphism.

**Example 4.1.2**: Let $L = L_5(2) \times S_3$ and $L_1 = L_9(5) \times G$ where $G = \langle g / g^3 = 1 \rangle$ be any two S-loop II. We see $L_1$ and L are S-isomorphic loops. We say two S-loop II are isomorphic. So for two S-loop II to be S-isomorphic it is not essential the loops must be of same order. So two S-loop II are isomorphic if the two loops have normal subgroups of same order which are isomorphic as groups.

**DEFINITION 4.1.5**: *Let L be a loop, L is said to be a Smarandache ARIF loop of level II (S-ARIF II) if L has a subset A which is a S-subloop II i.e. B ⊂ A is a normal subgroup of A, and A is a ARIF loop.*

**THEOREM 4.1.3**: *Every commutative ARIF loop is a S-ARIF loop II.*

*Proof*: Follows from the very definition and the nature of commutativity and the fact every ARIF loop is a diassociative.

We define S-representation of loops level II and the concept of S-principal isotopes of level II relative to S-subloops II. It is left for the reader to study and analyse these concepts in an analogous way as it was done for S-loops. Finally it is for the reader to solve these properties of S-loop II and differentiate it from S-loop I.

The natural or the classical question would be the study of lattice of subloops of S-loop II, lattice of S-subloops of S-loop II, the lattice of normal subgroups of S-loop II, and the lattice of subgroups of S-loop II. We first illustrate by examples.

***Example 4.1.3***: Let $L = L_5(2) \times S_3$ be the mixed product which is a S-loop II. The loop L has only one subloop given by $A = L_5(2) \times \{1\}$ where $\{1\}$ is the identity element of $S_3$. The lattice of subloops of L has the following figure:

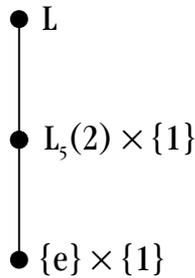

The lattice of S-subloops of L are $B_1 = L_5(2) \times \{1\}$, $B_2 = L_5(2) \times \{1, p_1\}$, $B_3 = L_5(2) \times \{1, p_2\}$, $B_4 = L_5(2) \times \{1, p_3\}$ and $B_5 = L_5(2) \times \{1, p_5, p_4\}$.

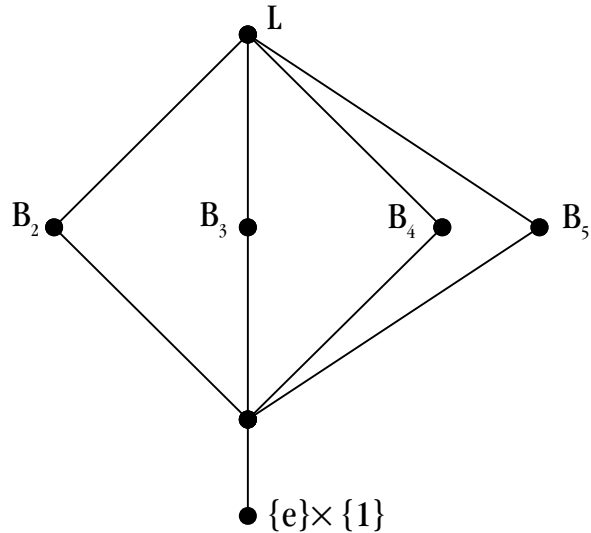

The lattice of subgroups of L is given by $C_1 = \{i, e\} \times \{1\}$, $C_2 = \{i, e\} \times \{1, p_1\}$, $C_3 = \{i, e\} \times \{1, p_2\}$, $C_4 = \{i, e\} \times \{i, p_3\}$, $C_5 = \{i, e\} \times \{1, p_4, p_5\}$ and $C_6 = \{1, e\} \times S_3$.

We have 36 subgroups of this form together with $L_5(2) \times S_3$. So we do not venture to draw the lattice diagram of these subgroups. We proceed on to find the normal subgroups of L, $D_0 = \{e\} \times \{1\}$, $D_1 = \{e\} \times \{1, p_4, p_5\}$ and $D_2 = \{e\} \times S_3$ having the following lattice representation:

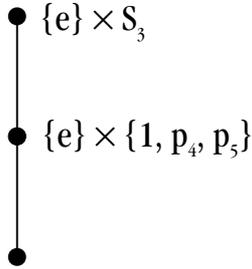

***Example 4.1.4***: To study only the subgroups of a S-loop we take the following simple example L, where $L = L_5(3) \times G$, here $G = \langle g / g^{19} = 1 \rangle$. Clearly L is a S-loop II.

The subgroups of L are: $A_0 = \{e\} \times \{1\}$, $A_1 = \{e, 1\} \times \{1\}$, $A_2 = \{e\} \times \{1\}$, $A_3 = \{e, 3\} \times \{1\}$, $A_4 = \{e, 4\} \times \{1\}$, $A_5 = \{e, 5\} \times \{1\}$, $A_6 = \{e\} \times G$, $A_7 = \{e, 1\} \times G$, $A_8 = \{e, 2\} \times G$, $A_9 = \{e, 3\} \times G$, $A_{10} = \{e, 4\} \times G$, $A_{11} = \{e, 5\} \times G$ together with the largest element L for the lattice is to be formed with 13 elements.

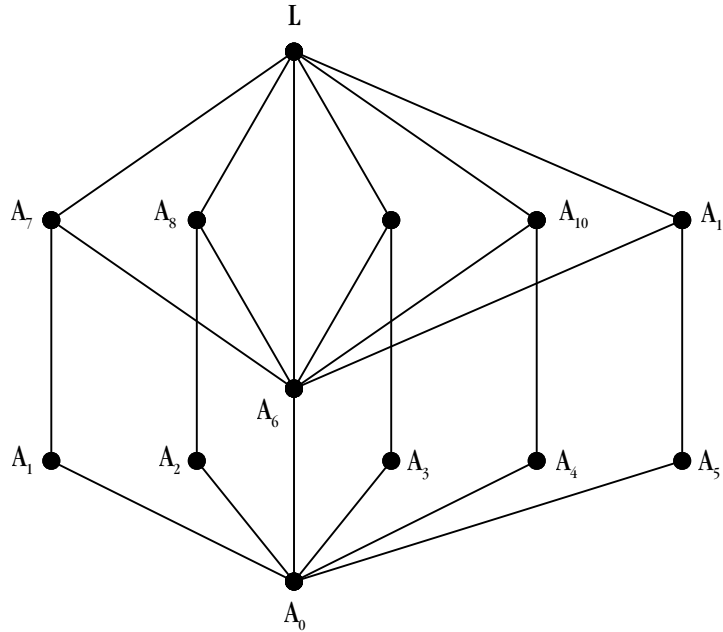

This lattice is also non-modular. Thus the study of lattice formation of subgroups, S-subloops, normal subgroups of the S-loop II is an interesting one as we cannot make any decision about the structure of these lattices.

**THEOREM 4.1.7**: *Let $L_n(m) \in L_n$ n a prime. If $L_n(m) \times G$ where G is a cyclic group of order p, p a prime, then we have $2p + 3 = 2n + 3$ subgroups forming a non-modular non-distributive lattice.*

*Proof*: We have n = p, p a prime. $A_i = \{e, i\} \times \{1\}$ for i = 1, 2, …, p gives p-subgroups. $B_i = \{e, i\} \times G$ for i = 1, 2, …, p gives another p-subgroups. Take $A = \{e\} \times \{1\}$, $B = \{e\} \times G$ and $L = L_n(m) \times G$. Thus we have 2p + 3 subgroups with the following lattice diagram:

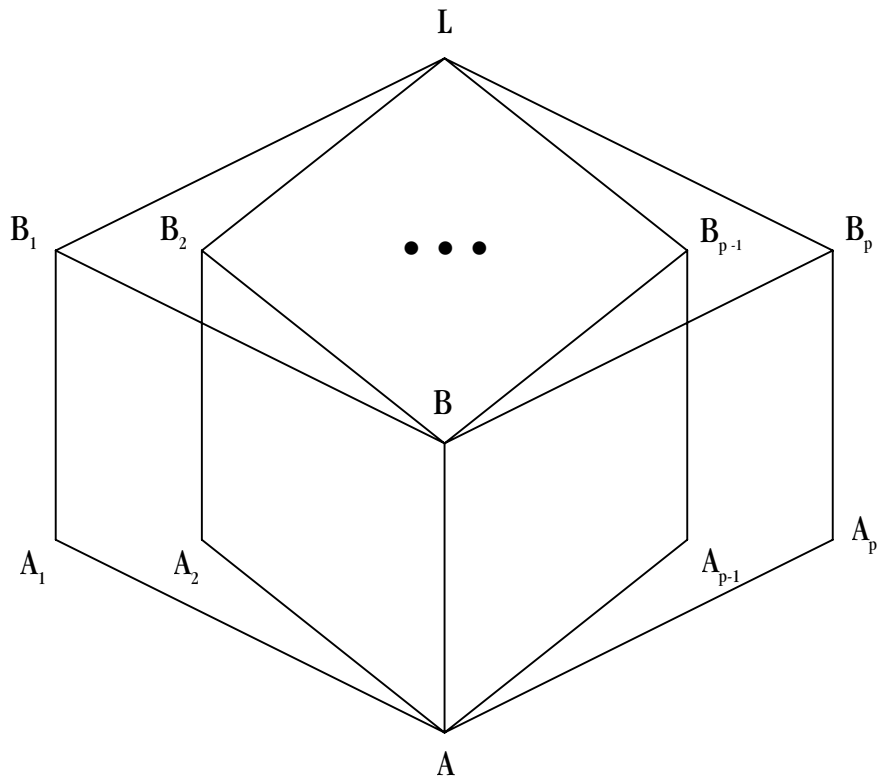

which is easily verified to be a non-distributive non-modular lattice.

**PROBLEMS**:

1. Can any loop in the class $L_n(m) \in L_n$ be a S-loop II? (We know all loops in $L_n$ are S-loops. Here n > 3 and n odd and (m, n) = 1 = (m − 1, n)).
2. Can the loop $L_7(2)$ be a S-loop II?
3. Is the loop $L_5(3)$ a S-loop II?
4. How many loops in $L_n(m) \in L_n$ are S-loops II. $\left(m \neq \dfrac{n+1}{2}\right)$?
5. Find a S-Moufang loop II that is not a Moufang loop.
6. Find a S-Bruck loop II of order 17.
7. Find a S-Bol loop II of order 10, which is not a Bol loop.
8. Give an example of a loop L which is not a S-Bol loop II but which is only a S-Bol loop I.

9. For a loop of order 18 find
   i. S-principal representation II.
   ii. S-isotope II.
10. Give example of a class of loops which are S-loop II and which are also S-Moufang II.

## 4.2 Properties of S- loop II

In this section we define Smarandache loop homomorphism of level II and study the concept of Nuclei, centre, Moufang centre for S-loops II, we also define Smarandache Lagrange criteria II and Smarandache Sylow criteria II in case of S-loops II. It is left for the reader to develop relevant analogue for S-loop II using S-loops.

**DEFINITION 4.2.1**: *Let L and L' be two S-loop II. We say a map $\phi : L$ to $L'$ is a Smarandache loop homomorphism II if $\phi: A \to A'$ where A and A' are normal subgroups of L and L' respectively and $\phi$ is a group homomorphism from A to A'.*

We define Smarandache isomorphism II and Smarandache automorphism II from L to L' in a similar way.

Note: We cannot get any S-loop homomorphism II from a S-loop II to S-loop if the S-loop is not a S-loop II. Here we proceed to define Moufang centre, nuclei etc in case of S-loop II.

**DEFINITION 4.2.2**: *Let L be a loop. A be a S-subloop II of L. Then Smarandache left Nucleus II of loop L is defined as $SS(N_\lambda) = \{a \in A / (a, x, y) = e$ for all $x, y \in A\}$ is a subloop of A.*

*Similarly we define $SS(N_\mu)$ and $SS(N_\rho)$ the Smarandache middle nucleus II and Smarandache right nucleus II respectively.*

**DEFINITION 4.2.3**: *Let L be a loop. SSN(L) the Smarandache nucleus II is defined to be $SSN(L) = SS(N_\lambda) \cap SS(N_\mu) \cap SS(N_\rho)$. If the loop L has no S-subloop II but L is itself a S-loop II then we replace the S-subloop by the S-loop II itself in which case we have $SS(N_\lambda) = N_\lambda$, $SS(N_\mu) = N_\mu$ and so on.*

**DEFINITION 4.2.4**: *Let L be a loop. The Smarandache Moufang centre II of L is defined relative to the S-subloop II A, as $SSC(L) = \{x \in A / xy = yx$ for all $y \in A\}$. If L has no S-subloop II but L is a S-loop II then we replace A by L itself in which case $SSC(L) = C(L)$.*

**DEFINITION 4.2.5**: *Let L be a loop, A ⊂ L be the S-subloop II, the centre of the loop L; SSZ(L) = SSC(L) ∩ SSN(L). If L does not have a S-subloop II but L is itself a S-loop II then replace A by L. Now we see in case of S-loop II we can have like S-loops several S-nuclei II, S-centre II, S-Moufang centre II but unlike loops which has a unique nucleus, centre, Moufang centre and so on.*

Now one more question is Smarandache Lagrange criteria II (S-Lagrange criteria II) and Smarandache Sylow criteria II. We define these two concepts and leave all the properties about S-loops II to the reader.

**DEFINITION 4.2.6**: *Let L be a finite S-loop II. L is said to satisfy Smarandache Lagrange criteria II if the order of every normal subgroup of L divides | L |.*

**DEFINITION 4.2.7**: *Let S be a finite S-loop II. L is said to satisfy Smarandache Sylow criteria II if for every prime p / | L | there is a normal subgroup of order p.*

Similarly we define other related concepts for S-loops II.

**PROBLEMS:**

1. Does the loop $L_{51}(26)$ satisfy
    i. S-Lagrange criteria II?
    ii. S-sylow criteria II?
2. Find the Smarandache Nuclei II for $L_{51}(26)$.
3. Is $SSN(L_{51}(26)) = SSC(L_{51}(26)) = SSZ(L_{51}(26))$? Justify your answer.
4. Does there exist a loop in $L_{15}$ for which SS(N) = S(N)?
5. Find a loop of odd order in which SSC(L) = S(C(L) = C(L).
6. Can we have a non-commutative loop L for which SSC(L) is the same for all S-subloop II?
7. Do we have a loop of odd order in which for all S-subloops we have the same SSZ(L)?
8. Find loop which has many S-subloops but a unique SN(L).
9. Find a S-loop isomorphism of two loops of order 18 and order 28.

## 4.3 Smarandache representation of a finite loop L

In this section we introduce the concept of Smarandache representation of a finite loop, this definition is possible only when the loop has a S-subloop; other wise we have the concept of Smarandache pseudo representations only when L is a S-loop. Thus when L is not a S-loop we do not have any Smarandache representation for them.

**DEFINITION 4.3.1**: *Let L be a finite loop. A be a S-subloop of L (A should not be a group). For $\alpha \in A$ we define a right multiplication $R_\alpha$ as a permutation of the subloop A as follows:*

*$R_\alpha: x \to x \bullet \alpha$ we will call the set $\{R_\alpha \mid \alpha \in A\}$ the Smarandache right regular representation (S-right regular representation) of the loop L or briefly the representation of L.*

*If a loop L has no proper S-subloops then the loop cannot have Smarandache right regular representations. Thus for Smarandache right regular representation to exist for a loop L, L must have S-subloops.*

**THEOREM 4.3.1**: *The class of loops $L_n(m) \in L_n$, n a prime are S-loops but have no S-subloops so these loops do not have S-right regular representation but has right regular representations.*

*Proof*: All loops $L_n(m) \in L_n$ when n is a prime are S-loops for $A_i = \{e, i\}$ is a subgroup for every $i \in \{1, 2, \ldots, n\}$. But this class of loops do not have even subloops so all the more they do not have S-subloops. Hence we cannot speak of S- right regular representation for these loops.

We illustrate the right regular representation of a loop from $L_n(m) \in L_n$ when n is a prime.

***Example 4.3.1***: Let $L_7(4)$ be a loop given by the following table:

|   | e | 1 | 2 | 3 | 4 | 5 | 6 | 7 |
|---|---|---|---|---|---|---|---|---|
| e | e | 1 | 2 | 3 | 4 | 5 | 6 | 7 |
| 1 | 1 | e | 5 | 2 | 6 | 3 | 7 | 4 |
| 2 | 2 | 5 | e | 6 | 3 | 7 | 4 | 1 |
| 3 | 3 | 2 | 6 | e | 7 | 4 | 1 | 5 |
| 4 | 4 | 6 | 3 | 7 | e | 1 | 5 | 2 |
| 5 | 5 | 3 | 7 | 4 | 1 | e | 2 | 6 |
| 6 | 6 | 7 | 4 | 1 | 5 | 2 | e | 3 |
| 7 | 7 | 4 | 1 | 5 | 2 | 6 | 3 | e |

The right regular representation of the loop $L_7(4)$ is

$$I$$
$$(e\ 1)\ (2\ 5\ 3)\ (4\ 6\ 7)$$
$$(e\ 2)\ (1\ 5\ 7)\ (3\ 6\ 4)$$
$$(e\ 3)\ (1\ 2\ 6)\ (4\ 7\ 5)$$
$$(e\ 4)\ (1\ 6\ 5)\ (2\ 3\ 7)$$

$$\begin{array}{l}(e\ 5)\ (1\ 3\ 4)\ (2\ 7\ 6)\\(e\ 6)\ (1\ 7\ 3)\ (2\ 4\ 5)\\(e\ 7)\ (1\ 4\ 2)\ (3\ 5\ 6)\end{array}$$

where I is the identity permutation of the loop $L_7(4)$. To overcome this problem we define Smarandache pseudo representation of loops.

**DEFINITION 4.3.2**: *Let L be a loop. Suppose L has no S-subloops only subgroups B then we define Smarandache pseudo representation (S-pseudo representation) of loops as the set $\{R_\alpha \mid \alpha \in B\}$.*

In the example 4.2.1 we see the loop has S-pseudo representation as every element generates a subgroup of order 2.

Now we see some interesting properties satisfied by these representations. Let $p_1$ = (e 1) (2 5 3) (4 6 7) clearly $\langle p_1, e \rangle$ generates a cyclic group of order 6. Consider $p_2$ = (e 2) (1 5 7) (3 6 4), $\langle p_2, e \rangle$ generates a cyclic group of order 6. Thus we are facing a situation in which each representation associated with each subgroup gives a cyclic group of order 6. Finally it can be checked that none of the permutation in the representation of $L_7(4)$ when producted gives a representation in the same set that is if $\{p_1, p_2, p_3, p_4, p_5, p_6, p_7, e\}$ = P is the set, $p_i \bullet p_j \notin$ P.

Thus Smarandache pseudo permutations of a loop are so pseudo even closure cannot be contemplated.

**THEOREM 4.3.2**: *Let $L_n(m) \in L_n$ and n an odd prime k is a fixed positive integer such that $(m - 1)^k \equiv (-1)^k \pmod{n}$, then $L_n(m)$ has only S-pseudo permutation and any permutation ($a \neq e$) in the representation of $L_n(m)$ is a product of a 2 cycle and t, k-cycles where $t = \dfrac{n-1}{k}$.*

*Proof*: The above example is a proof. Every S-pesudo representation is a product of a 2 cycles where $t = \dfrac{n-1}{k}$. Suppose take n = 11, m = 4 then m − 1 = 3 so in the loop $L_{11}(4)$ we have $(m - 1)^{10} + (-1)^9 \equiv 0 \pmod{11}$, so k = 10. Now we get a 2-cycle and 1 10-cycle. For example I is the identity permutation. (e 1) (2, 9, 10, 7, 5, 11, 4, 3, 6, 8) like wise we get with each (e 2) (e 3) and … (e 11).

The proof of theorem 4.3.2 is left as an exercise to the reader and the proof involves only number theoretic techniques.

We propose some problems for the reader.

**PROBLEMS:**

1. Find the number of S-representation of the loop $L_9(3)$.
2. How many S-pseudo representation exist for $L_{13}(4)$?
3. Find all S-representation of $L_{15}(8)$. How many are distinct?
4. Are all the S-representation in $L_{21}(11)$ distinct?
5. Find a loop L of order n, n odd so that it has only n, S-representations.
6. Can a loop of order n have lesser than n, S-representation?
7. Find a loop L of odd order which has only
    i. S-representation.
    ii. S-pseudo representation.

## 4.4 Smarandache isotopes of loops

In this section we introduce the concept of Smarandache isotopes in loops. The Smarandache isotopes exists only when $L_n(m) \in L_n$ are loops of order n + 1 where n is a non-prime. We see for a given loop we can have several S-isotope loops. Characterize those loops for which S-isotopes are only one.

**DEFINITION 4.4.1**: *Let (L, •) be a loop, we define Smarandache principal isotope of a S-subloop A, (A, \*) of the subloop (A, •) for any predetermined a, b $\in A$ for all x, y $\in A$, x \* y = X • Y where X • a = x and b • Y = y. If L has no S-subloops only subgroups then we do not get the concept of Smarandache principal isotope even if L is a S-loop having subgroups. If L is a S-loop having no S-subloops then we define for any subgroup Smarandache pseudo isotopes associated with L. In case the loop is a S-loop, having no S-subloop the concept of Smarandache pseudo isotope coincides with the principal isotope of the loop.*

**DEFINITION 4.4.2**: *Let L be a loop having a S-subloop. L is said to be a Smarandache G-loop (S-G-loop) if it is isomorphic to all of its Smarandache principal isotopes (S-principal isotopes) on the same S-subloop. Thus we see unlike in the case of principal isotopes there are many depending on each of the S-subloop.*

**THEOREM 4.4.1**: *No loop in the class $L_n$ is a S-G-loop.*

*Proof*: Left for the reader to verify using number theoretic techniques as the reader is expected to have a good background of number theory.

This theorem is also left as an exercise for the reader.

**THEOREM 4.4.2**: *Let $L_n(m) \in L_n$. Then any S-principal isotope with respect to the pair (a, a) ($a \in A \subset L_n(m)$), A a S-subloop is commutative if and only if A is commutative.*

**THEOREM 4.4.3**: *The principal isotope of a commutative S-subloop A can also be a strictly non-commutative loop.*

*Proof*: Left as an exercise for the reader using number theoretical methods.

**PROBLEMS:**

1. Find for the loop $L_{51}$, a S-subloop A, the S-principal isotope.
2. Find for the loop $L_{45}(8)$; how many S-principal isotopes exist?
3. Find for the loop $L_{33}(17)$ a S-subloop, whose S-principal isotope is strictly non-commutative.
4. Find for the loops in the class of loops $L_{33}$: Loops which have several S-principal isotopes and a loop which has only one S-principal isotope.
5. Show $L_{19}$ has no S-principal isotopes only S-pseudo isotopes.

## 4.5 Smarandache hyperloops

The study of hypergroups was introduced only by De Mario Francisco [18] for the group $Z_n$ of integers modulo n under addition as $x * y = x + y$, $x + y + q$ for $q \in Z_n$ and denoted by $(Z_n, q)$. He proved that the hypergroups of $Z_n$ partitions $Z_n \times Z_n$. The definition of hyperloops has meaning only for loop built using $Z_n$. Thus the new class of loops defined in Chapter II alone has the relevance for hyperloops. We proceed on to define Smarandache hyperloops, Smarandache hyperloops II, Smarandache A-hyperloops and Smarandache A-hyperloops II.

We just recall the definition of hyperloops.

**DEFINITION [64]**: *Let $(L_n(m), \bullet)$ be a loop such that $(m, n) = 1$ and $(m - 1, n) = 1$ using modulo integers. Let '$*$' be a binary operation defined on $L_n(m)$ by the following rule:*

*For all $x, y \in L_n(m)$, $x * y = (x \bullet y, (x \bullet y) \bullet q)$ where $q \in L_n(m)$. Then $(L_n(m), *, q)$ is a hyperloop. Clearly $(L_n(m), *, q)$ is a subset of $L_n(m) \times L_n(m)$.*

**Example 4.5.1**: Let $L_5(4) = \{e, 1, 2, 3, 4, 5\}$. The hyperloops of $(L_5(4), \bullet)$ is given by
$(L_5(4), *, 5) = \{(e, 5), (1, 2), (2, 4), (3, 1), (4, 3), (5, e)\}$
$(L_5(4), *, 4) = \{(e, 4), (1, 3), (2, 5), (3, 2), (4, e), (5, 1)\}$
$(L_5(4), *, 3) = \{(e, 3), (1, 4), (2, 1), (3, e), (4, 5), (5, 2)\}$

$(L_5(4), *, 2) = \{(e, 2), (1, 5), (2, e), (3, 4), (4, 1), (5, 3)\}$
$(L_5(4), *, 1) = \{(e, 1), (1, e), (2, 3), (3, 5), (4, 2), (5, 4)\}$
$(L_5(4), *, e) = \{(e, e), (1, 1), (2, 2), (3, 3), (4, 4), (5, 5)\}$

From this example we see that the hyperloops in general need not be a loop. For the simple reason for the hyperloops $\{L_5(4), *, n\}$ when $n \neq e$ is not even closed.

Now one of the natural questions would be, why should one take in the definition. (x * y) * q we could also take x * (y * q); when we take (x * y) * q instead of x * (y * q) the nature of partition is affected. So only in case of hyperloops we define hyperloop defined using x * (y * q) as A-hyperloops.

**DEFINITION [64]**: *Let $(L_n(m), \bullet)$ be a loop. Let $*$ be a binary operation defined on $L_n(m)$ as $x * y = (x \bullet y, x \bullet (y \bullet q))$ for all $x, y \in L_n(m)$ and $q \in L_n(m)$. $L_n(m)$ together with the binary operation $*$ is called the A-hyperloop and is given by $\{(L_n(m), *, q)_A = (x \bullet y, x \bullet (y \bullet q) / q \in L_n(m)\}$.*

The A-hyperloop for the $L_5(4)$ is given by the following example:

***Example 4.5.2***: Let $(L_5(4), \bullet)$ be the loop. The A-hyperloops of $(L_5(4), \bullet)$ is as follows: $(L_5(4), *, 5)_A = \{(e, 5), (1, 2), (2, 4), (3, 1), (4, 3), (5, 3), (5, e), (4, e),$ $(3, 4), (2, 1), (3, e), (1, 3), (5, 1), (4, 2), (5, 4), (2, e), (1, e), (5, 2)\}$. Similarly we get 18 elements in each of the cases $(L_5(4), *, 4)_A, \ldots, (L_5(4), *, 1)_A$. But in case of $(L_5(4), *, e)_A = \{(e, e), (1, 1), (2, 2), (3, 3), (4, 4), (5, 5)\}$.

Clearly A-hyperloops and hyperloops behave in a different way on the same loop which is evident from these examples.

For more about these loops please refer [64].

**DEFINITION 4.5.1**: *Let $L_n(m) \in L_n$ be a loop. The Smarandache hyperloop (S-hyperloop) of $L_n(m)$ is $(L_n(m), q)$ is a subset of $L_n(m) \times L_n(m)$ and q is an element of a S-subloop. A of $L_n(m)$ that is $x * y = (x \bullet y, (x \bullet y) \bullet q)$ note if $L_n(m)$ has no S-subloops but is a S-loop then we replace A by $L_n(m)$ in this case the Smarandache hyperloop coincides with the hyperloop.*

**DEFINITION 4.5.2**: *Let $L_n(m) \in L_n$ be a loop $(L_n(m), q)$ is defined as the Smarandache hyperloop II if $q \in A$, where A is a S-subloop II. $x * y = (x \bullet y, (x \bullet y) \bullet q)$.*

Note in both the definitions 4.5.1 and 4.5.2 if we take x • (y • q) instead of (x ∗ y) • q we get a different S-hyperloops. Like in definition 4.5.1 if L has no S-subloop II but L is a S-loop II then S-hyperloop is defined by replacing A by L.

**THEOREM 4.5.1**: *Let $L_n(m) \in L_n$. If n is a prime then S hyperloop of $L_n(m)$ is the same as hyperloop of $L_n(m)$.*

*Proof*: If n is a prime $L_n(m)$ is only a S-loop and $L_n(m)$ has no S-subloops.

**DEFINITION 4.5.3**: *Let $L_n(m) \in L_n$ we say the Smarandache A-hyperloop (S-A-hyperloop) of $L_n(m)$ is defined as $\{L_n(m), q\}$ where q is in a S-subloop B of L. If $L_n(m)$ has no S-subloop then we replace B by L provided L is a S-loop and we say in this case the A-hyperloop and the S-A-hyperloop coincide. $(L_n(m), q, *) = \{x \bullet y, x \bullet (y \bullet q) \mid q \in B\}$.*

**DEFINITION 4.5.4**: *Let $L_n(m) \in L_n$. The Smarandache A-hyperloop II (S-A-hyperloop II) is the set $(L_n(m), q, *) = \{x \bullet y, x \bullet (y \bullet q) \mid q \in B$ where B is a S-subloop II$\}$ If $L_n(m)$ has no S-subloop II but is a S-loop II then we replace B by $L_n(m)$.*

Now we will illustrate by 2 examples one for each.

***Example 4.5.3***: Let $L_{45}(8) \in L_{45}$, A = {e, 1, 6, 11, 16, 21, 26, 31, 36, 41} is a subloop of $L_{45}(8)$ which is a S-loop. The S-hyperloop with respect to A. Find $(L_{45}(8), *, q) = \{x \bullet y, x \bullet (y \bullet q) \mid q \in A\}$ It is only simple calculation to find the S-A-hyperloop of A. Find $(L_{45}(8), *, q) = \{x \bullet y, x \bullet (y \bullet q) \mid q \in A\}$. Compare them.

***Example 4.5.4***: Let $L_9(5)$ be the loop given by the following table:

|   | e | 1 | 2 | 3 | 4 | 5 | 6 | 7 | 8 | 9 |
|---|---|---|---|---|---|---|---|---|---|---|
| e | e | 1 | 2 | 3 | 4 | 5 | 6 | 7 | 8 | 9 |
| 1 | 1 | e | 6 | 2 | 7 | 3 | 8 | 4 | 9 | 5 |
| 2 | 2 | 6 | e | 7 | 3 | 8 | 4 | 9 | 5 | 1 |
| 3 | 3 | 2 | 7 | e | 8 | 4 | 9 | 5 | 1 | 6 |
| 4 | 4 | 7 | 3 | 8 | e | 9 | 5 | 1 | 6 | 2 |
| 5 | 5 | 3 | 8 | 4 | 9 | e | 1 | 6 | 2 | 7 |
| 6 | 6 | 8 | 4 | 9 | 5 | 1 | e | 2 | 7 | 3 |
| 7 | 7 | 4 | 9 | 5 | 1 | 6 | 2 | e | 3 | 8 |
| 8 | 8 | 9 | 5 | 1 | 6 | 2 | 7 | 3 | e | 4 |
| 9 | 9 | 5 | 1 | 6 | 2 | 7 | 3 | 8 | 4 | e |

This is a loop of order 10 commutative and is a S-loop. The S-subloop of $L_9(5)$ is A = {e, 2, 5, 8} as {e, 2} is a subgroup of $L_9(5)$. S-hyperloop $\{L_9(3), *, e\} = \{(e,e),(1,1),$ (2, 2), ... , (9, 9)}. $\{L_9(3), *, 2\}$, $\{L_9(3), *, 5\}$ and $\{L_9(3), *, 8\}$ can be calculated and test whether the S-hyperloops partition $L_9(5) \times L_9(5)$; on similar lines one can find S-A-hyperloop, S-hyperloop II and S-A hyperloop II and compare them and test whether they partition or not.

**PROBLEMS:**

1. Find a S-hyperloop of $L_{15}(8)$.
2. Find S-A-hyperloop of $L_{15}(8)$ and compare it with problem 1.
3. Can $L_{15}(8)$ have S-hyperloop II. Justify your answer.
4. Compare all the four types of hyperloops for the loop $L_{15}(8)$ (whenever it exists).
5. Find a general method of finding the S-hyperloop of $L_{49}(9)$.
6. Find a S-hyperloop of $L_{25}(7)$.
7. Can S-hyperloop II exist in $L_{25}(7)$? If it exist find it.
8. Find the S-hyperloop and S-A-hyperloop of $L_{27}(14)$. Compare them, which of them partition $L_{27}(14) \times L_{27}(14)$.
9. Prove $L_n(m)$ when $m = \dfrac{n+1}{2}$ has S-hyperloops and S-A-hyperloop (n is not a prime).
10. Find all S-hyperloops for $L_{21}(11)$. Compare it with S-A-hyperloops.

**Chapter five**

# RESEARCH PROBLEMS

This chapter proposes open research problems to researchers / students / algebraists and above all to those mathematicians who form the growing community of researching the Smarandache notions. Since already several researchers are working on Smarandache algebraic structures we by writing this book add to the attraction of more and more students / researchers to study as this Smarandache structures paves way for analysis of any structure is an exemplary way which cannot be done otherwise. We list here the fifty two research problems which will be a major attraction to all Smarandache mathematicians.

1. Let L be a loop. Can we prove the notion of u.p and t.u.p are equivalent on loops? (In 1980, Strojnowski, A. [62] has proved that in the case of groups, t.u.p and u.p coincide.)

2. Let L be a S-loop. Can we prove the notion of S.u.p and S.t.u.p are equivalent at least on S-loops?

3. Characterize those class of loops which are

    i) u.p loops,
    ii) S.u.p loops,
    iii) t.u.p loops,
    iv) S.t.u.p loops.

4. Characterize those commutative loops which are strongly semi-right commutative.

5. Characterize those Moufang loops which are semi-right commutative but non-commutative.

6. Let L be a non-commutative loop.

    i. What is the relation between P(L) and SP(L)?
    ii. When is P(L) = SP(L)? (that is characterize those loops L on which the equality holds good)

7. Let L be a non-associative loop.

    i. Find a relation between SPA(L) and PA(L).
    ii. Find conditions on the loop L so that PA(L) = SPA(L).

8. Can a Jordan loop be a Bruck loop? Moufang loop? Bol loop?

9. Find conditions on loops L so that both the first normalizer is equal or identical with the second normalizer for all subloops.

10. Characterize those loops L for which the first and the second normalizer are distinct for all subloops of L.

11. Does there exists a A-loop, which is not Moufang?

12. Does there exist a C-loop, which is not a S-loop? or is every C-loop a S-loop?

13. Characterize those loops L (where L is not a loop in the class $L_n$ with n a prime), which are S-subgroup loop.

14. Does there exist a loop L that has normal subloops but they are not S-normal subloops?

15. Can the class of loops $L_n$ for any odd n, n > 3 have a subgroup whose order is greater than four?

16. Can we ever have a class of loops, which contains the group $S_n$ as its subgroup? (Does not include the S-mixed direct product of loops).

17. Can we have an extension of Cayleys theorem for at least S-loops.(Hint: If the solution to the above problem is true certainly a solution to this problem exists or we have Cayleys theorem to be true for S-loops).

18. Does there exists a S-strongly commutative loop, which is not a S-strongly cyclic loop?

19. Characterize those loops L which have only one S-commutator subloop.

20. Characterize those loops L that has always the S-commutator subloop to be coincident with the commutator subloop. (Not loops from the class $L_n$).

21. Characterize those loops L, which has n-S-subloop and n distinct S-commutator subloops, associated with them. (Is this possible?)

22. Characterize those loops L that has many distinct S-subloops but one and only one S-commutator subloop.

23. Characterize those loops L that has always the S-associator subloop to be coincident with the associator subloop of L. (Not loops from the class $L_n$ as it has been already studied.)

24. Characterize those loops L, which has n-S-subloops and n distinct S-associator subloops, associated with them.

25. Characterize those loops L that has one and only one S-associator subloop (L has S-subloops).

26. Characterize those loops L in which $L^A = PA(L^S) = SPA(L^S) = A(L)$.

27. Can loop $L_n(m)$ where n is not a prime have proper subloops which are

    i. Moufang loops?
    ii. Bruck loops?
    iii. Bol loops?

    (This will in turn prove $L_n(m)$ when n is not prime is a S-Moufang loop, S-Bol loop and S-Bruck loop)

28. Can we say every Smarandache strong Moufang triple (Smarandache strong Bol triple or Smarandache strong Bruck triple) generate a Moufang subloop (Bol subloop or Bruck subloop) which is a S-subloop of the given loop.

29. Characterize those loops L, which are not Moufang but in which every S-subloop is a Moufang loop that is L is a S-strongly Moufang loop.

30. A similar problem in case of

    i. non-Bol loops which are S-strongly Bol loops.
    ii. non-Bruck loops which are S-strong Bruck loops.
    iii. non-WIP loops which are S-strong WIP loops.
    iv. non-diasscoiative loops which are S-strong diassociative loops.
    v. non-power associative loops which are S-strong power associative loops.

    Note: If a loop L has only one S-loop and no other S-subloop which has the stipulated property L becomes a S-strongly loop having that property.

31. Characterize those loops which have many S-subloops still

    i. S-Moufang centre is unique.
    ii. S-nucleus $SN_\lambda SN_\mu SN_\rho$ is unique.
    iii. S-centre is unique.

iv. S-first normalizer is unique.
v. S-second normalizer is unique.

32. Characterize those loops which has several S-subloops having distinct S-Moufang centre, S-centre, S-nucleus, S-first and S-second normalizer.

33. Characterize those loops for which every S-subloop has $SN_1 = SN_2$.

34. Characterize those loops for which every distinct pair of S-subloops $SN_1 \neq SN_2$.

35. Does there exist examples of S-loops of prime order in which every subgroup is a normal subgroup?

36. Characterize those loops L in which every subgroup is a normal subgroup.

37. Obtain loops L of odd order not got from the S-mixed direct product of loops but which satisfy S-Lagrange criteria.

    i. Does there exist such loops?
    ii. Can you characterize such loops of odd order.

38. Characterize those loops which satisfy S-Sylow criteria or equivalently those loops, which do not satisfy S-Sylow criteria.

39. Suppose a loop L satisfies S-Lagrange's criteria can we say L satisfies S-Sylow criteria? Characterize those loops, which satisfy both!

40. For any S-loop L when will every subgroup A give

    i. $\bigcup_{i=1}^{n} Ax_i = L \ (Ax_i \neq Ax_j, Ax_i \cap Ax_j = \phi)$

    ii. $\bigcup_{i=1}^{n} x_i A = L; \ (x_i \in L; x_i A \neq x_j A, x_i A \cap x_j A = \phi)$

    iii. Can the same elements $\{x_1, x_2, x_3, \ldots, x_n\} \subset L$ serves the propose?

    Characterize those S-loops for which i, ii and iii is true.

41. Prove for all non-commutative loops $L_n(m) \in L_n$, n a prime for every subgroup $A_i = \{e, i\}$, $i = 1, 2, \ldots, n$ we have the S-right coset decompositions $L_n(m) = \bigcup_{i=1}^{\frac{n+1}{2}} A_j x_i \ (Ax_i \cap Ax_k = \phi, i \neq k)$, $L_n(m) = \bigcup_{i=1}^{\frac{n+1}{2}} A_j y_i \ (Ay_i \cap Ay_k = \phi, i \neq k)$ where the set $X = \{x_1, \ldots, x_{\frac{n+1}{2}}\}$ and $Y = \{y_1, y_2, \ldots, y_{\frac{n+1}{2}}\}$ are such that $X \cup Y = L_n(m)$ and $X \cap Y = \phi$.

42. Study problem 41 in case S-left coset representation for the same class of loops.

43. What happens to the coset representation for the loop $L_n\left(\frac{n+1}{2}\right)$ which is the only commutative loop in the class $L_n$? Can we have problem 42 and 43 to be true? If so illustrate with examples and give characterization theorem for such loops.

44. Let L be a S-loop of odd order n. Study the coset representation when

    i. n is a prime.
    ii. n is an odd prime.

45. Let $L_n(m) \in L_n$. Suppose B is a subgroup of order greater than 2. When will

    i. $L_n(m) = \bigcup Bx_i$, $x_i \in L_n(m)$ (with $Bx_i \cap Bx_j = \phi$, if $x_i \neq x_j$)
    ii. $L_n(m) = \bigcup x_iB$, $x_i \in L_n(m)$ (with $x_iB \cap x_jB = \phi$, if $x_i \neq x_j$) possible? Characterize such loops in $L_n$.

46. Find those loops in $L_n(m) \in L_n$ where n is a power of a prime which have

    i. S-hyperloops which partition $L_n(m)$.
    ii. S-A hyperloop which partition $L_n(m)$.
    iii. S-hyperloop II of $L_n(m)$ for a suitable m.
    iv. S-A hyperloop II of $L_n(m)$ for a suitable m.

    Derive some interesting results and relations about these S-hyperloops.

47. Do we have loops in $L_n(m) \in L_n$, $n = n = p_1^{\alpha_1} \ldots p_k^{\alpha_k}$ $k > 2$ which are such that their S-hyperloops partition them?

48. Find some nice characterization about S-hyperloops and S-A-hyperloops. (As all loops in $L_n(m) \in L_n$ which are non-commutative are simple we cannot have S-hyperloops II).

49. Characterize those loops for which all S principal isotopes II are isomorphic with each other.

50. Characterize those loops L for which S-representations II exists. Can we find a method to prove or disprove the number of S-representation II is equal to order of L, or equal to order of $A \subset L$ where A is a S-loop II?

51. Find all S-loop II which satisfy
    i. S-Sylow criteria II
    ii. S-Lagrange criteria II.

52. Let $L_n(m)$ be a loop where $n = p_1^{\alpha_1} p_2^{\alpha_2} \ldots p_k^{\alpha_k}$, $\alpha_i > 1$ and $p_1, p_2, \ldots p_k$ are distinct primes. Find the lattice diagram of

    1. S-subloops
    2. subgroups.

# INDEX

### A

### B

### C

## M

## N

## O

## P

## Q

## ABOUT THE AUTHOR

Dr. W. B. Vasantha is an Associate Professor in the Department of Mathematics, Indian Institute of Technology Madras, Chennai, where she lives with her husband Dr. K. Kandasamy and daughters Meena and Kama. Her current interests include Smarandache algebraic structures, fuzzy theory, coding/ communication theory. In the past decade she has completed guidance of seven Ph.D. scholars in the different fields of non-associative algebras, algebraic coding theory, transportation theory, fuzzy groups, and applications of fuzzy theory to the problems faced in chemical industries and cement industries. Currently, six Ph.D. scholars are working under her guidance. She has to her credit 241 research papers of which 200 are individually authored. Apart from this she and her students have presented around 262 papers in national and international conferences. She teaches both undergraduate and post-graduate students at IIT and has guided 41 M.Sc. and M.Tech projects. She has worked in collaboration projects with the Indian Space Research Organization and with the Tamil Nadu State AIDS Control Society. She is currently authoring a ten book series on Smarandache Algebraic Structures in collaboration with the American Research Press.

She can be contacted at vasantha@iitm.ac.in
You can visit her on the web at: http://mat.iitm.ac.in/~wbv